\documentclass[USenglish]{article}
\usepackage[utf8]{inputenc}
\usepackage[small]{dgruyter}
\usepackage{makeidx}
\usepackage{amssymb,amsmath}
\usepackage{graphicx}
\usepackage{booktabs}
\def\bPsi{\boldsymbol\Psi}
\def\bPhi{\boldsymbol\Phi}
\renewcommand{\leq}{\leqslant}
\renewcommand{\geq}{\geqslant}
\DeclareMathOperator{\Rea}{Re}

\bibpunct{[}{]}{,}{n}{}{;} 
\def\newblock{\hskip .11em plus .33em minus .07em}
\newtheorem{proposition}{Proposition}
\begin{document}
  \author[1]{Alexander Zlotnik}
  \author[2]{Ilya Zlotnik}
  \runningauthor{A. Zlotnik and I. Zlotnik}
  \affil[1]{Department of Higher Mathematics at Faculty of Economics,National Research University Higher School of Economics,Myasnitskaya 20, 101000 Moscow, Russia}
  \affil[2]{Settlement Depository Company, 2-oi Verkhnii Mikhailovskii proezd 9, building 2, 115419 Moscow, Russia}
  \title{On the Richardson Extrapolation in Time of Finite Element Method with Discrete TBCs for the Cauchy Problem for the 1D Schr\"odinger Equation}
\runningtitle{On the Richardson Extrapolation of FEM with DTBCs}

\abstract{We consider the Cauchy problem for the 1D generalized Schr\"odinger equation on the whole axis.
To solve it, any order finite element in space and the Crank-Nicolson in time method with the discrete transpa\-rent boundary conditions (TBCs) has recently been constructed.
Now we engage the Richardson extrapolation to improve significantly the accuracy in time step.
To study its properties, we give results of numerical experiments and enlarged practical error analysis for three typical examples.
The resulting method is able to provide high precision results in the uniform norm for reasonable computational costs that is unreachable by more common 2nd order methods in either space or time step.
Comparing our results to the previous ones, we obtain \textit{much more} accurate results using \textit{much less} amount of both elements and time steps.
}

  \journalyear{2014}
  \startpage{1}
  \aop
\maketitle

\section{Introduction}
The time-dependent Schr\"odinger equation is the key one in many fields from quantum mechanics to wave physics. It should be often solved in unbounded space domains.
A number of approaches were developed to deal with such prob\-lems using approximate transparent boundary conditions (TBCs) at the artificial boundaries, see review \cite{AABES08}.
\par Among the best methods are those using the so-called discrete TBCs, see \cite{A98,EA01,M04} and \cite{DZ06}-\cite{DZZ09}, \cite{IZ13}, remarkable by clear mathematical background and the corres\-pon\-ding rigorous stability theory as well as complete absence of spurious reflections in practice.
Higher order methods of such kind are of special interest due to their practical efficiency.
To solve the 1D generalized Schr\"odinger equation on the axis or half-axis, any order finite element in space and the Crank-Nicolson in time method with the discrete TBCs has recently been constructed, studied and verified \cite{ZZ12,IZ13}.
\par In this paper, we present results on stability of the method in two norms and engage the Richardson extrapolations of increasing orders to improve significantly the accuracy in time step.
To demonst\-rate its nice practical error properties in various respects, we present enlarged results of the error analysis in numerical expe\-ri\-ments on the propaga\-tion of the Gaussian wave package for three rather standard examples: the free propagation in Example 1, tunneling through a rectangular barrier in Example 2 and a double barrier stepped quantum well in Example 3 (the last is the most complicated one in \cite{AABES08}).
The method is truly able to provide high precision results in the uniform norm (required in some problems in quantum mechanics) for reasonable computational costs that is unreachable by the 2nd order methods in either space or time step and not demonstrated previously.
\par Comparing our results to the previous ones, we obtain \textit{much more} accurate results using \textit{much less} amount of both space elements $J$ and time steps $M$.
In particular, concerning Example 3, we achieve the relative uniform in time and $L^2$ in space error $e=4E$$-6$
using only $J=36$ (!) and $M=2016$ for the 9th degree finite elements and the Richardson 6th order extrapolation method versus the best $e=4E$$-4$, $J=6000$ and $M=16000$ presented in \cite{AABES08}.

\section{The Cauchy problem and numerical methods}

We deal with the Cauchy problem for the 1D time-dependent generalized Schr\"o\-din\-ger equation on the whole axis
\begin{gather}
 i\hbar\rho D_t\psi={\mathcal H}\psi:=-\textstyle{\frac{\hbar^{\,2}}{2}}D(B D\psi)+V\psi\ \
 \text{on}\ \ \mathbb{R}\times\mathbb{R}^+,
\label{eq:se}\\
 \left. \psi\right|_{t=0}=\psi^0(x)\ \ \mbox{on}\ \ \mathbb{R}.
\label{eq:ic}
\end{gather}
Hereafter $\psi=\psi(x,t)$ is the complex-valued unknown wave function,
$i$ is the imaginary unit and $\hbar>0$ is a physical constant.
The $x$-depending coefficients $\rho,B,V\in L^\infty(\mathbb{R})$
are real-valued and satisfy $\rho(x)\geq \underline{\rho}>0$ and $B(x)\geq \underline{B}>0$.
Additionally $D_t=\frac{\partial}{\partial t}$ and $D=\frac{\partial}{\partial x}$ are the partial derivatives.
\par We also assume that, for some (sufficiently large) $X_0>0$,
 \begin{equation}
 \rho(x)=\rho_{\infty},\ \ B(x)=B_{\infty}>0,\ \ V(x)=V_{\infty}
 \ \ \mbox{and}\ \ \psi^0(x)=0\ \ \mbox{for}\ \ |x|\geq X_0.
 \label{eq:s21}
 \end{equation}
More generally, it could be assumed that $\rho$, $B$ and $V$ have different constant values for $x\leq -X_0$ and for $x\geq X_0$.
Let $\Omega=(-X,X)$ for some $X>X_0$.
\par We consider the weak solution $\psi\in C([0,\infty); H^1(\mathbb{R}))$ having
$D_t\psi\in C([0,\infty); L^2(\mathbb{R}))$ and satisfying the integral identity
\begin{equation}
 i\hbar(\rho D_t\psi(\cdot,t),\varphi)_{L^2(\mathbb{R})}
 =\mathcal{L}_{\mathbb{R}}(\psi(\cdot,t),\varphi)\ \ \text{for any}\
 \varphi\in H^1(\mathbb{R}),
\label{weak1}
\end{equation}
for any $t>0$. Hereafter we use the standard complex Lebesgue and Sobolev spaces and a Hermitian-symmetric sesquilinear form related to $\mathcal{H}$:
\[
 \mathcal{L}_I(w,\varphi)
 :=\textstyle{\frac{\hbar^{\,2}}{2}}(BDw,D\varphi)_{L^2(I)}+(Vw,\varphi)_{L^2(I)},\ \
 \text{with}\ \ I=\mathbb{R}\ \ \text{or}\ \ \Omega.
\]
\par Let $\ldots<x_{-J}=-X<x_{-J+1}<\ldots<x_J=X<\ldots$ be a mesh on $\mathbb{R}$ and $\Delta_j:=[x_{j-1}, x_j]$ be elements, for any integer $j$.
We set $h_j:=x_j-x_{j-1}$ and assume that $x_{-J+1}\leq -X_0$, $x_{J-1}\geq X_0$ and $h_j=h$ for $j\leq -J+1$ or $j\geq J$.
Let $h_{\rm max}=\max_jh_j$.
\par For $n\geq 1$, let $H_{h,\infty}^{(n)}$ be the finite element space of (piecewise polynomial) functions
$\varphi\in H^1(\mathbb{R})$
such that $\varphi|_{\Delta_j}$ are complex polynomials of the degree no more than $n$,
for any integer $j$.
Let $H_h^{(n)}$ be the restriction of $H_{h,\infty}^{(n)}$ to $\bar{\Omega}$.
\par Let $\overline{\omega}^{\,\tau}_M$ be the uniform mesh in $[0,T]$, for some $T>0$, with nodes $t_m=m\tau$, $0\leq m\leq M$, and $\tau=\frac{T}{M}$.
Let $\overline{\partial}_tY^m:= \frac{Y^m-Y^{m-1}}{\tau}$ and $\overline{s}_tY^m:= \frac{Y^{m-1}+Y^m}{2}$.
\par We introduce the FEM-Crank-Nicolson approximate solution $\Psi$: $\overline{\omega}^{\,\tau}_M\to H_{h,\infty}^{(n)}$ satisfying the integral identity
\begin{equation}
 i\hbar(\rho\overline{\partial}_t\Psi^m,\varphi)_{L^2(\mathbb{R})}
 =\mathcal{L}_{\mathbb{R}}(\overline{s}_t\Psi^m,\varphi)
 \ \ \text{for any}\
 \varphi\in H_{h,\infty}^{(n)} \ \text{and}\ 1\leq m\leq M,
\label{p43}
\end{equation}
compare with \eqref{weak1},
and the initial condition
$\Psi|_{t=0}=\Psi^0\in H_{h,\infty}^{(n)}$, where $\Psi^0$ approxi\-ma\-tes $\psi^0$.
This method is well defined and stable as it follows from \cite{ZZ12}. But it cannot be practically imple\-men\-ted since the number of unknowns is infinite at each time level. Nevertheless it is possible to restrict its solution from $\mathbb{R}$ to $\bar{\Omega}$ by imposing the discrete TBCs at $x=\pm X$ provided that $\Psi^0(x_j)=0$ for $|j|\geq J-1$.
\par This restriction $\Psi=\Psi^{(\tau)}$: $\overline{\omega}^{\,\tau}_M\to H_h^{(n)}$ obeys the integral identity \cite{ZZ12}
\begin{gather}
 i\hbar(\rho \overline{\partial}_t\Psi^m,\varphi)_{L^2(\Omega)}
 =\mathcal{L}_{\Omega}(\overline{s}_t\Psi^m,\varphi)
\nonumber\\
 -\textstyle{\frac{\hbar^{\,2}}{2}}B_{\infty}(\mathcal{S}_{\rm ref}^{(n)\,m}\bPsi_X^m)\varphi^*(X)
 +\textstyle{\frac{\hbar^{\,2}}{2}}B_{\infty}(\mathcal{S}_{\rm ref}^{(n)\,m}\bPsi_{-X}^m)\varphi^*(-X)
\label{eq:ii2}
\end{gather}
for any $\varphi\in H_h^{(n)}$ and $1\leq m\leq M$,
and the initial condition $\Psi|_{t=0}=\Psi^0|_{\bar{\Omega}}\in H_h^{(n)}$.
Here $\bPsi_{\pm X}^m:=\{\Psi^l|_{x=\pm X}\}_{l=0}^m$ and $\varphi^*$ is the complex conjugate of $\varphi$.
The key point is that the operator ${\mathcal S_{\rm ref}^{(n)\,m}}$ has the discrete convolution form
\[
 {\mathcal S_{\rm ref}^{(n)\,m}}\bPhi^m
 =c_n\,\sum_{l=0}^{m} K_{\rm ref}^{(n),\,l}\,\Phi^{m-l}\ \ \text{for}\ \ \bPhi^m:=\{\Phi^l\}_{l=0}^m.
\]
The analytical calculation of the kernel $K_{\rm ref}^{(n)}$ (defined in turn as an $n$-multiple discrete convolution) and the constant $c_n$ is far from being simple and is presented in \cite{ZZ12}; we omit the explicit expressions here.
To compute the kernel, we apply the fast algorithm for computing discrete convolutions based on FFT, for example see \cite{TAL97}.
\par Let $\ell^m(\varphi)$ be a conjugate linear functional on $H_h^{(n)}$ that we add to the right-hand side of \eqref{eq:ii2} to study stability in more detail.
\begin{proposition}
Let $\ell^m(\varphi)=(F^m,\varphi)_{L^2(\Omega)}$ with $F^m\in L^2(\Omega)$ for $1\leq m\leq M$.
Then the following first stability bound holds
\begin{equation}
 \max_{0\leq m\leq M} \left\|\sqrt{\rho}\,\Psi^m \right\|_{L^2(\Omega)}
 \leq \left\|\sqrt{\rho}\,\Psi^0 \right\|_{L^2(\Omega)}
 +\frac{2}{\hbar}\sum_{m=1}^M \left\|\frac{F^m}{\sqrt{\rho}}\right\|_{L^2(\Omega)}\tau.
 \label{p12}
\end{equation}
\label{prop:1}
\end{proposition}
\par We introduce the ``energy'' norm such that
\begin{gather}
\hspace{-15pt}\|w\|_{{\mathcal H}+\hat{v}\rho;\,\Omega}^2
:={\mathcal L}_{\Omega}(w,w)+\hat{v}\left\|\sqrt{\rho}\,w\right\|^2_{L^2(\Omega)}>0
\,\ \text{for any}\ w\in H^1(\Omega),\ w\not\equiv 0,
\label{p20}
\end{gather}
for some real number $\hat{v}$.
In particular, for $\hat{v}$ so large that $V+\hat{v}\rho> 0$, \eqref{p20}
is clearly valid.
We define also the corresponding dual mesh depending norm
\[
 \|w\|^{(-1)}_h
 := \max_{\varphi\in H_h^{(n)}:\ \|\varphi\|_{{\mathcal H}+\hat{v}\rho;\,\Omega}=1}
 |\langle w,\varphi\rangle_{\Omega}|
 \leq c\|w\|_{H^{-1}(\Omega)},\ \ H^{-1}(\Omega)=[H^1(\Omega)]^*,
\]
where $\langle w,\varphi\rangle_{\Omega}$ is the conjugate duality relation on
$H^{-1}(\Omega)\times H^1(\Omega)$.
\begin{proposition}
Let $\ell^m(\varphi)=\langle F^m,\varphi\rangle_{\Omega}$
with $F^m\in  H^{-1}(\Omega)$ for $1\leq m\leq M$ and $F^0\in H^{-1}(\Omega)$ be arbitrary.
Then the following second stability bound holds
\begin{gather}
\max_{0\leq m\leq M}\left\|\Psi^m\right\|_{\mathcal{H}+\hat{v}\rho;\,\Omega}
 \leq \left\| \Psi^0 \right\|_{{\mathcal H}+\hat{v}\rho;\,\Omega}
\nonumber\\
+4\sum_{m=1}^M \Bigl(\frac{|\hat{v}|}{\hbar}\,\|F^m\|^{(-1)}_h
 + \left\| \overline{\partial}_tF^m\right\|^{(-1)}_h\Bigr)\tau
 +4\left\| F^0\right\|^{(-1)}_h.
\label{p22}
\end{gather}
\label{prop:2}
\end{proposition}
\par Propositions \ref{prop:1} and \ref{prop:2} are proved similarly to \cite{ZZ12}.
\par For sufficiently smooth $\psi$, the error of the described method $O(\tau^{2}+h_{\rm max}^{n+1})$ is of the order $n+1$ (i.e., any) in $h_{\rm max}$ but only the 2nd in $\tau$. To remove this drawback, we further engage the classical Richardson extrapolation in time \cite{G08}.
To this end, we assume that the following error expansion holds
\begin{gather}
 \psi^m-\Psi^{(\tau),\,m}=\sum_{k=1}^{r-1}g_k^m\tau^{2k}+O(\tau^{2r}+h_{\rm max}^{\tilde{n}+1}),\ \ 0\leq m\leq M,
\label{eq:errexp}
\end{gather}
for $r=2,3$ or 4,
with some functions $g_k$ independent of the space-time mesh, and $0\leq\tilde{n}\leq n$ depending on the space smoothness of $\psi$. Then, for $r=2,3$ and $4$ and $0\leq rm\leq M$, we can exploit the following Richardson extrapolations
\begin{align}
& \Psi_{2R}^{2m}=\frac43\Psi^{(\tau),\,2m}-\frac13\Psi^{(2\tau),\,m},
\label{eq:rich2}\\[1mm]
& \Psi_{3R}^{3m}=\frac{81}{40}\Psi^{(\tau),\,3m}-\frac{16}{15}\Psi^{(3\tau/2),\,2m}
 +\frac{1}{24}\Psi^{(3\tau),\,m},
\label{eq:rich3}\\[1mm]
& \Psi_{4R}^{4m}=\frac{1024}{315}\Psi^{(\tau),\,4m}-\frac{729}{280}\Psi^{(4\tau/3),\,3m}+\frac{16}{45}\Psi^{(2\tau),\,2m}
 -\frac{1}{360}\Psi^{(4\tau),\,m}.
\label{eq:rich4}
\end{align}
It is supposed that $M$ is multiple of 2 (i.e., even), 6 and 12 respectively in \eqref{eq:rich2}, \eqref{eq:rich3} and \eqref{eq:rich4}.
The coefficients in these formulas are specific numbers
being uniquely found so that expansion \eqref{eq:errexp} implies
the higher order error bound
\begin{equation}
 \psi^{rm}-\Psi_{rR}^{rm}=O(\tau^{2r}+h_{\rm max}^{\tilde{n}+1}),\ \ r=2,3,4.
\label{eq:errrich}
\end{equation}
The more higher order Richardson extrapolations could be introduced as well.
\par It is not difficult to derive the Cauchy problems for the functions $g_k$ in \eqref{eq:errexp}. They are similar to \eqref{eq:se}, \eqref{eq:ic} but with recurrently defined additional free terms and zero initial function. For example, we have
\[
 i\hbar\rho D_tg_1-{\mathcal H}g_1=-\frac{1}{24}i\hbar\rho D_t^3\psi+\frac{1}{8}D_t^2{\mathcal H}\psi\ \
 \text{on}\ \ \mathbb{R}\times\mathbb{R}^+,\ \ g_1|_{t=0}=0
\]
similarly to  \cite{G08}. Clearly the right-hand side of the equation can be rewritten shorter as $\frac{1}{12}i\hbar\rho D_t^3\psi$.
\par Notice that the computation of $\Psi^{(\tau)}$ needs asymptotically $aJM+bM^2$ arith\-me\-tic operations ($bM^2$ is due to the discrete convolutions), for some $a>0$ and $b>0$. It is easy to check that then the computation of $\Psi_{rR}$ requires totally
\begin{equation}
 \frac32 aJM+\frac54 bM^2,\ \ 2aJM+\frac{14}{9} bM^2,\ \ \frac52 aJM+\frac{15}{8} bM^2
\label{eq:richcost}
\end{equation}
arithmetic operations respectively for $r=2,3$ and $4$. So the additional costs for implementing the Richardson extrapolation are less than $(r-1)\cdot50\%$. See also the corresponding practical results in Table \ref{tab:EX01r:ExecutionTime} below.
\par The Richardson extrapolations allow to achieve much better accuracy than the basic Crank-Nicolson discretization for the same mesh and with that:
(i)~ they inherit stability properties;
(ii) they deal with the same discrete TBC;
(iii) they exploit the same code for passing from the current time level to the next one (repeatedly for several time steps).

\section{Numerical experiments and error analysis}
\label{FREE}

In our numerical experiments, we intend to study in detail the practical error behavior for the Richardson extrapolations.
We choose $\hbar=1$, $\rho(x)\equiv 1$ and $B(x)\equiv 2$ (in Examples 1 and 2) or $B(x)\equiv 1$ (In Example 3) (the atomic units) and use the finite uniform space mesh $x_j=jh$, $|j|\leq J$, with the step $h=\frac{X}{J}$.
Let $\Psi^0\in H_h^{(n)}$ be simply the interpolant of $\psi^0$.
\smallskip
\par\textbf{3.1}. In Example 1, we rely upon the known exact solution (the scaled Gaussian wave package) for the Cauchy problem \eqref{eq:se}, \eqref{eq:ic}
\[
 \psi=\psi_G(x,t)\equiv
 \frac{1}{\sqrt[4]{2\pi\alpha}\sqrt[+]{1+i\,\frac{t}{\alpha}}}
 \exp\left\{ik(x-x^{(0)}-kt)-\frac{(x-x^{(0)}-2kt)^2}{4(\alpha+it)}\right\},
\]
where $x^{(0)}$, $k$ and $\alpha>0$ are the real parameters, in the case $V(x)\equiv 0$ (the free propagation of the wave). Thus the initial function takes the form
\begin{equation}
\label{FEM:s62}
 \psi^0(x)=\psi_G(x,0)=\frac{1}{\sqrt[4]{2\pi\alpha}}
 \exp\left\{ik(x-x^{(0)})-\frac{(x-x^{(0)})^2}{4\alpha}\right\}.
\end{equation}
It satisfies the property $\|\psi^0\|_{L_2(\mathbb{R})}=1$.
Though formally $\psi^0(x)\neq 0$ for any $x$, it decays rapidly as $|x-x^{(0)}|\to\infty$.
\par We choose the parameters  $x^{(0)}=0$, $k=100$, $\alpha=\frac{1}{120}$ and $X=0.8$ thus ensuring
$\left|\psi^0(x)\right|<1E$$-8$ for $|x|\geq X$.
Notice that this limits from below the least error that can be achieved (if required, it can be easily improved by small increasing of $X$).
Since $\max_{t\geq 0}|\psi_G(-X,t)|<1E$$-8$ as well, we can simply pose the zero Dirichlet boundary condition $\Psi|_{x=-X}=0$ instead of the discrete TBC at
$x=-X$ (as in \cite{DZZ09,ZZ12}).
Let also $T=0.006$.
Almost the same data were taken in several papers including \cite{EA01,DZZ09,ZZ12,IZ13}.
\par
\begin{figure}[htbp]
    \begin{minipage}[h]{0.3\linewidth}\center{
        \includegraphics[width=1\linewidth]{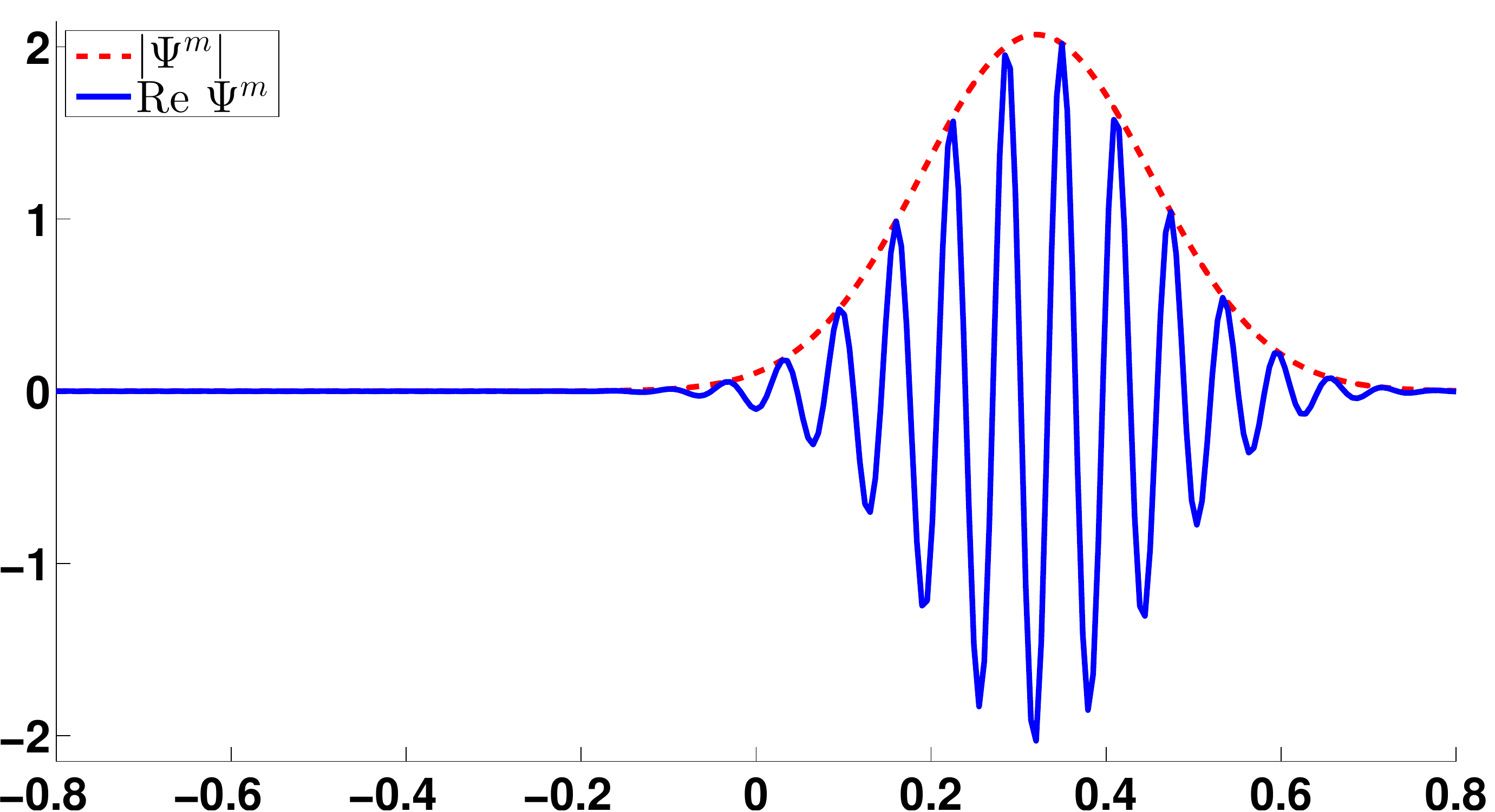}
    } \small{(a) $t_m=0.0016$, $m=80$} \\
    \end{minipage}\hfill
    \begin{minipage}[h]{0.3\linewidth}\center{
        \includegraphics[width=1\linewidth]{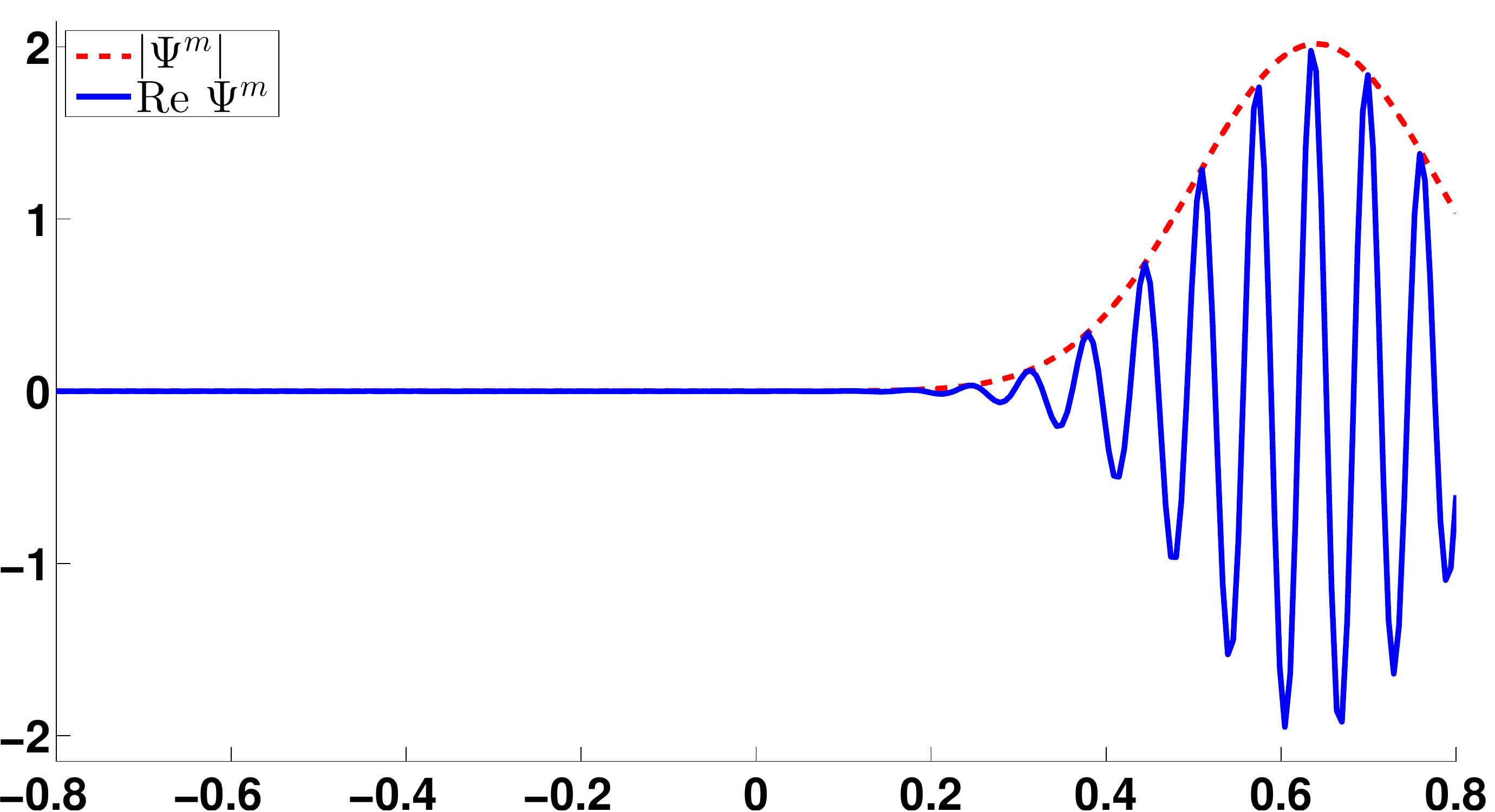}
    } \small{(b) $t_m=0.0032$, $m=160$} \\
    \end{minipage}\hfill
    \begin{minipage}[h]{0.3\linewidth}\center{
        \includegraphics[width=1\linewidth]{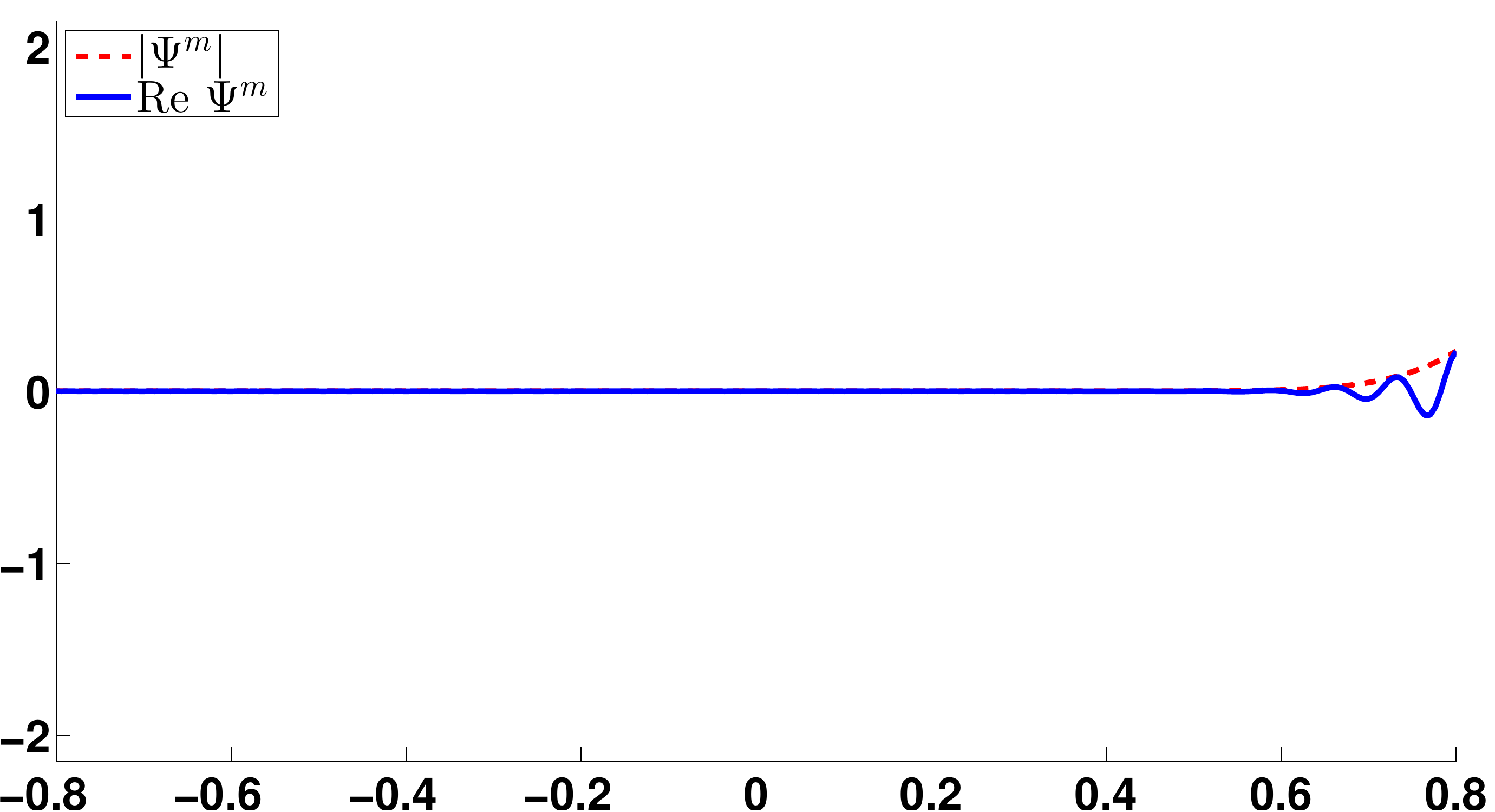}
    } \small{(c) $t_m=0.0056$, $m=280$} \\
    \end{minipage}\vfill
\caption{\small{Example 1. $|\Psi_{4R}^{m}|$ and $\Rea\Psi_{4R}^{m}$ for $n=9$ and $(J,M)=(30,300)$}
\label{fig:EX01r:Solution}}
\end{figure}

 On Fig.~\ref{fig:EX01r:Solution}, the solution is briefly represented by $\Psi_{4R}^m$ for high $n=9$ but $(J,M)=(30,300)$ only, with a suitable uniform accuracy (see Table \ref{tab:EX01r} below). The wave moves to the right, spreads slightly and leaves the domain $\bar{\Omega}$.
We emphasize that hereafter imposing of the discrete TBC does not produces any spurious reflections from the artificial boundary $x=X$ (as usual).
\par For any error $e^m$, we compute the mesh $L^2$-norm $\|e^m\|_{L_h^2}$
by applying the compound Newton-Cotes quadrature formula to the integral in $\|e^m\|_{L^2(\Omega)}$ (each element is divided into $n$ equal parts) and the mesh uniform norm $\|e^m\|_{C_h}$ (especially interesting in practice)
over the uniform mesh with the step $h/n$ in $\bar{\Omega}$.
Looking ahead, notice that though the theory concerns mainly $L^2$ or $H^1$-like norms, fortunately in general the practical error behavior in $C$ norm is close to $L^2$ one; this is not obvious at all in advance.
\par
\begin{figure}[htbp]
    \begin{minipage}[h]{0.49\linewidth}\center{
        \includegraphics[width=1\linewidth]{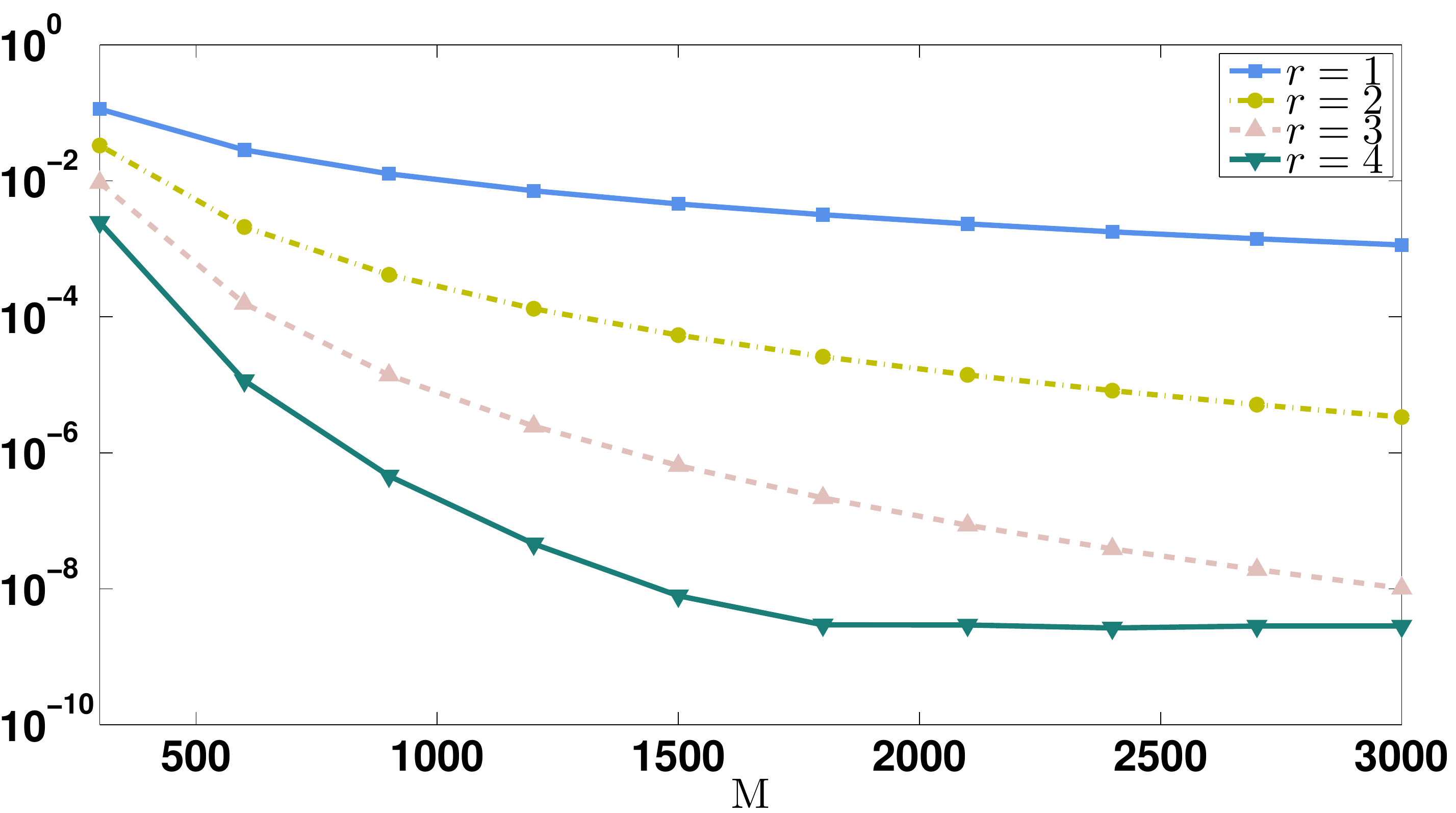}
    } \small{(a) in $L^2$ space norm} \\
    \end{minipage}\hfill
    \begin{minipage}[h]{0.49\linewidth}\center{
        \includegraphics[width=1\linewidth]{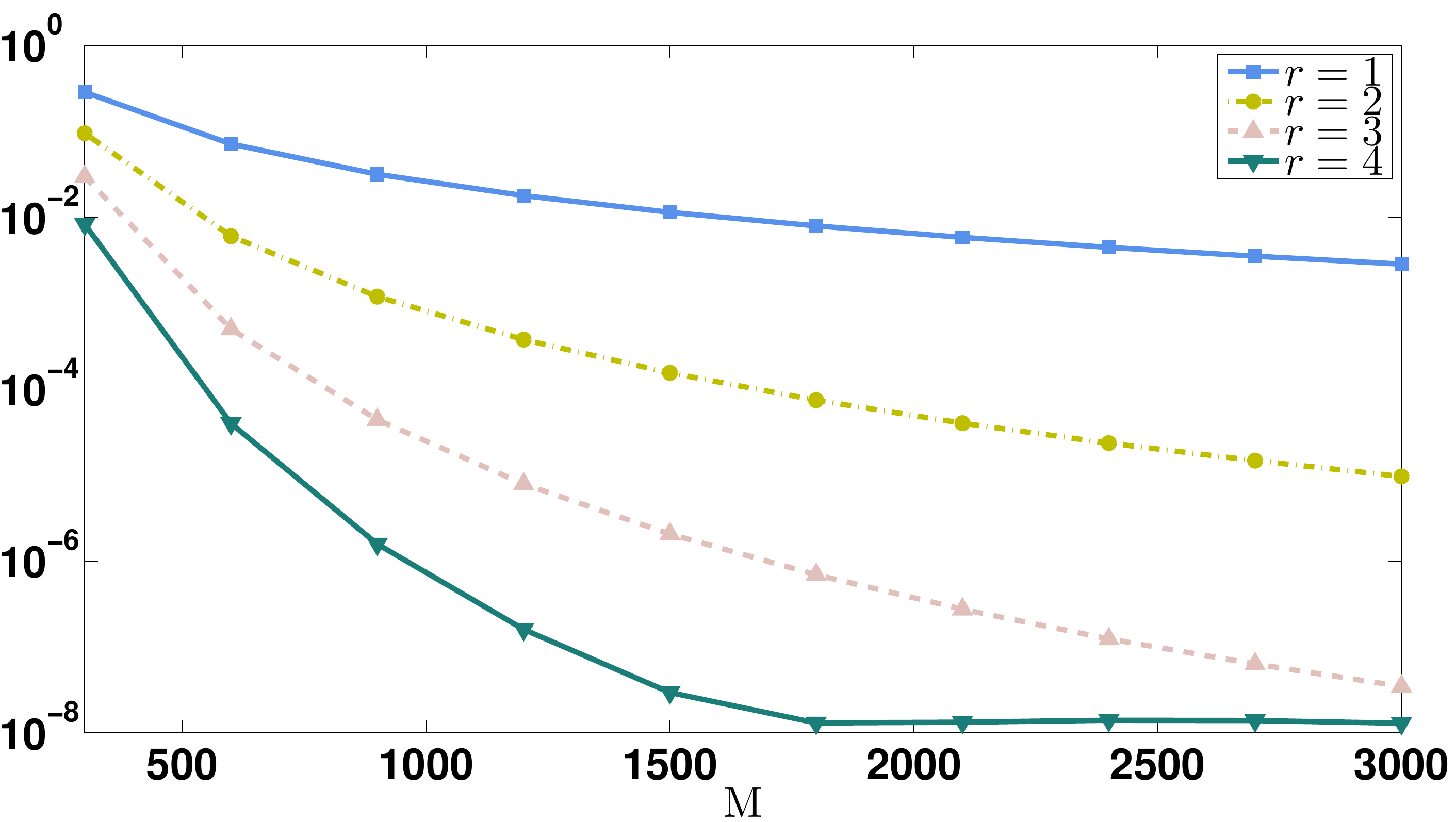}
    } \small{(b) in $C$ space norm} \\
    \end{minipage}
\caption{\small{Example 1. The errors $\max_{0\leq rm\leq M}\|\psi^{rm}-\Psi_{rR}^{rm}\|$,
for $n=9$ and $J=90$, in dependence with $r=1,2,3,4$ and $M=300,600,\dots,3000$}
\label{EX01r:FEM:N=9:J=90:MaxAbsError}}
\end{figure}

 On Fig. \ref{EX01r:FEM:N=9:J=90:MaxAbsError}, we present the errors
$\max_{0\leq rm\leq M}\|\psi^{rm}-\Psi_{rR}^{rm}\|$,
for $n=9$ and $J=90$, in dependence with $r=1,2,3,4$ and $M=300,600,\dots,3000$,
where we set $\Psi_{1R}=\Psi^{(\tau)}$ for convenience. For $r=1$, i.e. without the extrapolation, the errors decay too slowly. They decay faster and faster as $r$ grows excepting the case $r=4$ and $M\geq M_1=1800$, where the errors stabilize since their lowest levels have already been achieved.
Moreover, the error values decrease remarkably as $r$ grows: the ratio
\begin{equation}
 \max_{0\leq m\leq M}\|\psi^{m}-\Psi^{m}\|/\max_{0\leq rm\leq M}\|\psi^{rm}-\Psi_{rR}^{rm}\|,
\label{eq:errrat}
\end{equation}
for example, for $M=600$, equals (approximately) 13.6, 181 and 2473 in $L^2$ norm as well as 11.8, 143 and 1790 in $C$ norm whereas the ratio
\[
 \bigl.\max_{0\leq m\leq M}\|\psi^{m}-\Psi^{m}\|\bigr|_{M=3000}/
 \bigl.\max_{0\leq rm\leq M}\|\psi^{rm}-\Psi_{rR}^{rm}\|\bigr|_{M=600}
\]
equals 0.54, 7.26 and 98.9 in $L^2$ norm as well as 0.47, 5.70 and 71.6 in $C$ norm respectively for $r=2,3$ and $4$.
For $M=1800$, ratio \eqref {eq:errrat} equals already $122$, $14588$ and $1068678$ (!) in $L^2$ norm as well as $106$, $11451$ and $605095$ in $C$ norm respectively for $r=2,3$ and $4$. For the final $M=3000$, it is even much larger: 339 and 111688 in $L^2$ norm as well as 295 and 82000 in $C$ norm respectively for $r=2$ and $3$.
\begin{table}\centering{
\begin{tabular}{rcccccccc}
$M$&$J=30$&$J=40$&$J=50$&$J=60$&$J=70$&$J=80$&$J=90$\\ \midrule %
$300$&$2.45E{-}3$&$2.45E{-}3$&$2.45E{-}3$&$2.45E{-}3$&$2.45E{-}3$&$2.45E{-}3$&$2.45E{-}3$\\%
$600$&$1.53E{-}4$&$1.15E{-}5$&$1.15E{-}5$&$1.15E{-}5$&$1.15E{-}5$&$1.15E{-}5$&$1.15E{-}5$\\%
$900$&$1.41E{-}4$&$7.68E{-}6$&$6.79E{-}7$&$4.62E{-}7$&$4.57E{-}7$&$4.56E{-}7$&$4.56E{-}7$\\%
$1\,200$&$1.33E{-}4$&$6.73E{-}6$&$6.89E{-}7$&$1.11E{-}7$&$4.83E{-}8$&$4.63E{-}8$&$4.59E{-}8$\\%
$1\,500$&$1.19E{-}4$&$6.35E{-}6$&$6.70E{-}7$&$1.12E{-}7$&$2.85E{-}8$&$8.46E{-}9$&$7.91E{-}9$\\%
$1\,800$&$1.93E{-}5$&$6.54E{-}6$&$6.41E{-}7$&$1.12E{-}7$&$2.66E{-}8$&$8.34E{-}9$&$2.96E{-}9$\\
\end{tabular}%
\\[2mm]
\small{(a) in $L^2$ space norm}\\[2mm]
\begin {tabular}{rcccccccc}
$M$&$J=30$&$J=40$&$J=50$&$J=60$&$J=70$&$J=80$&$J=90$\\ \midrule %
$300$&$8.30E{-}3$&$8.27E{-}3$&$8.27E{-}3$&$8.27E{-}3$&$8.27E{-}3$&$8.27E{-}3$&$8.27E{-}3$\\%
$600$&$8.26E{-}4$&$4.27E{-}5$&$3.97E{-}5$&$3.97E{-}5$&$3.97E{-}5$&$3.97E{-}5$&$3.97E{-}5$\\%
$900$&$6.98E{-}4$&$2.70E{-}5$&$2.03E{-}6$&$1.63E{-}6$&$1.60E{-}6$&$1.58E{-}6$&$1.57E{-}6$\\%
$1\,200$&$7.26E{-}4$&$2.67E{-}5$&$1.82E{-}6$&$3.98E{-}7$&$1.81E{-}7$&$1.64E{-}7$&$1.60E{-}7$\\%
$1\,500$&$7.30E{-}4$&$1.97E{-}5$&$1.71E{-}6$&$3.42E{-}7$&$1.17E{-}7$&$3.69E{-}8$&$2.95E{-}8$\\%
$1\,800$&$6.41E{-}4$&$2.42E{-}5$&$1.68E{-}6$&$3.28E{-}7$&$1.09E{-}7$&$3.75E{-}8$&$1.31E{-}8$\\
\end {tabular}%
\\[2mm]
\small{(b) in $C$ space norm}}
\caption{\small Example 1. The errors $\max_{0\leq 4m\leq M}\|\psi^{4m}-\Psi_{4R}^{4m}\|$, for $n=9$, in dependence with
$J$ and $M$}
\label{tab:EX01r}
\end{table}
\par Table \ref{tab:EX01r} contains the errors $E_{J,M}:=\max_{0\leq 4m\leq M}\|\psi^{4m}-\Psi_{4R}^{4m}\|$, for $n=9$, in dependence with $J=20,30,\dots,90$ and $M=300,600,\dots,1800$, and is rich in information.
Clearly the values decrease as $J$ or $M$ increases though they (almost) stabilize as $J$ increases and $M$ is fixed or, vice versa, $J$ is fixed and $M$ increases.
Next, for example, for $(J,M)=(30,1200)$, the ratio $E_{J,M}/E_{2J,M}$ equals (approximately)
1198 in $L^2$ norm and 1824 in $C$ norm
that corresponds to $2^{n+1}=1024$
whereas, for $(J,M)=(90,600)$, the ratio $E_{J,M}/E_{J,2M}$ equals
250 in $L^2$ norm and 248 in $C$ norm
that agrees well to $2^{2r}=256$, see \eqref{eq:errrich}.
Also, for $(J,M)=(40,600)$, the ratio $E_{J,M}/E_{2J,2M}$ equals
2600 in $L^2$ norm and 2483 in $C$ norm that shows the rapid decay of the error.
\par
\begin{figure}[htbp]
    \begin{minipage}[h]{0.49\linewidth}\center{
        \includegraphics[width=1\linewidth]{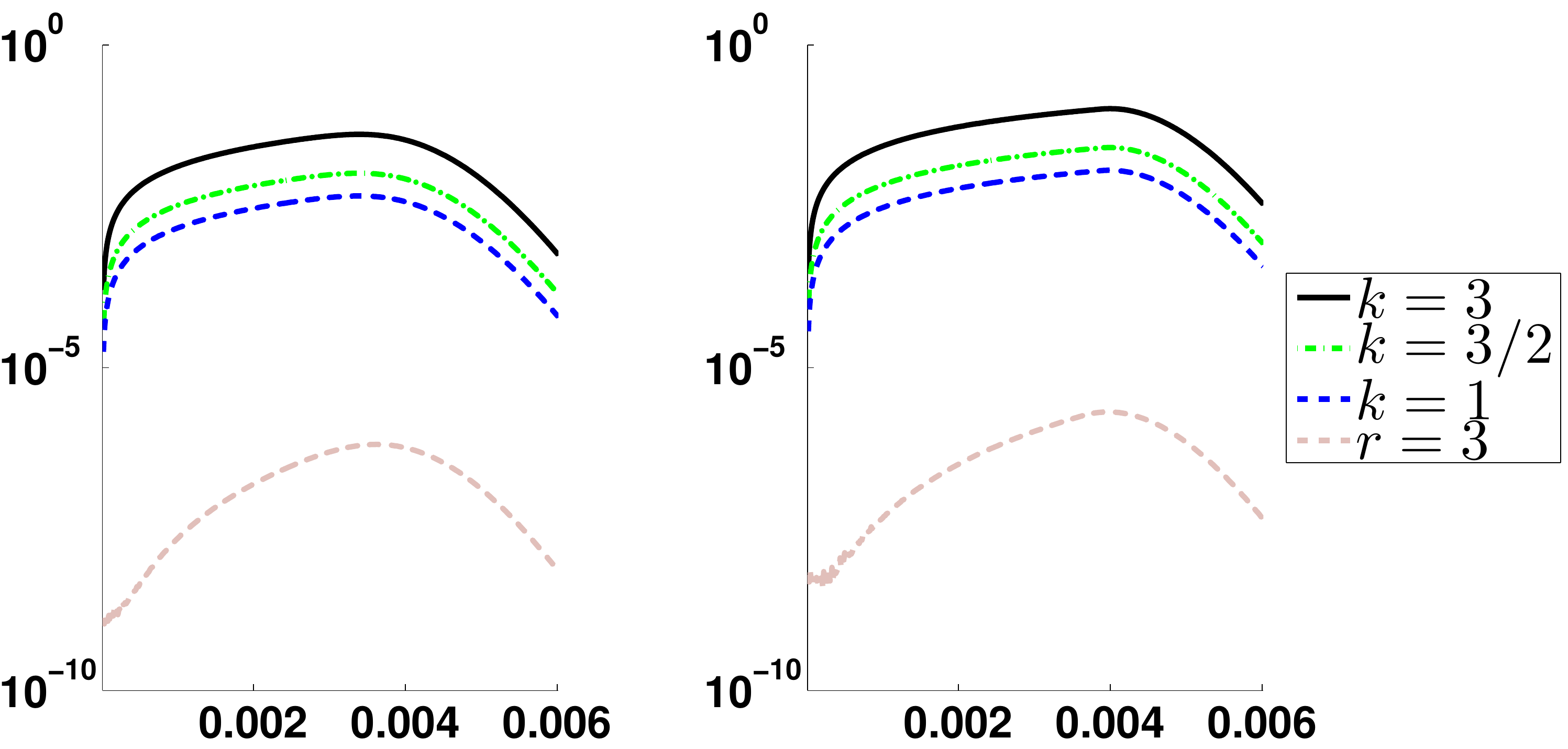}
    } \small{(a) in $L^2$ (left) and $C$ (right) norms, $r=3$}\\
    \end{minipage}\hfill
    \begin{minipage}[h]{0.49\linewidth}\center{
        \includegraphics[width=1\linewidth]{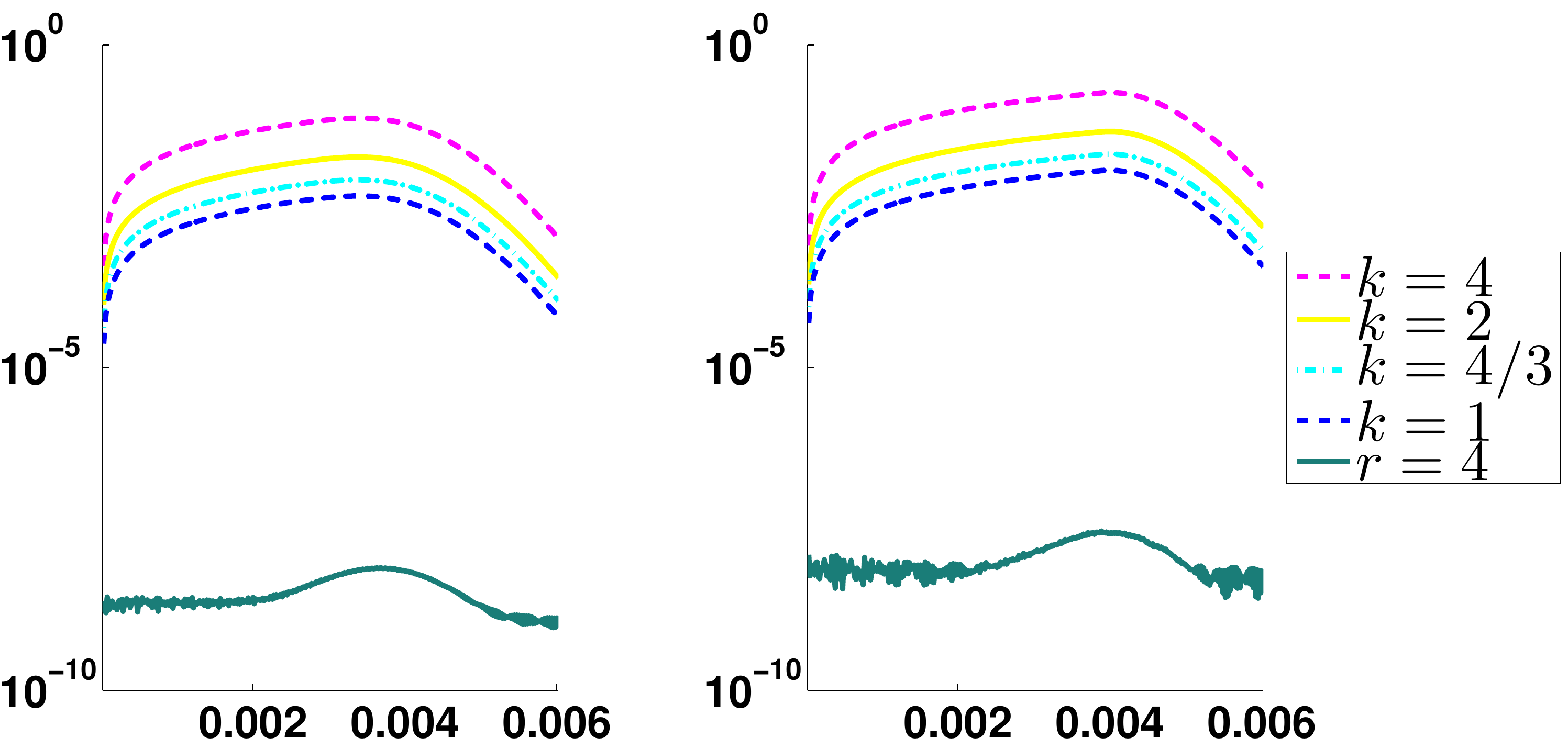}
    } \small{(b) in $L^2$ (left) and $C$ (right) norms, $r=4$}\\
    \end{minipage}\vfill
\caption{\small{Example 1.
The errors: (a) $\|\psi^{3m}-\Psi^{(k\tau),\,3m}\|$ for $k=1,3/2,3$
and $\|\psi^{3m}-\Psi_{3R}^{3m}\|$ ($r=3$), and
(b) $\|\psi^{4m}-\Psi^{(k\tau),\,4m}\|$ for $k=1,4/3,2,4$
and $\|\psi^{4m}-\Psi_{4R}^{4m}\|$ ($r=4$),
both for $n=9$ and $(J,M)=(90,1500)$, in dependence with $t_m$}
\label{EX01r:FEM:J=90:fig:Error}}
\end{figure}

 Fig. \ref{EX01r:FEM:J=90:fig:Error} demonstrates that the error $\|\psi^{rm}-\Psi^{(k\tau),\,rm}\|$ (corresponding to the summands of the Richardson extrapolation in \eqref{eq:rich3}, \eqref{eq:rich4}) decays monotonically but very slowly as $k$ decreases whereas the error of the Richardson extrapolation $\|\psi^{rm}-\Psi_{rR}^{rm}\|$, for $r=3$ and $4$, diminishes abruptly (by several orders of magnitude), for any $0\leq rm\leq M$.

\par
\begin{figure}[htbp]
    \begin{minipage}[h]{0.49\linewidth}\center{
        \includegraphics[width=1\linewidth]{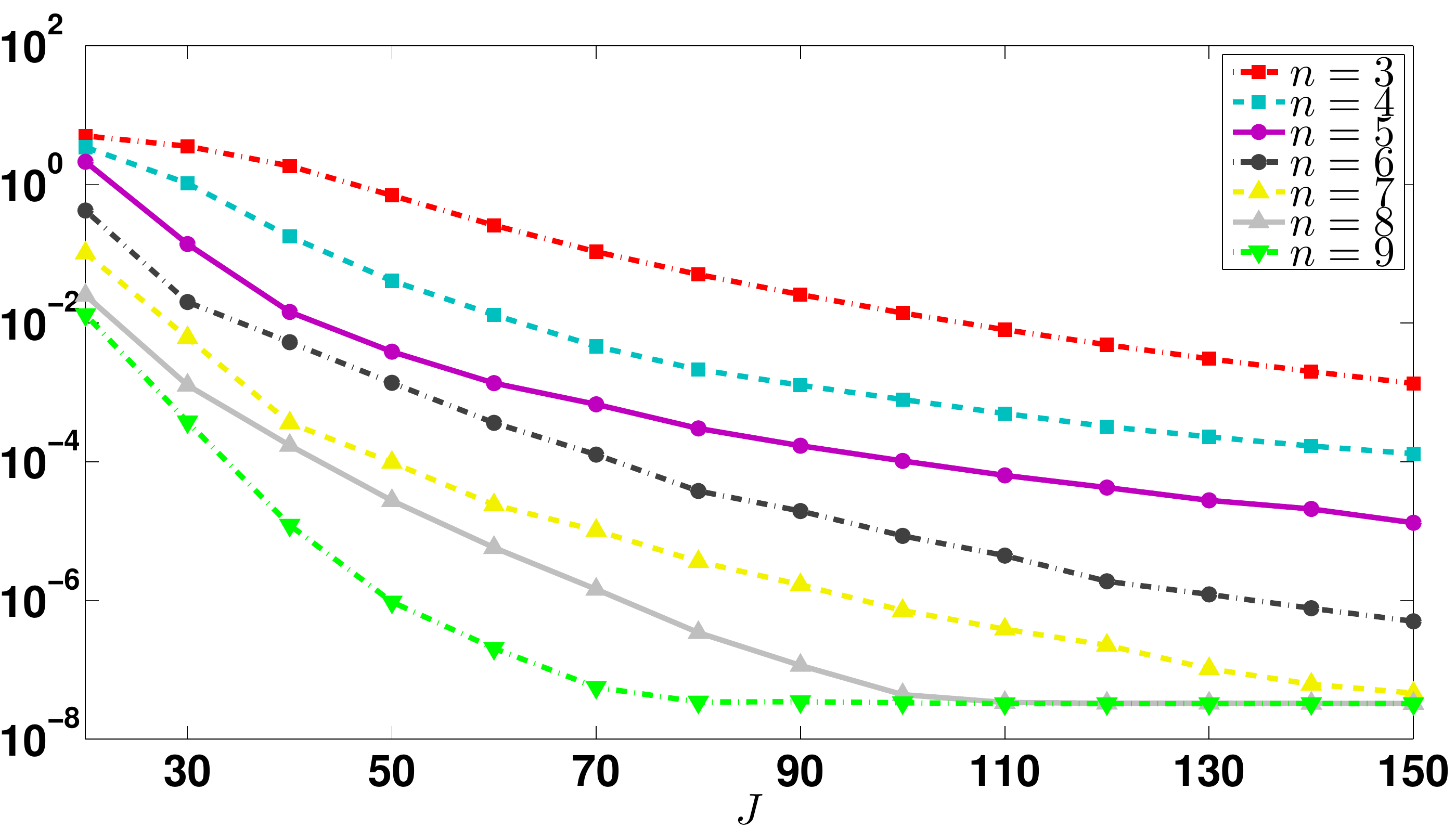}
    } \small{(a) $n=3,\dots,9$ and $J=20,30,\dots,150$} \\
    \end{minipage}\hfill
    \begin{minipage}[h]{0.49\linewidth}\center{
        \includegraphics[width=1\linewidth]{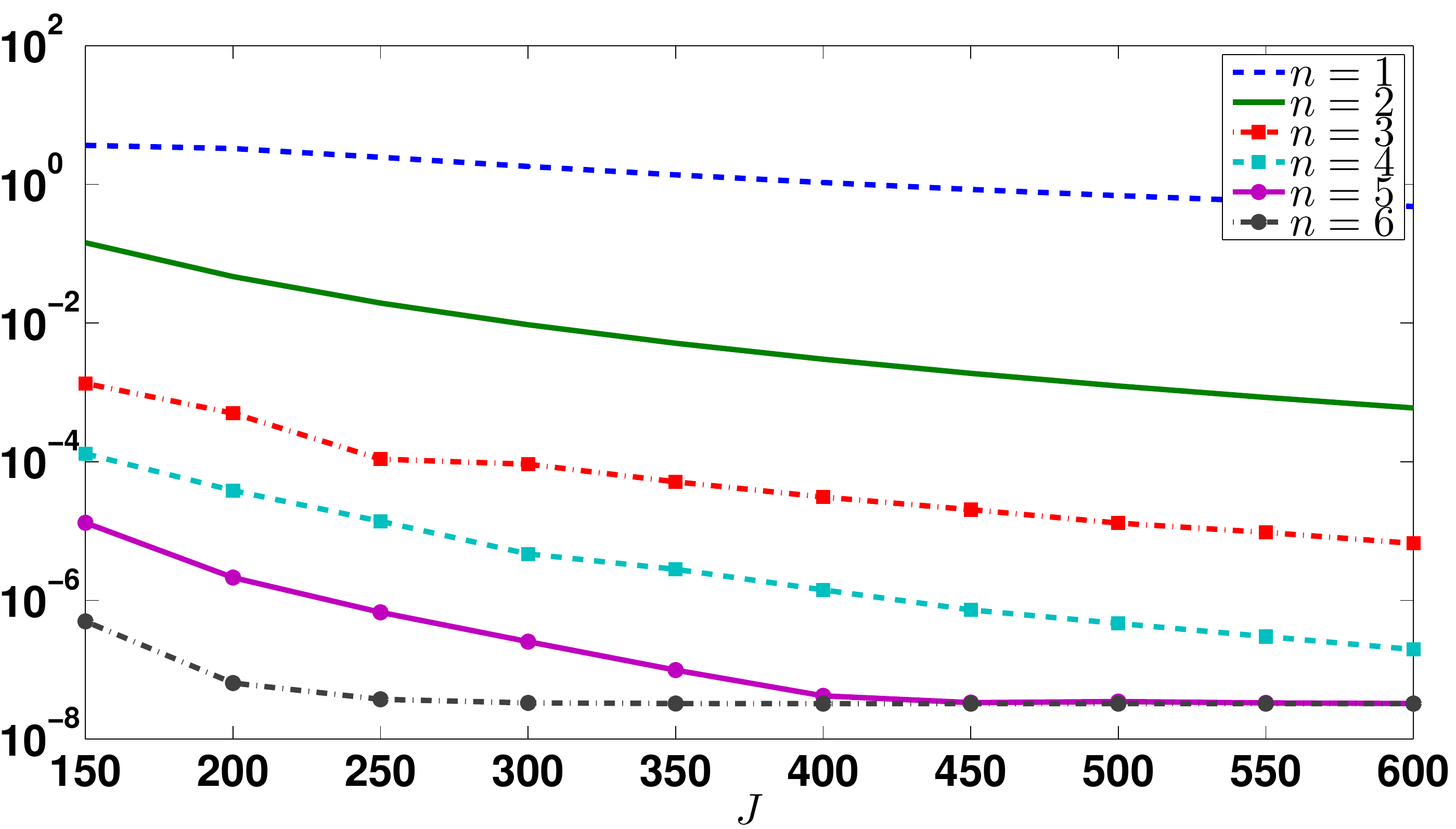}
    } \small{(b) $n=1,\dots,6$ and $J=150,200,\dots,600$} \\
    \end{minipage}
\caption{\small{Example 1. The errors $\max_{0\leq 3m\leq M}\|\psi^{3m}-\Psi_{3R}^{3m}\|_{C_h}$,
for $M=3000$, in dependence with $n$ and $J$}
\label{EX01r:3:FEM:MaxAbsError:C}}
\end{figure}

 Fig.~ \ref{EX01r:3:FEM:MaxAbsError:C} exhibits the behavior of the error
$\max_{0\leq 3m\leq M}\|\psi^{3m}-\Psi_{3R}^{3m}\|_{C_h}$,
for large $M=3000$, in dependence with $n=1,2,\ldots,9$ and $J=30,40,\ldots,150,$ $200,\ldots,600$. Similarly to the case without the extrapolation \cite{ZZ12}, the errors decrease monotonically and faster and faster in $J$ as $n$ grows. For the simplest and common in practice case $n=1$ (linear elements), unfortunately decreasing is especially slow and the error is unacceptable.
The advantage of the high degree elements over the low degree ones is obvious.
Once again the errors stabilize (in both norms) for $J\geq J_1(n)$ as soon as their lowest levels have been achieved, with $J_1(5)=450$, $J_1(6)=250$, $J_1(7)=150$, $J_1(8)=110$ and $J_1(9)=80$; clearly $J_1(n)$ decreases rather rapidly as $n$ grows. (Of course, the errors ultimately stabilize also for smaller $n$ but for much larger $J_1(n)>600$ absent on the figures.) The behavior of the similar errors in $L^2$ norm is quite close, with slightly better minimal values, and we omit their graphs.

\par In Table \ref{tab:EX01r:ExecutionTime}, we put the additional costs (in percents) that are required to compute $\Psi_{rR}$, $r=2,3,4$, in comparison with $\Psi^{(\tau)}$, for $n=5$ and several $J$ and $M$; for our computations,
the code in MATLAB R2013a is used on a quad-core processor PC.
The data in all the rows (except two for $J=300, 600$ and $M=3000$) are close
to above theoretical upper bound $(r-1)\cdot50\%$; some of them are slightly more than the bound (note that expressions \eqref{eq:richcost} do not take into account some costs like the computation of the stiffness and mass matrices and the discrete convolution kernel as well as details of exploiting PC hardware, etc.). The data in the exceptional two rows are essentially less than the bound that is also in agreement with costs \eqref{eq:richcost} for $J\ll M$.
\begin{table}\centering{
\begin {tabular}{rrccc}
\multicolumn {1}{c}{$J$}&
\multicolumn {1}{c}{$M$}&
\multicolumn {1}{c}{$r=2$}&
\multicolumn {1}{c}{$r=3$}&
\multicolumn {1}{c}{$r=4$}\\
\midrule %
$120$&$300$&$53.1\%$&$106.5\%$&$161.4\%$\\%
$120$&$600$&$51.3\%$&$102.2\%$&$153.6\%$\\%
$120$&$3\,000$&$45.9\%$&$\phantom{1}92.7\%$&$139.1\%$\\%
$300$&$300$&$52.2\%$&$106.9\%$&$162.5\%$\\%
$300$&$600$&$53.6\%$&$106.5\%$&$159.7\%$\\%
$300$&$3\,000$&$28.3\%$&$\phantom{1}63.7\%$&$102.1\%$\\%
$600$&$300$&$56.4\%$&$112.6\%$&$168.9\%$\\%
$600$&$600$&$54.9\%$&$108.4\%$&$163.0\%$\\%
$600$&$3\,000$&$27.6\%$&$\phantom{1}62.2\%$&$100.9\%$\\%
\end {tabular}%
}\\[2mm]
\caption{\small Example 1. The additional costs for computing $\Psi_{rR}$ versus $\Psi^{(\tau)}$} 
\label{tab:EX01r:ExecutionTime}
\end{table}
\smallskip\par\textbf{3.2}.
In Example 2, we treat the Cauchy problem \eqref{eq:se}, \eqref{eq:ic} for the piecewise constant potential $V=800\chi_I$, where $\chi_I$ is the characteristic function of the interval $I=(0.5,0.6)$,
and the initial function $\psi^0$ of form \eqref{FEM:s62}
with $x^{(0)}=-0.5$, $k=30$ and $\alpha=\frac{1}{120}$ now.
Thus tunneling through the discontinuous rectangular barrier is studied.
\par We choose $X=1.5$ and $T=0.09$. Now $|\psi^0(x)|<1E$$-13$ outside $\Omega$, and the discrete TBCs are posed at the both artificial boundaries $x=\pm X$.
A close example was considered in \cite{ZZ12}.
Looking ahead, notice that though the solution is not smooth in this and the next examples owing to the discontinuity of the potential,
nevertheless the Richardson extrapolation works well.
\par
\begin{figure}[htbp]
    \begin{minipage}[h]{0.3\linewidth}\center{
        \includegraphics[width=1\linewidth]{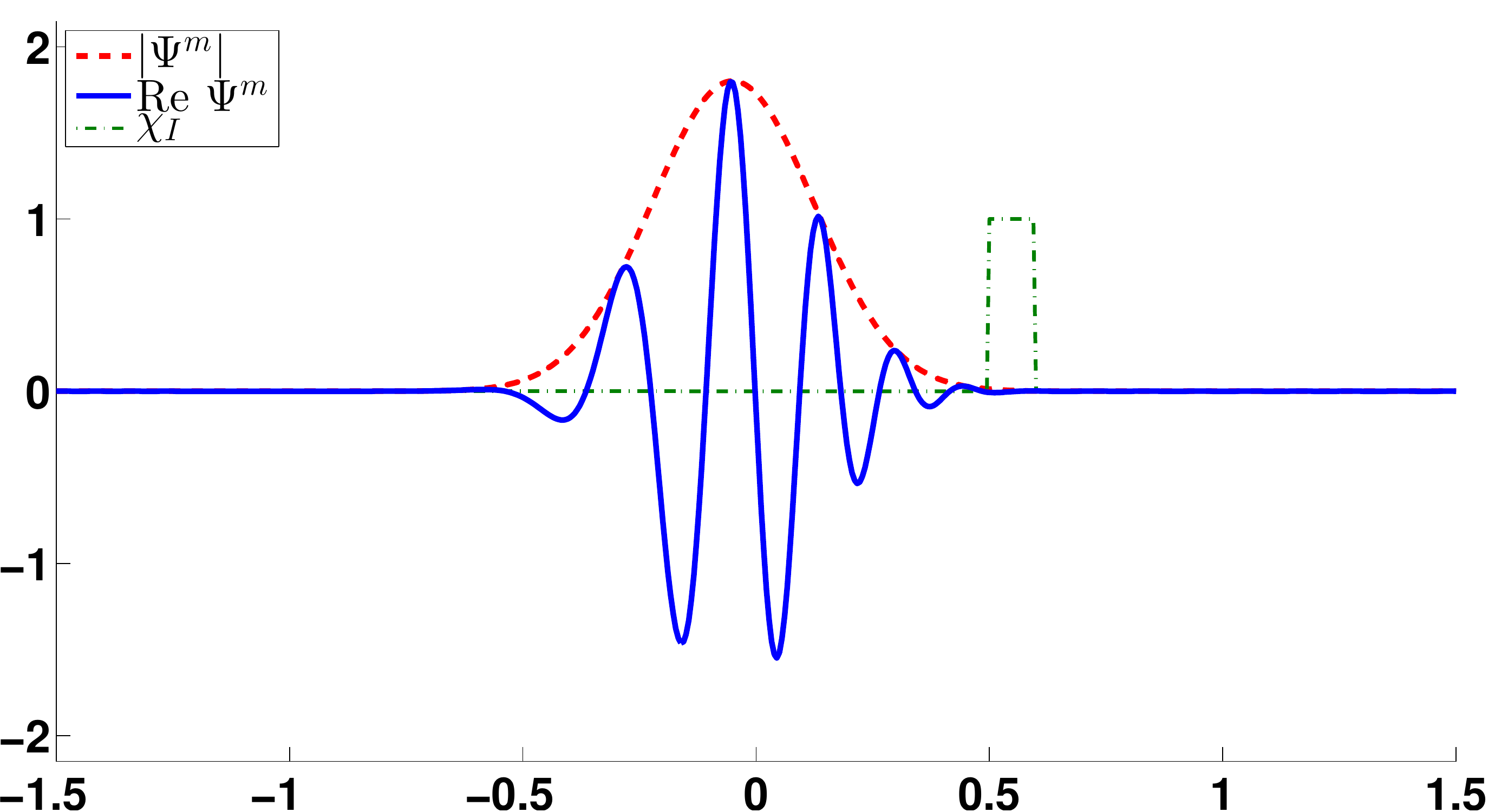}
    } \small{(a) $t_m= 0.0075$, $m=48$} \\
    \end{minipage}\hfill
    \begin{minipage}[h]{0.3\linewidth}\center{
        \includegraphics[width=1\linewidth]{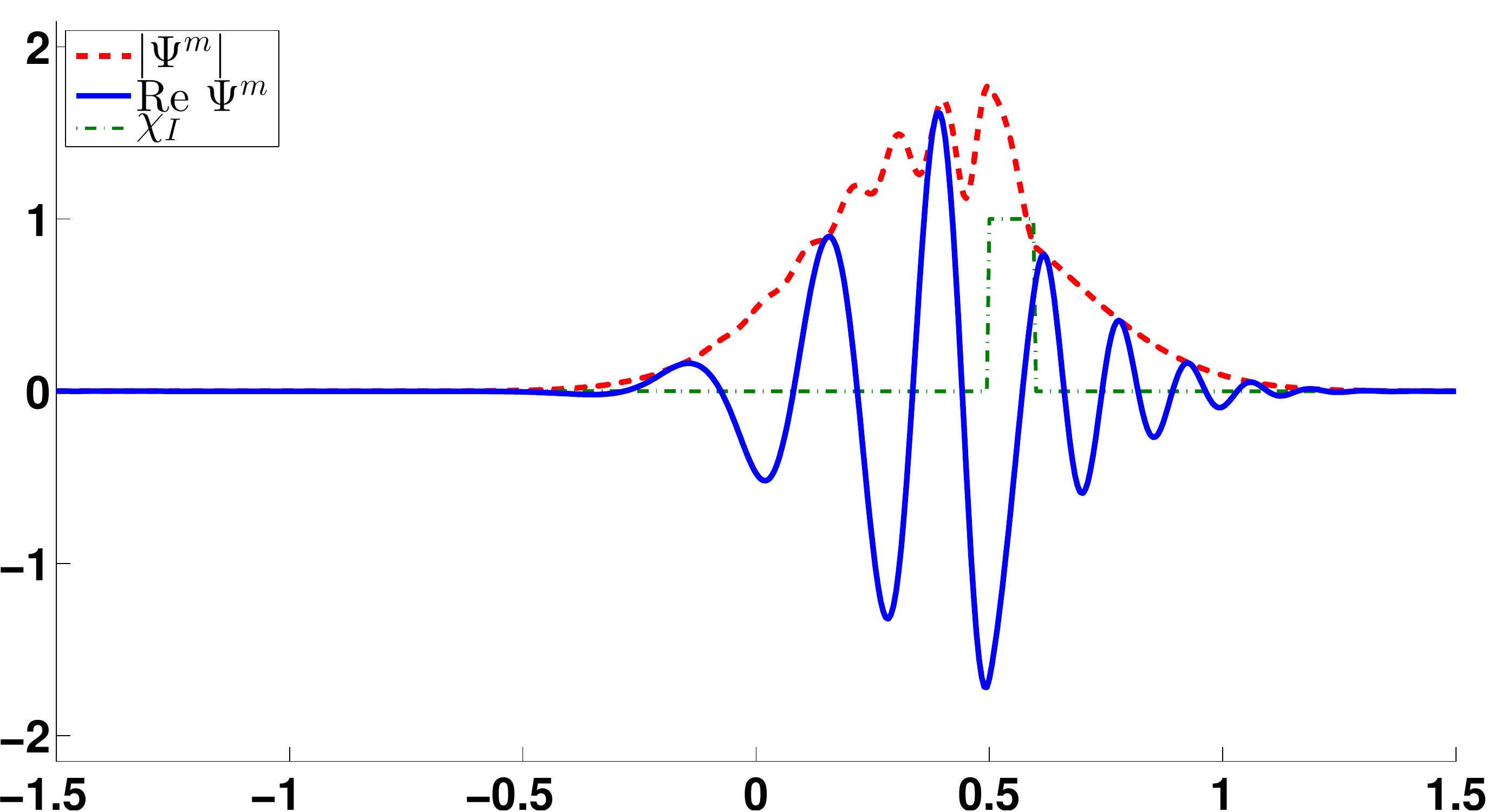}
    } \small{(b) $t_m=0.015$, $m=96$} \\
    \end{minipage}\hfill
    \begin{minipage}[h]{0.3\linewidth}\center{
        \includegraphics[width=1\linewidth]{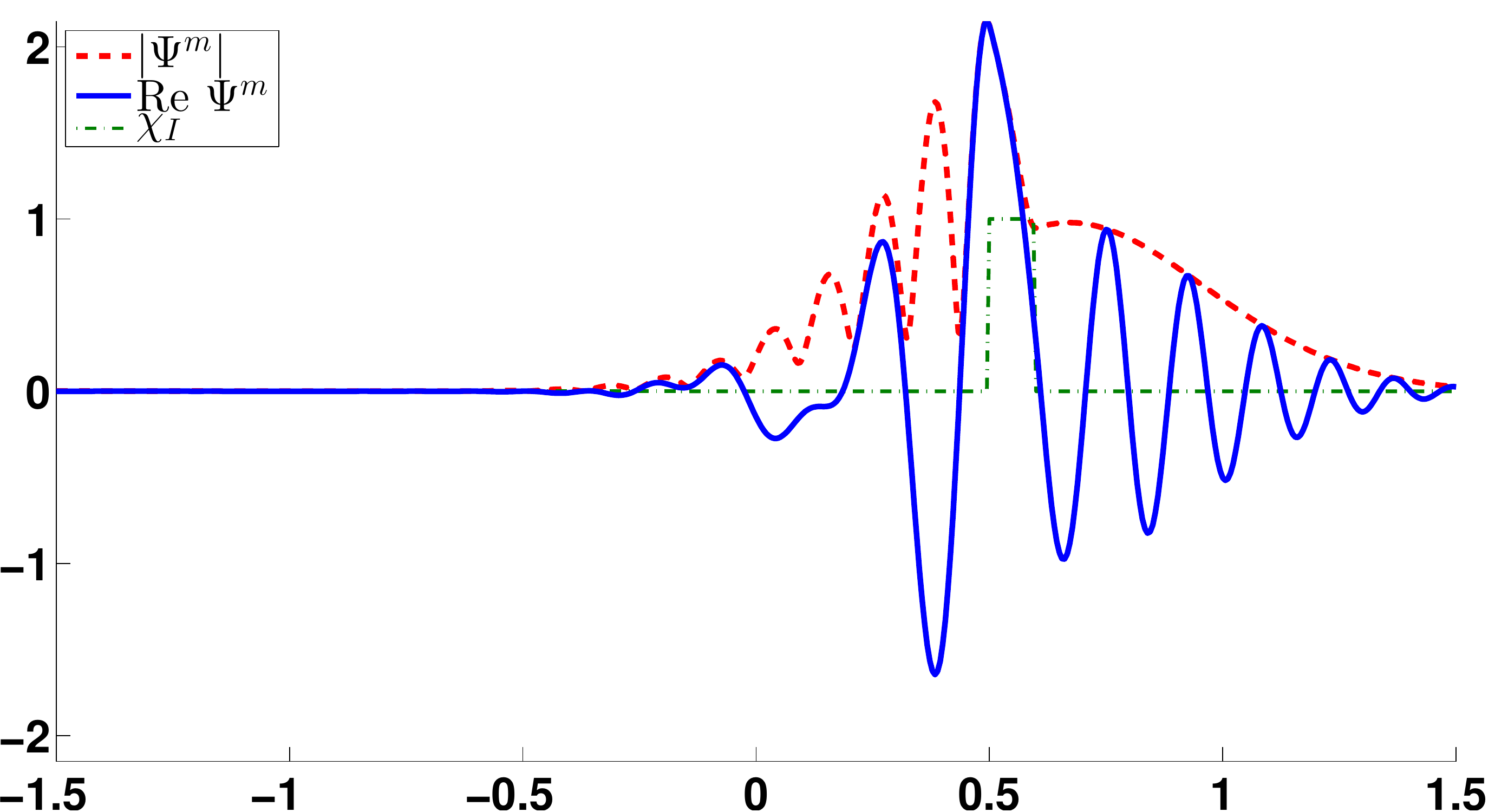}
    } \small{(c) $t_m=0.0187$, $m=120$} \\[3mm]
    \end{minipage}\vfill
    \begin{minipage}[h]{0.3\linewidth}\center{
        \includegraphics[width=1\linewidth]{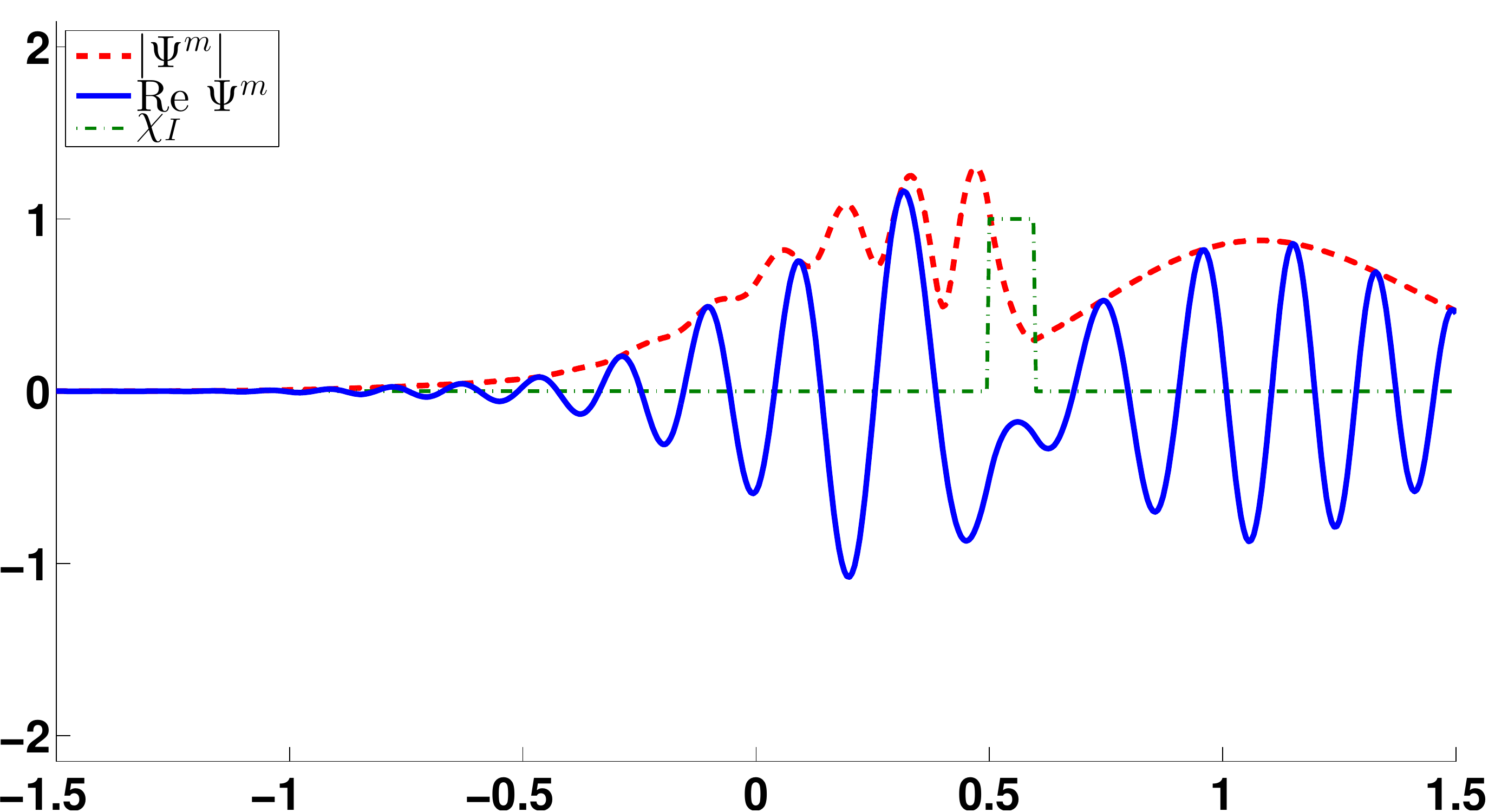}
    } \small{(d) $t_m=0.0250$, $m=160$} \\
    \end{minipage}\hfill
    \begin{minipage}[h]{0.3\linewidth}\center{
        \includegraphics[width=1\linewidth]{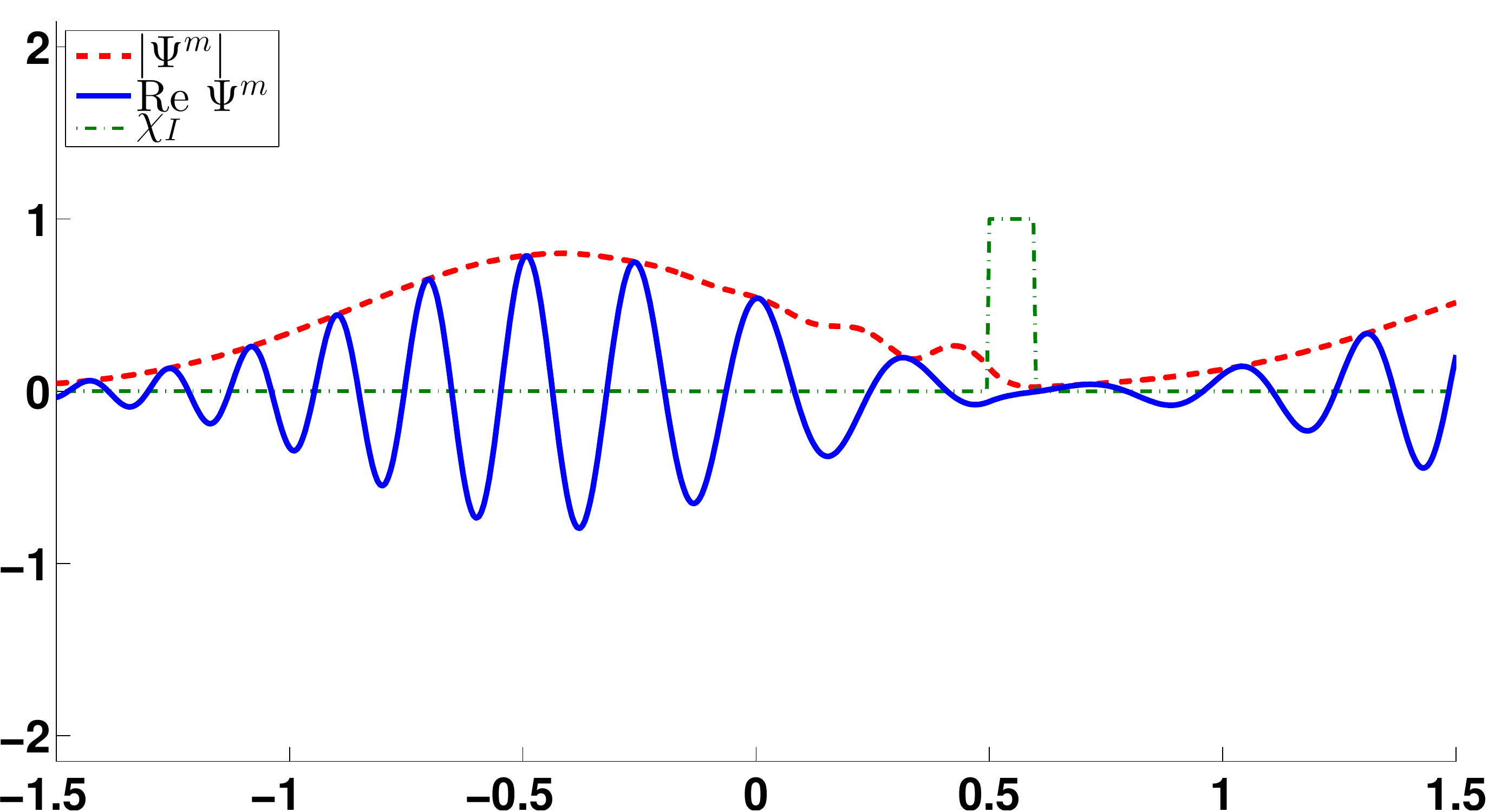}
    } \small{(e) $t_m=0.0375$, $m=240$} \\
    \end{minipage}\hfill
    \begin{minipage}[h]{0.3\linewidth}\center{
        \includegraphics[width=1\linewidth]{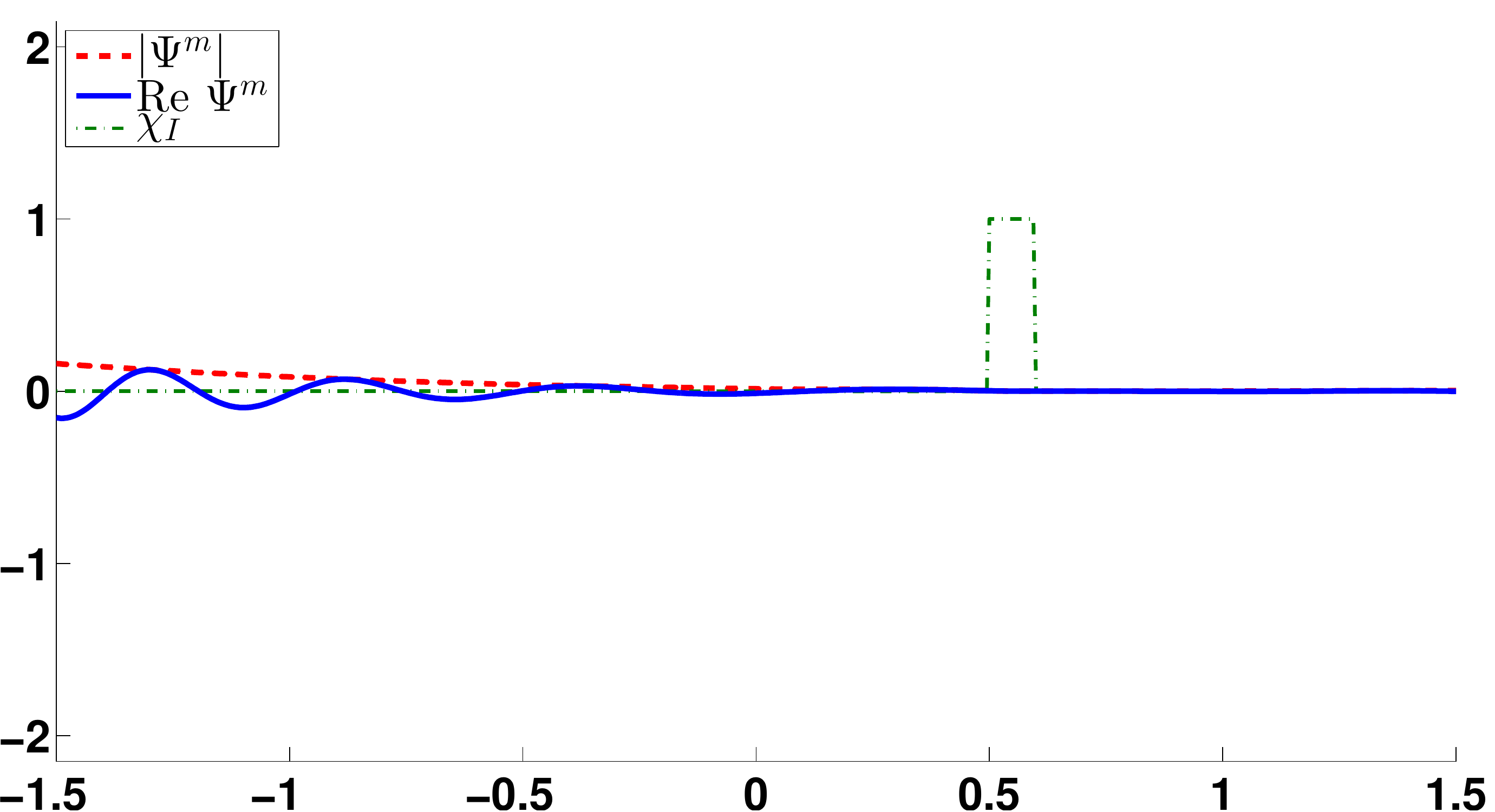}
    } \small{(f) $t_m=0.009$, $m=576$} \\[3mm]
    \end{minipage}\vfill
\caption{\small{Example 2. $|\Psi_{4R}^{m}|$ and $\Rea\Psi_{4R}^{m}$, for $n=9$ and $(J,M)=(60,576)$, and the scaled $V$}
\label{fig:EX22r:Solution}}
\end{figure}

 The behavior of $|\psi|$ and $\Rea\psi$ is shown on Fig.~\ref{fig:EX22r:Solution}. The wave moves to the right toward the barrier, interacts with it and then is divided into two comparable reflected and transmitted parts moving in the opposite directions.
The solution is represented by $\Psi_{4R}^{4m}$, for high $n=9$ but $(J,M)=(60,576)$ only, with a suitable uniform accuracy (see Table \ref{tab:EX22r} below).
\par Note that, for $J=30$,  $\bar{I}$ consists of exactly one element; that is why below our $J$ are multiples of 30. Since any simple analytical form of the exact solution $\psi$ is not known, below its role is played by the pseudo-exact solution $\Psi_{4R}^{4m}$ computed for high $n=9$, $J=150$ and large $M=36864$.
\par
\begin{figure}[htbp]
    \begin{minipage}[h]{0.49\linewidth}\center{
        \includegraphics[width=1\linewidth]{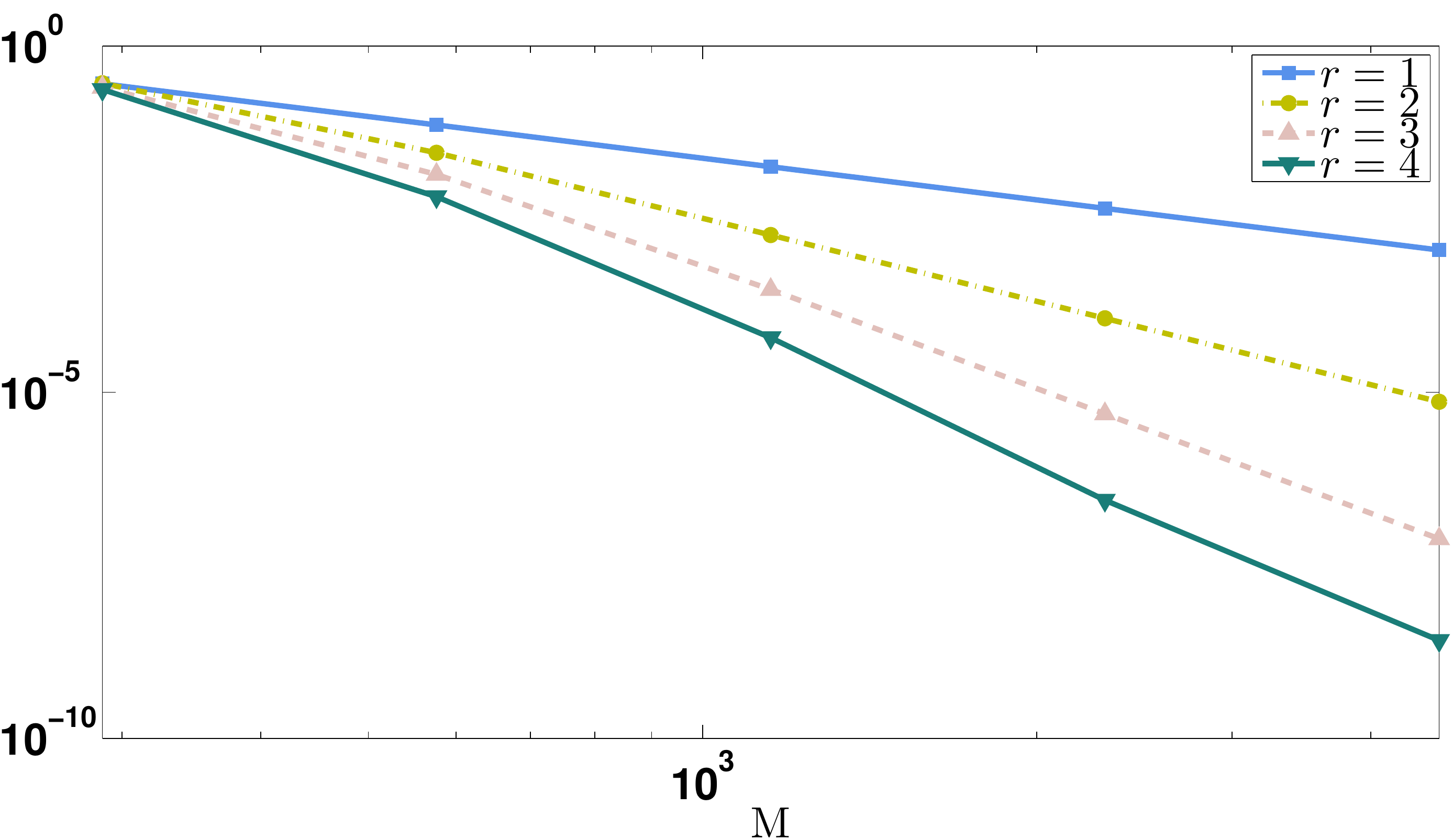}
    } \small{(a) in $L^2$ space norm} \\
    \end{minipage}\hfill
    \begin{minipage}[h]{0.49\linewidth}\center{
        \includegraphics[width=1\linewidth]{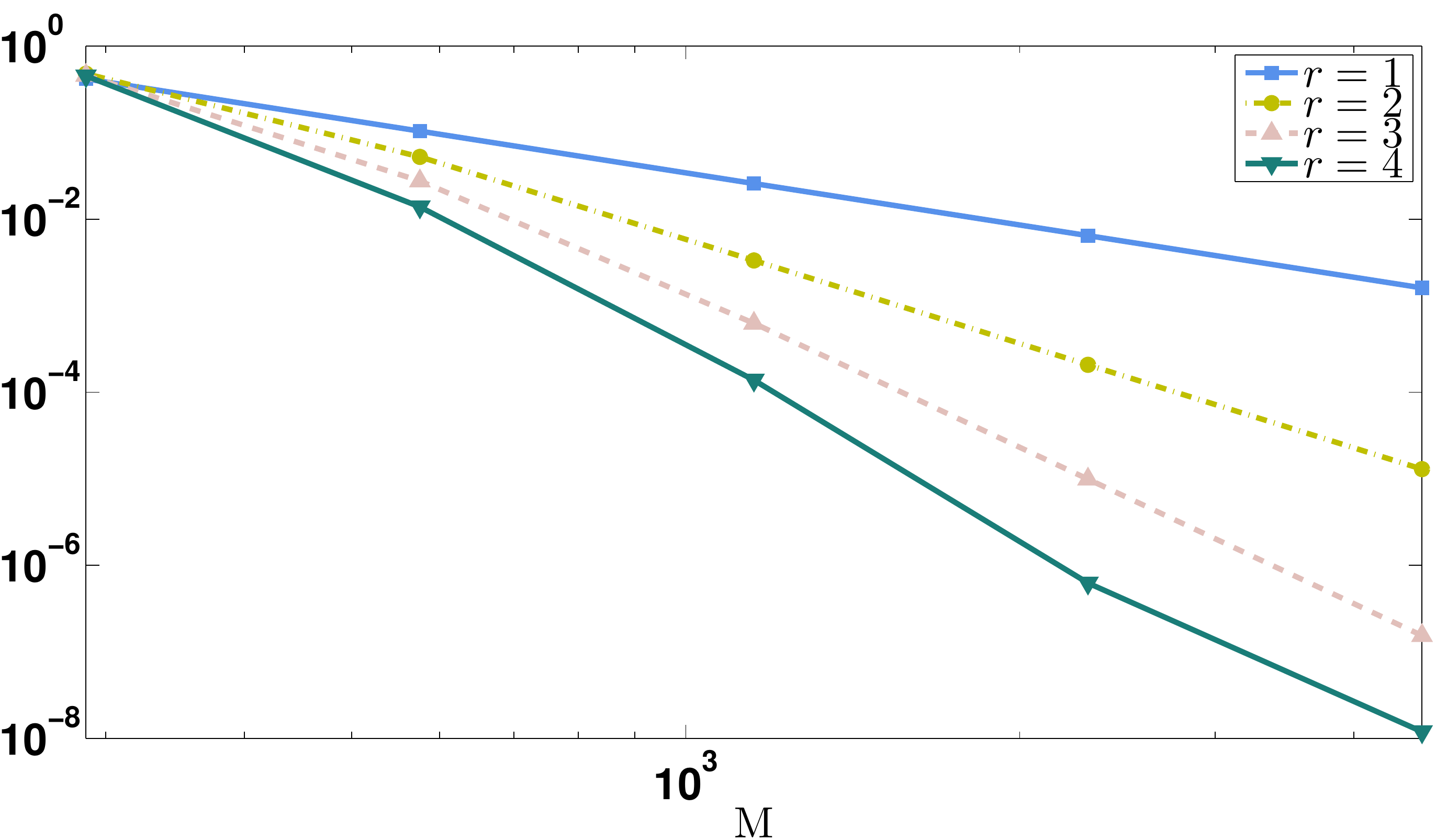}
    } \small{(b) in $C$ space norm} \\
    \end{minipage}
\caption{\small{Example 2. The errors $\max_{0\leq rm\leq M}\|\psi^{rm}-\Psi_{rR}^{rm}\|$,
for $n=9$ and $J=60$, in dependence with $r=1,2,3,4$ and $M=288\cdot 2^q$, $q=0,1,2,3,4$}
\label{fig:EX22r:N=9:J=60:MaxAbsError}}
\end{figure}

 On Fig. \ref{fig:EX22r:N=9:J=60:MaxAbsError}, we present the errors
$E_r^M=\max_{0\leq rm\leq M}\|\psi^{rm}-\Psi_{rR}^{rm}\|$,
for $n=9$ and $J=60$, in dependence with $r=1,2,3,4$ and $M=M_q=288\cdot 2^q$, $q=0,1,2,3,4$
(recall that $\Psi_{1R}=\Psi^{(\tau)}$).
Once again, for $r=1$, the errors decay most slowly. They decay faster and faster as $r$ grows.
Notice that the behavior of the corresponding relative errors
$\max_{0\leq rm\leq M}\|\psi^{rm}-\Psi_{rR}^{rm}\|/\|\psi^{rm}\|$ is quite similar.
\par The same data as on Fig. \ref{fig:EX22r:N=9:J=60:MaxAbsError} are put into
Table \ref{tab:EX22r}
together with the corresponding ratios $E_r^{M_q}/E_r^{M_{q-1}}$ for clarity.
Note that all the errors are large for the smallest $M=288$.
Once again the errors decay most slowly for $r=1$ but they decay faster and faster as $r$ grows.
We can compare the ratios with their theoretically predicted values $2^{2r}=(4,16,64,256)$ respectively for $r=1,2,3,4$ (see \eqref{eq:errexp} and \eqref{eq:errrich}). We see their closeness for all
$q$ for $r=1$, $q\geq 2$ for $r=2$,
$q=3$ and $4$ for $r=3$ as well as $q=3$ for $r=4$;
in any case, the ratios grows significantly as $r$ increases for fixed $q$ (excepting the last value for $q=r=4$ in $C$ space norm).
\begin{table}\centering{
\begin{tabular}{rrcrcrcrc}
\multicolumn {1}{c}{$M$}&
\multicolumn {2}{c}{$r=1$}&
\multicolumn {2}{c}{$r=2$}&
\multicolumn {2}{c}{$r=3$}&
\multicolumn {2}{c}{$r=4$}\\
\midrule %
$288$&$0.28$&--&$0.29$&--&$0.25$&--&$0.24$&--\\%
$576$&$7.23E{-}2$&$3.94$&$2.86E{-}2$&$10.11$&$1.39E{-}2$&$18.07$&$6.64E{-}3$&$35.42$\\%
$1\,152$&$1.81E{-}2$&$4$&$1.86E{-}3$&$15.39$&$3.03E{-}4$&$45.89$&$6.12E{-}5$&$108.51$\\%
$2\,304$&$4.51E{-}3$&$4$&$1.16E{-}4$&$16$&$4.86E{-}6$&$62.43$&$2.75E{-}7$&$222.38$\\%
$4\,608$&$1.13E{-}3$&$4$&$7.26E{-}6$&$16.01$&$7.60E{-}8$&$63.95$&$2.59E{-}9$&$106.18$\\%
\end {tabular}%
\\[2mm]
\small{(a) in $L^2$ space norm}\\[2mm]
\begin{tabular}{rrcrcrcrc}
\multicolumn {1}{c}{$M$}&
\multicolumn {2}{c}{$r=1$}&
\multicolumn {2}{c}{$r=2$}&
\multicolumn {2}{c}{$r=3$}&
\multicolumn {2}{c}{$r=4$}\\
\midrule %
$288$&$0.42$&--&$0.48$&--&$0.46$&--&$0.45$&--\\%
$576$&$0.1$&$4.01$&$5.22E{-}2$&$9.17$&$2.75E{-}2$&$16.6$&$1.38E{-}2$&$32.58$\\%
$1\,152$&$2.57E{-}2$&$4.02$&$3.32E{-}3$&$15.69$&$6.24E{-}4$&$44.12$&$1.38E{-}4$&$100.38$\\%
$2\,304$&$6.42E{-}3$&$4.01$&$2.07E{-}4$&$16.09$&$9.89E{-}6$&$63.12$&$6.23E{-}7$&$221.26$\\%
$4\,608$&$1.61E{-}3$&$4$&$1.29E{-}5$&$16.03$&$1.54E{-}7$&$64.15$&$1.19E{-}8$&$52.36$\\%
\end {tabular}%
\\[2mm]
\small{(b) in $C$ space norm}}\\[2mm]
\caption{\small Example 2. The errors $E_r^M=\max_{0\leq rm\leq M}\|\psi^{rm}-\Psi_{rR}^{rm}\|$
and their ratios $E_r^{M_q}/E_r^{M_{q-1}}$, for $n=9$ and $J=60$, in dependence with
$r$ and $M_q=288\cdot 2^q$, $q=0,1$,...,$4$}
\label{tab:EX22r}
\end{table}

\par
\begin{figure}[htbp]
    \begin{minipage}[h]{0.49\linewidth}\center{
        \includegraphics[width=1\linewidth]{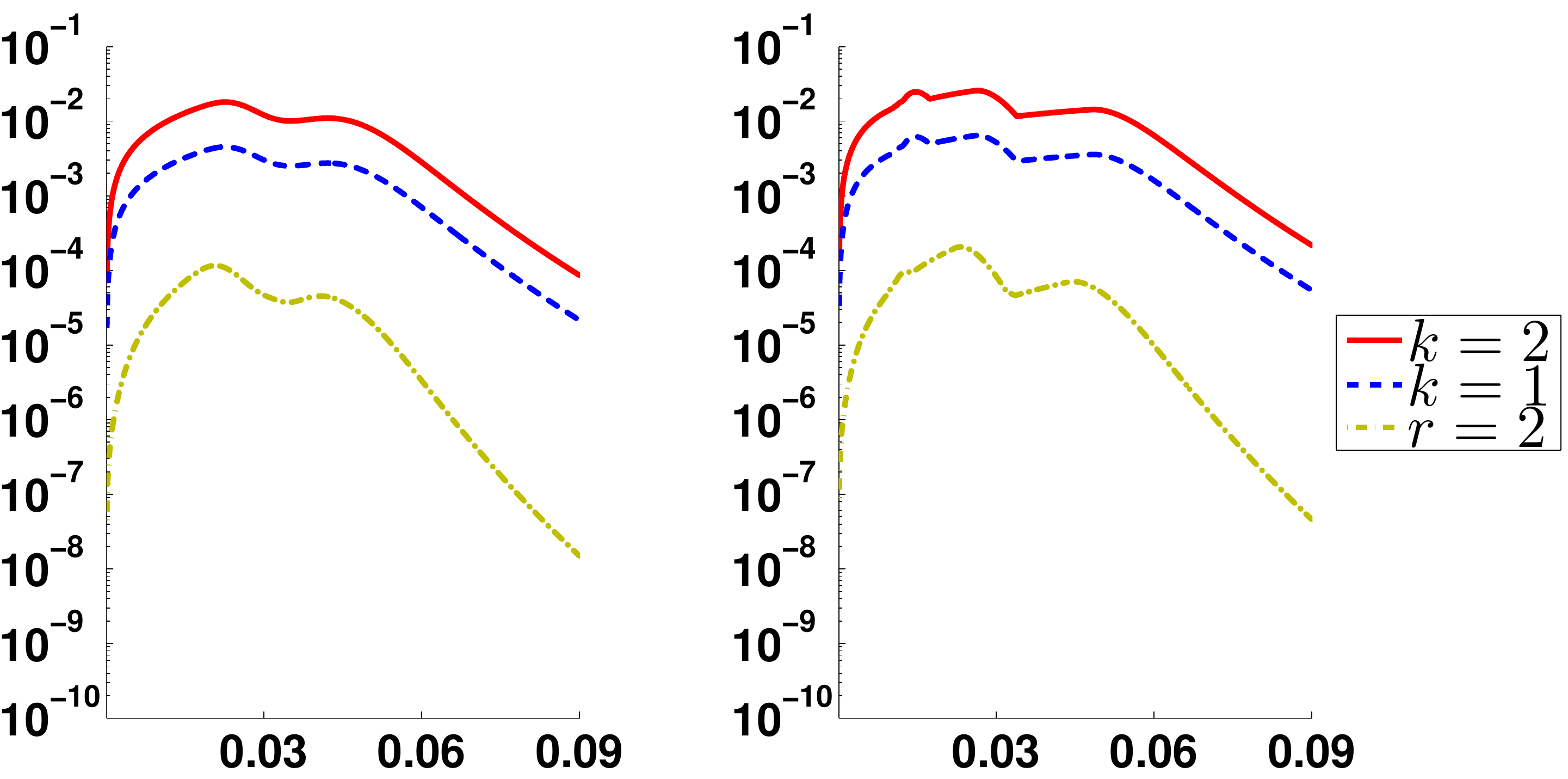}
    } \small{(a) in $L^2$ (left) and $C$ (right) norms, $r=2$} \\
    \end{minipage}\hfill
    \begin{minipage}[h]{0.49\linewidth}\center{
        \includegraphics[width=1\linewidth]{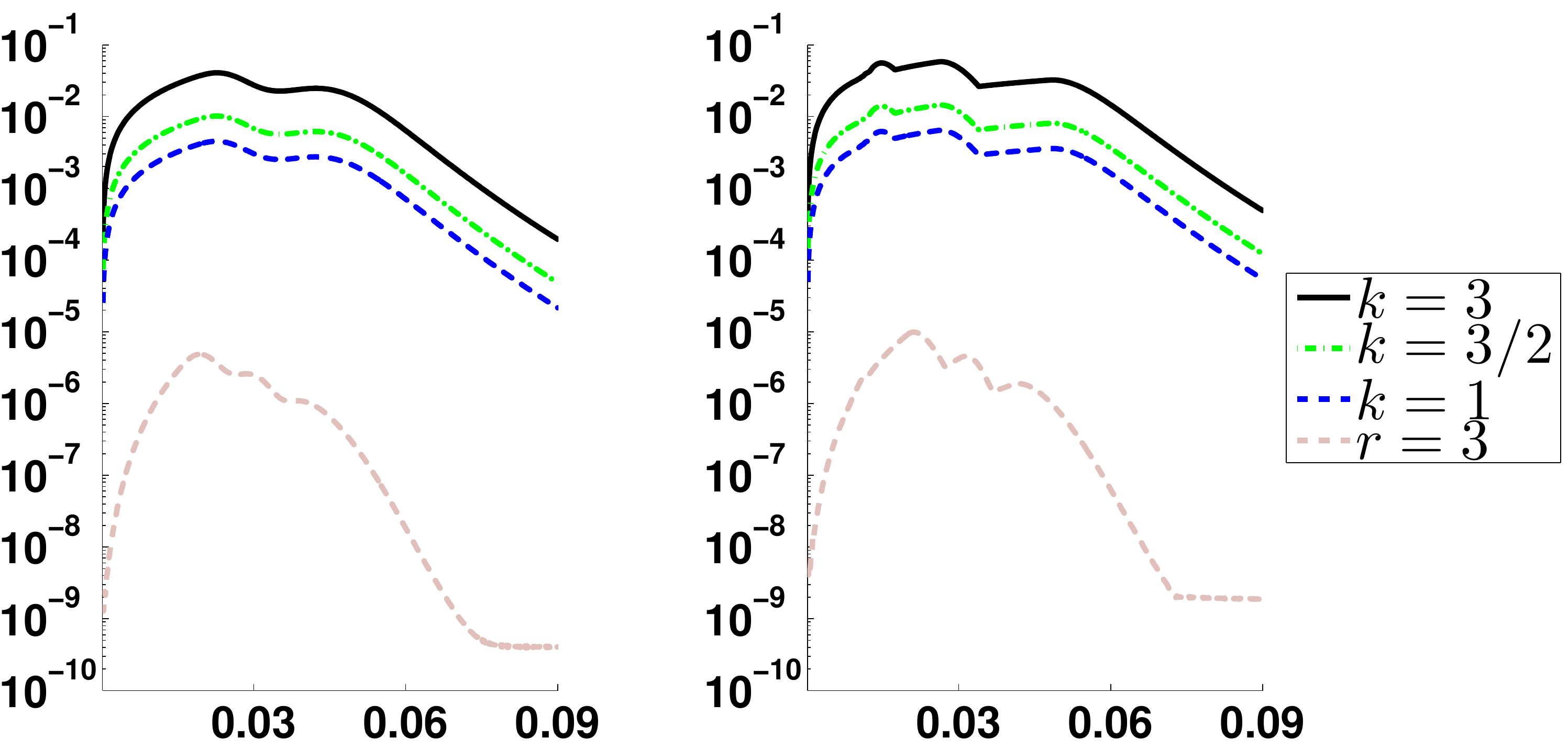}
    } \small{(b) in $L^2$ (left) and $C$ (right) norms, $r=3$} \\
    \end{minipage}
\caption{\small{Example 2. The errors: (a) $\|\psi^{2m}-\Psi^{(k\tau),\,2m}\|$ for $k=1,2$
and $\|\psi^{2m}-\Psi_{2R}^{2m}\|$ ($r=2$), and
(b) $\|\psi^{3m}-\Psi^{(k\tau),\,3m}\|$ for $k=1,3/2,3$
and $\|\psi^{3m}-\Psi_{3R}^{3m}\|$ ($r=3$),
both for $n=9$ and $(J,M)=(60,2304)$, in dependence with $t_m$}
\label{fig:EX22r:Error}}
\end{figure}

 Fig. \ref{fig:EX22r:Error} demonstrates the slow monotone decay of the error $\|\psi^{rm}-\Psi^{(k\tau),\,rm}\|$ as $k$ decreases and much less error (by several orders of magnitude) of the Richardson extrapolation $\|\psi^{rm}-\Psi_{rR}^{rm}\|$, for $r=2$ and especially for $r=3$, on the whole time segment $[0,T]$.

\par
\begin{figure}[htbp]
    \begin{minipage}[h]{0.49\linewidth}\center{
        \includegraphics[width=1\linewidth]{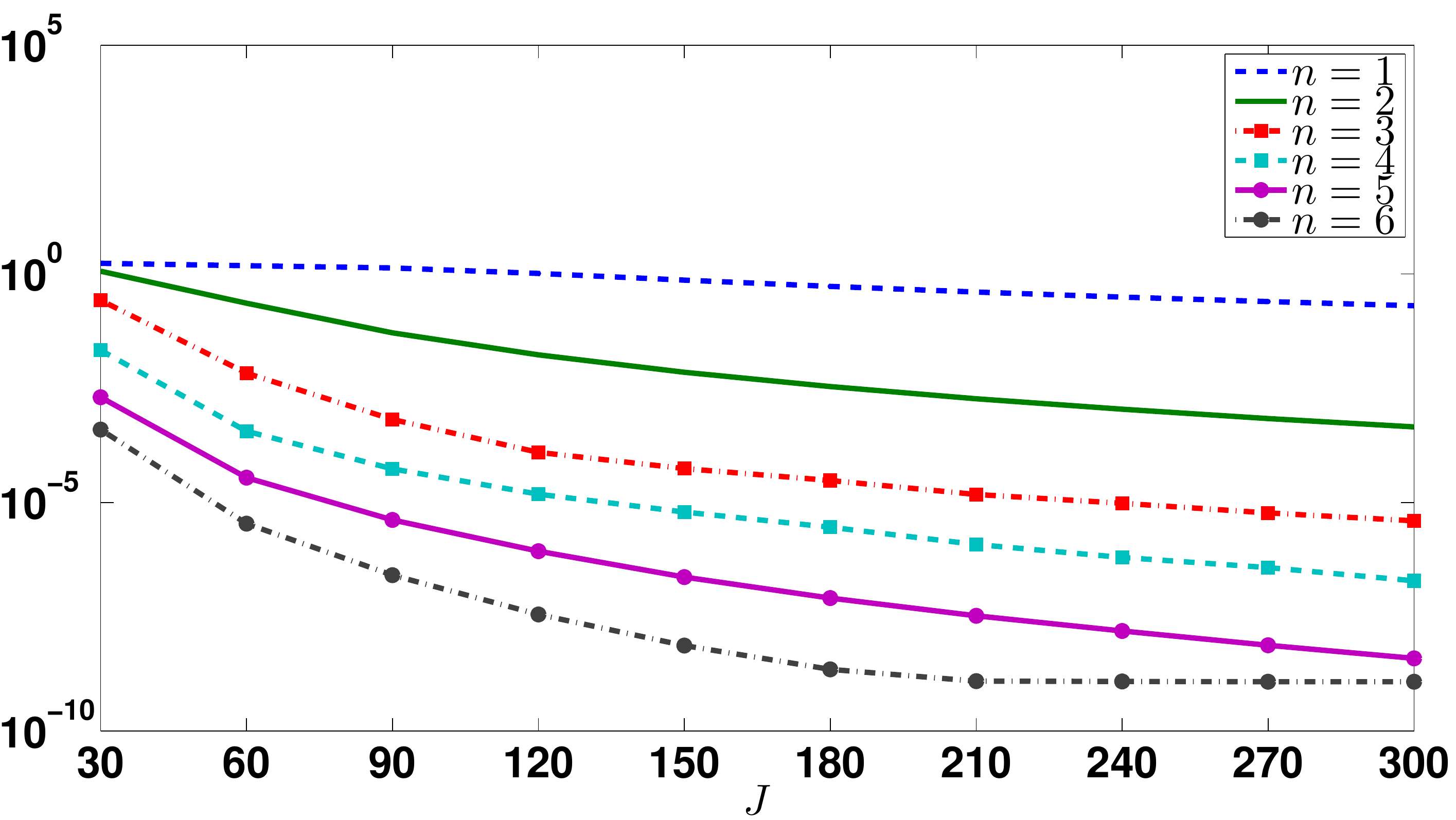}
    } \small{(a) in $L_2$ space norm}
    \end{minipage}\hfill
    \begin{minipage}[h]{0.49\linewidth}\center{
        \includegraphics[width=1\linewidth]{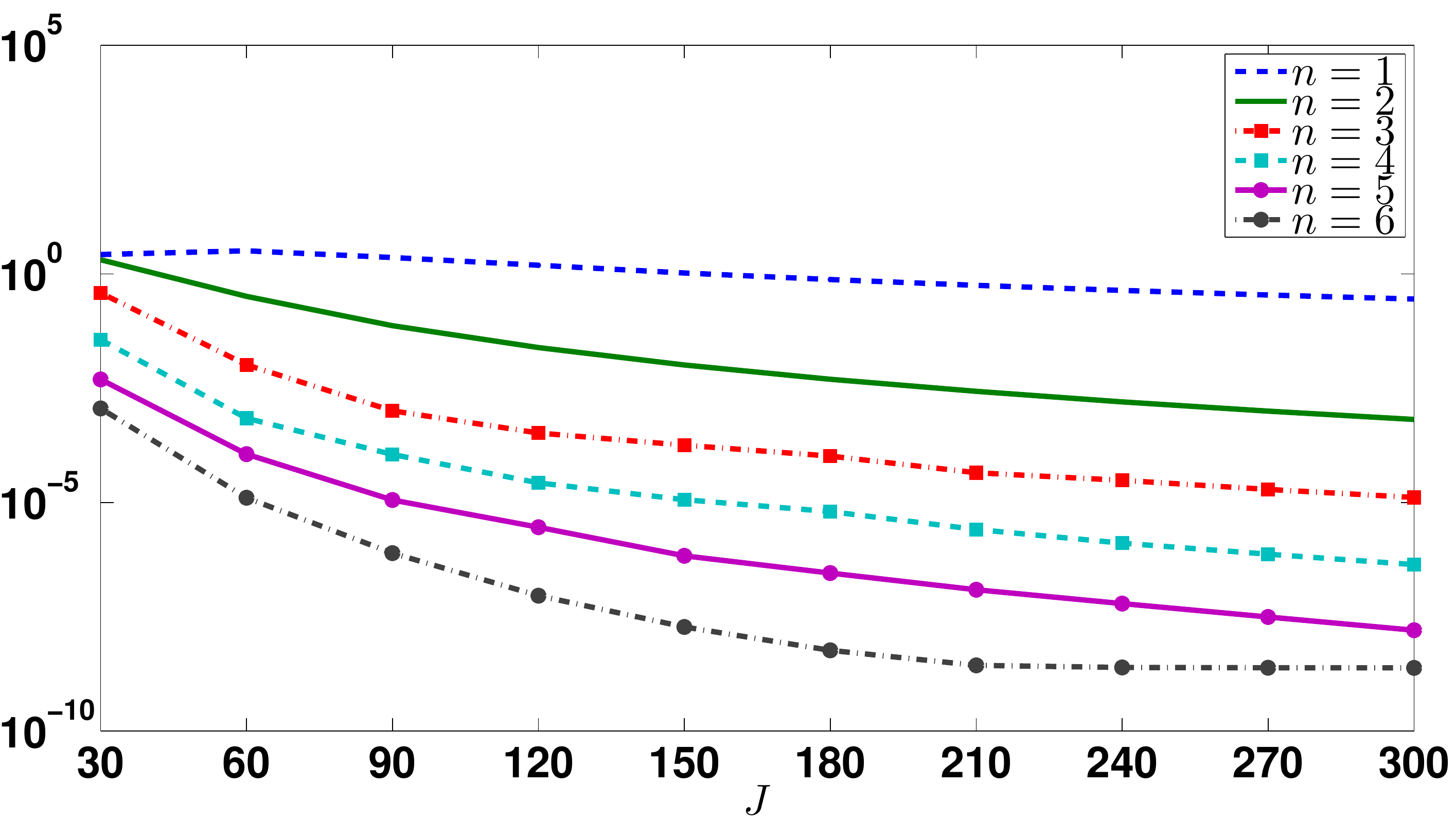}
    } \small{(b) in $C$ space norm}
    \end{minipage}
\caption{\small{Example 2. The errors $\max_{0\leq 3m\leq M}\|\psi^{3m}-\Psi_{3R}^{3m}\|$, for $M=9216$, in dependence with $n=1,2,\ldots,6$ and $J=30,60,\ldots,300$}
\label{fig:EX22r:MaxAbsError:r=3}}
\end{figure}

\begin{figure}[htbp]
    \begin{minipage}[h]{0.49\linewidth}\center{
        \includegraphics[width=1\linewidth]{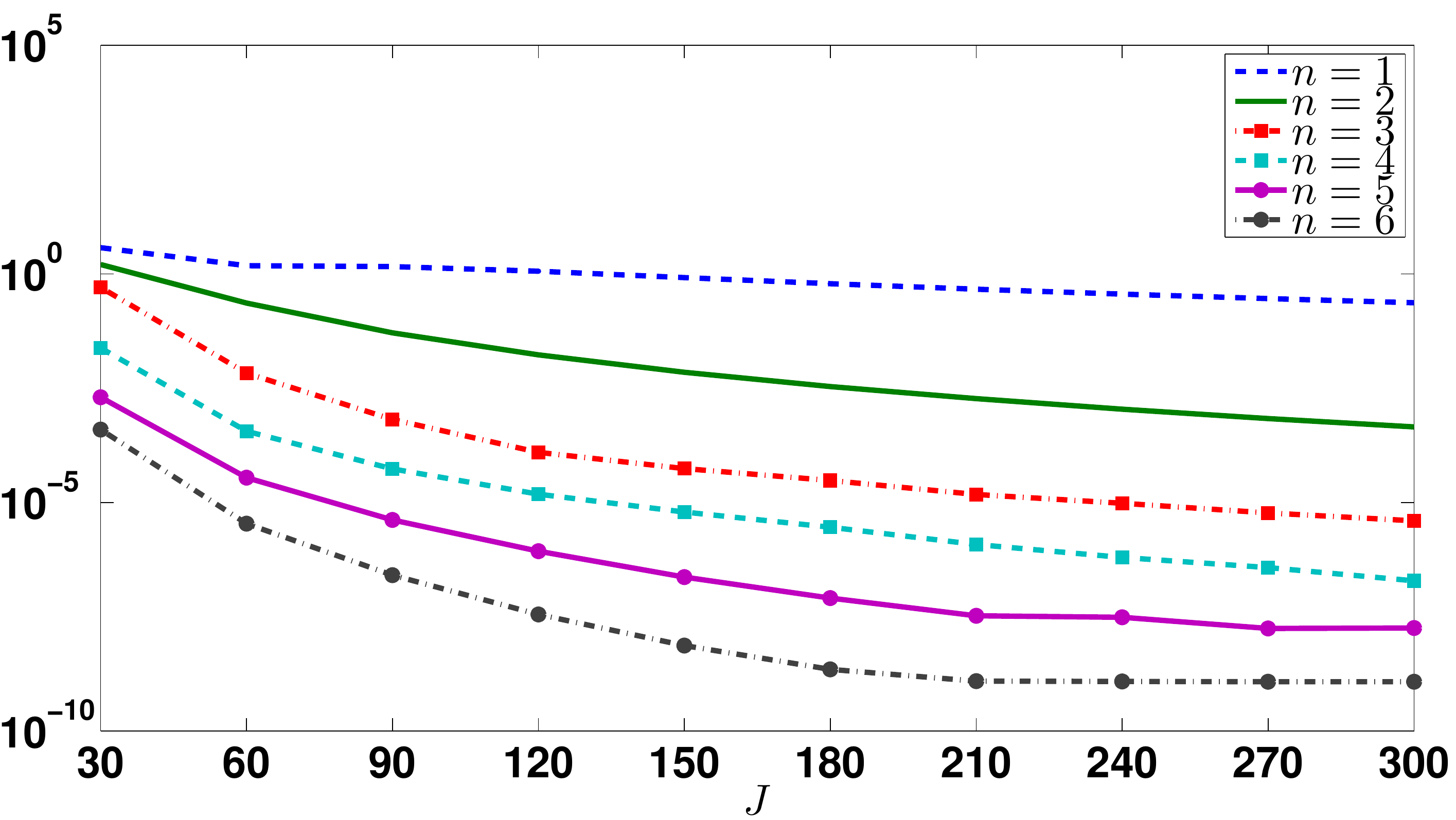}
    } \small{(a) in $L_2$ space norm}
    \end{minipage}\hfill
    \begin{minipage}[h]{0.49\linewidth}\center{
        \includegraphics[width=1\linewidth]{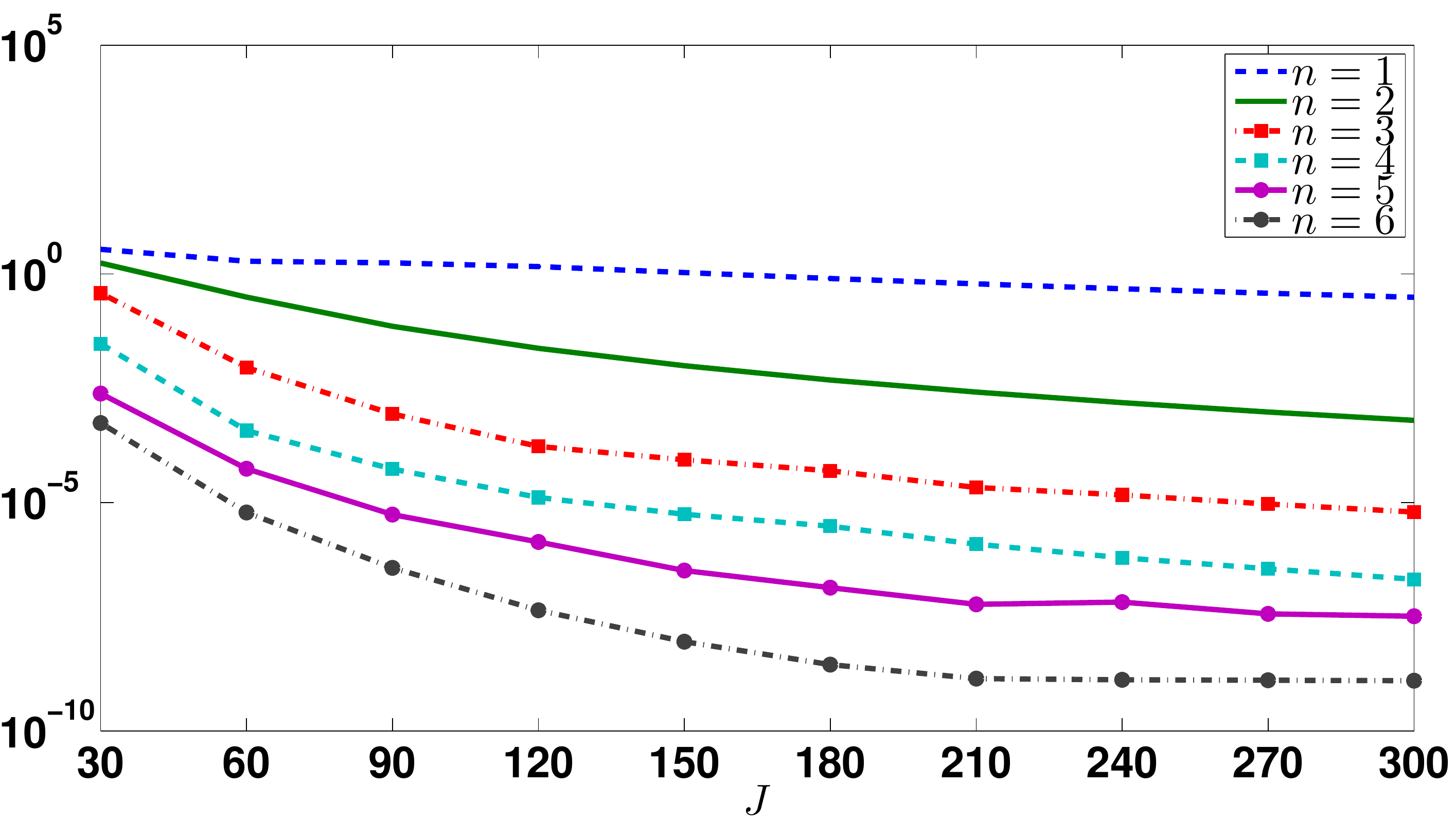}
    } \small{(b) in $C$ space norm}
    \end{minipage}
\caption{\small{Example 2. The relative errors $\max_{0\leq 3m\leq M}\|\psi^{3m}-\Psi_{3R}^{3m}\|/\|\psi^{3m}\|$, for $M=9216$, in dependence with $n=1,2,\ldots,6$ and $J=30,60,\ldots,300$}
\label{fig:EX22r:MaxRelError:r=3}}
\end{figure}

 Fig.~\ref{fig:EX22r:MaxAbsError:r=3} exhibits the behavior of the errors
$\max_{0\leq 3m\leq M}\|\psi^{3m}-\Psi_{3R}^{3m}\|$, in $L^2$ and $C$ norms
for large $M=9216$, in dependence with $n=1,2,\ldots,6$ and $J=30,60,\ldots,300$.
As above at the absence of the potential, the errors decrease monotonically as $J$ grows. They also decrease rapidly as $n$ grows whereas, for $n=1$ (linear elements), decreasing is very slow and the error is still unacceptable.
The errors stabilize (now due to the fixed value of $M$), for $n=6$ and $J\geq J_1(6)=210$.
Interestingly, the behavior of the corresponding relative errors is essentially quite similar (except for $n=5$), see
Fig. \ref{fig:EX22r:MaxRelError:r=3}.
\smallskip
\par
\begin{figure}[htbp]
    \begin{minipage}[h]{0.3\linewidth}\center{
        \includegraphics[width=1\linewidth]{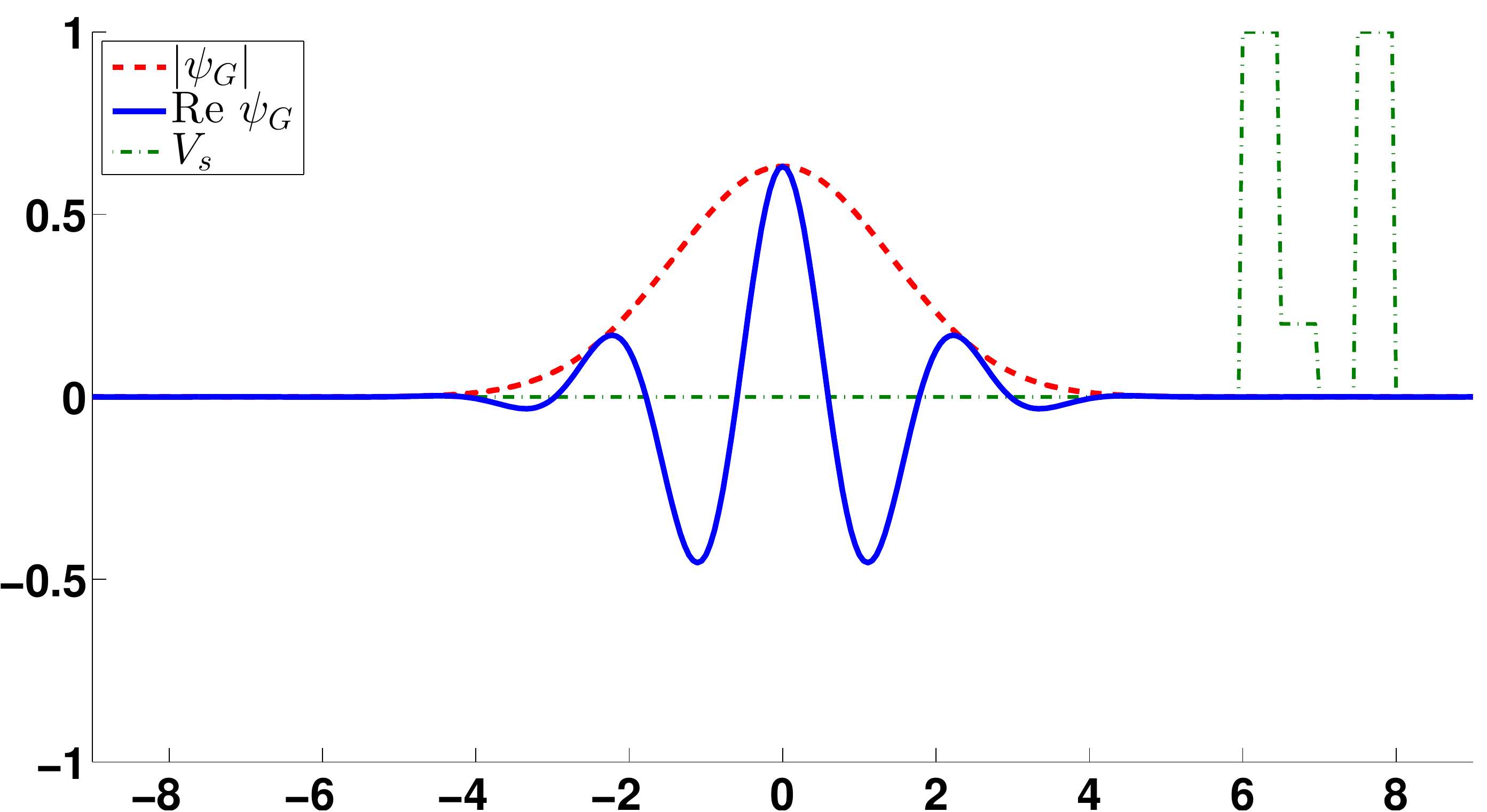}
    } \small{(a) $t_m=0$, $m=0$} \\
    \end{minipage}\hfill
    \begin{minipage}[h]{0.3\linewidth}\center{
        \includegraphics[width=1\linewidth]{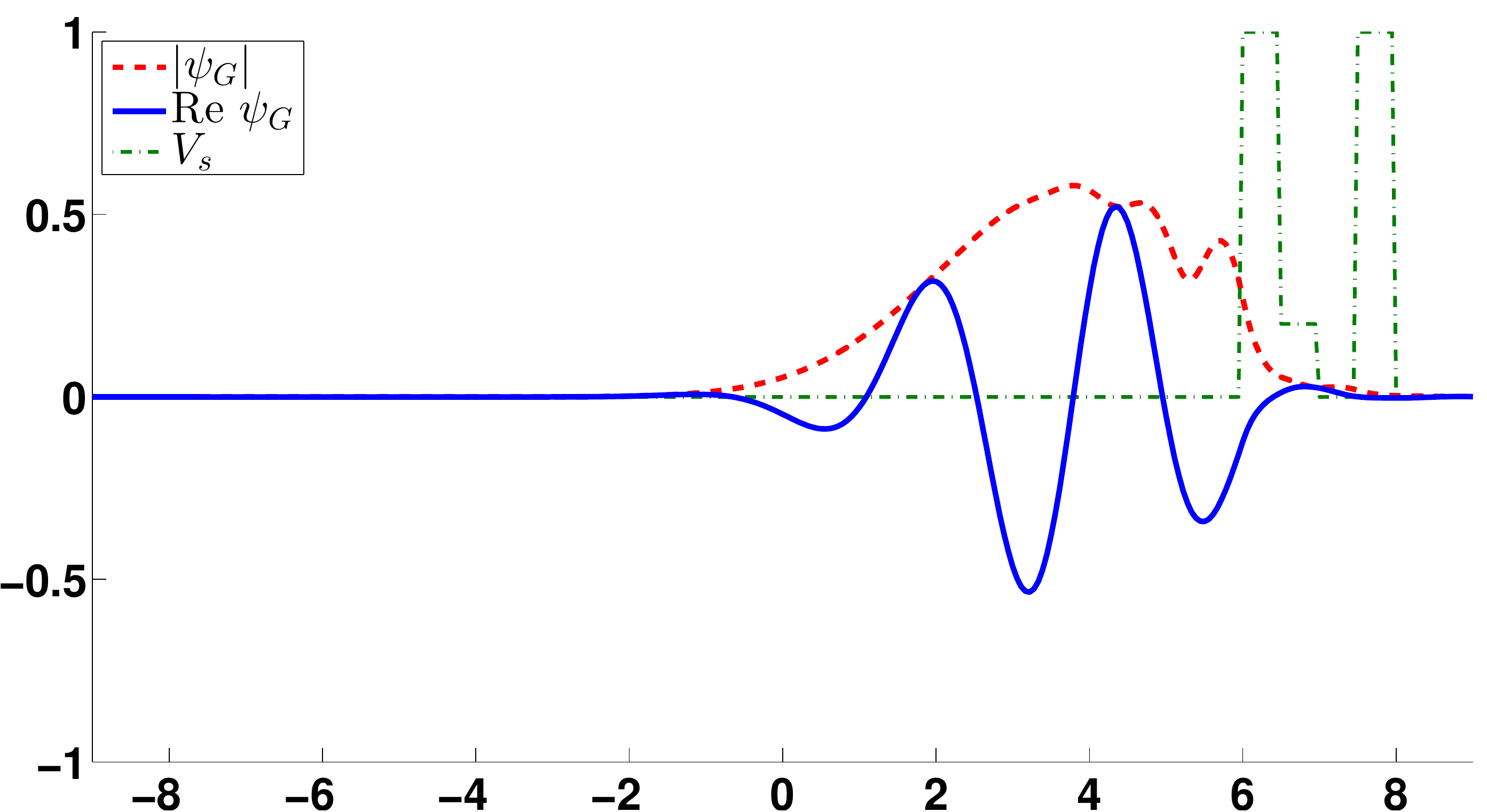}
    } \small{(b) $t_m=1.4286$, $m=45$} \\
    \end{minipage}\hfill
    \begin{minipage}[h]{0.3\linewidth}\center{
        \includegraphics[width=1\linewidth]{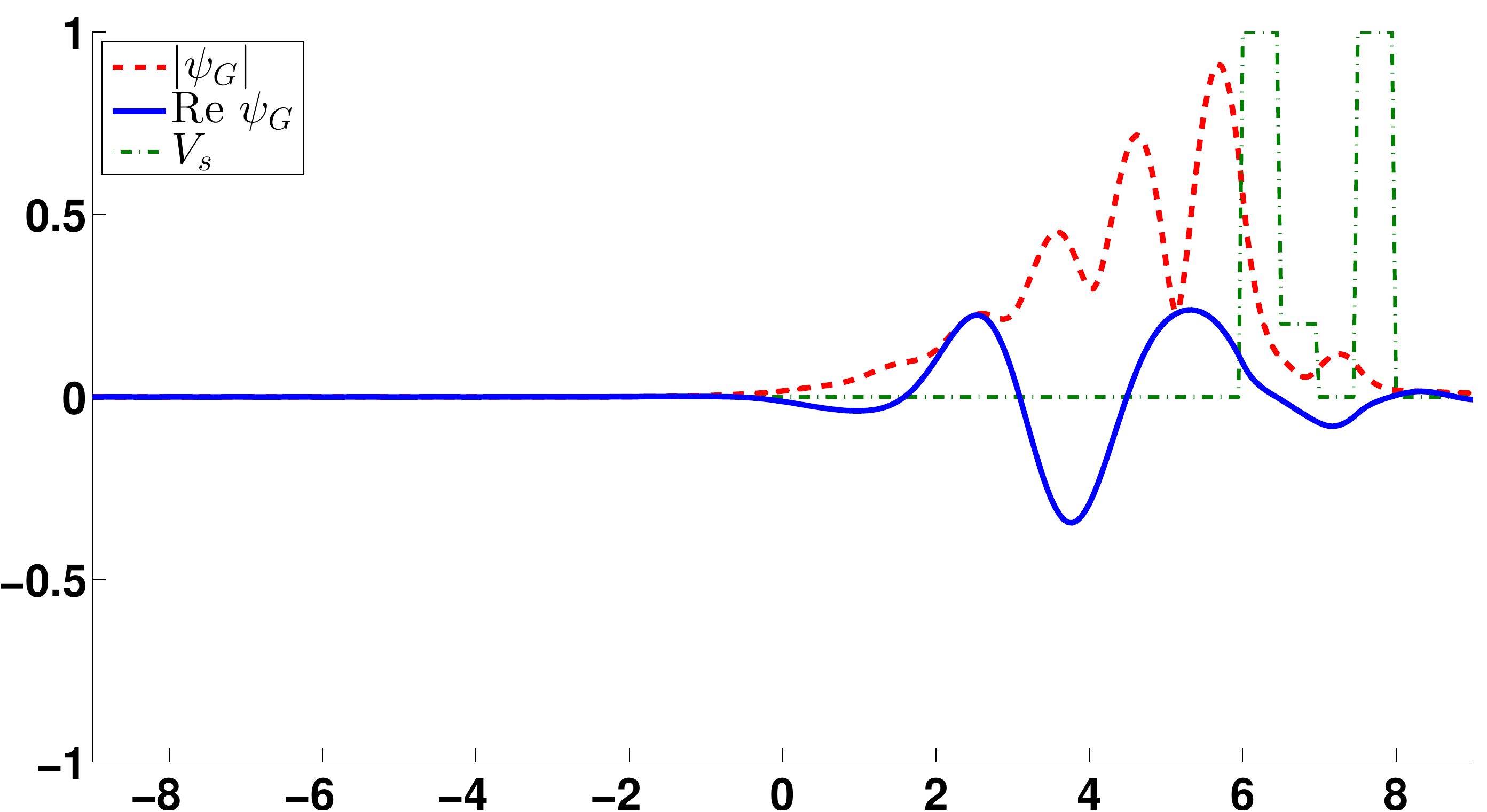}
    } \small{(c) $t_m=2$, $m=63$} \\[2mm]
    \end{minipage}\vfill
    \begin{minipage}[h]{0.3\linewidth}\center{
        \includegraphics[width=1\linewidth]{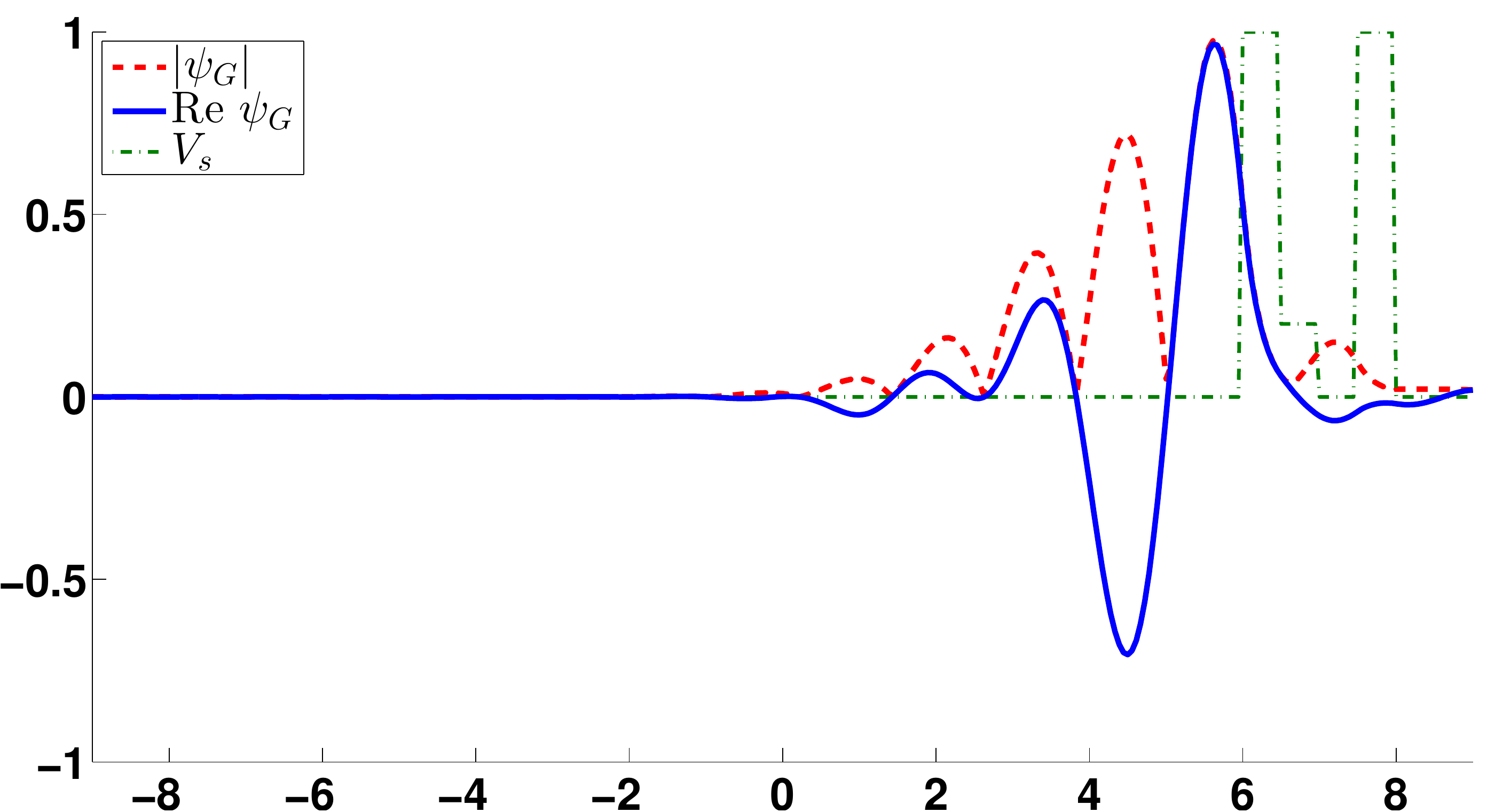}
    } \small{(d) $t_m=2.3810$, $m=75$} \\
    \end{minipage}\hfill
    \begin{minipage}[h]{0.3\linewidth}\center{
        \includegraphics[width=1\linewidth]{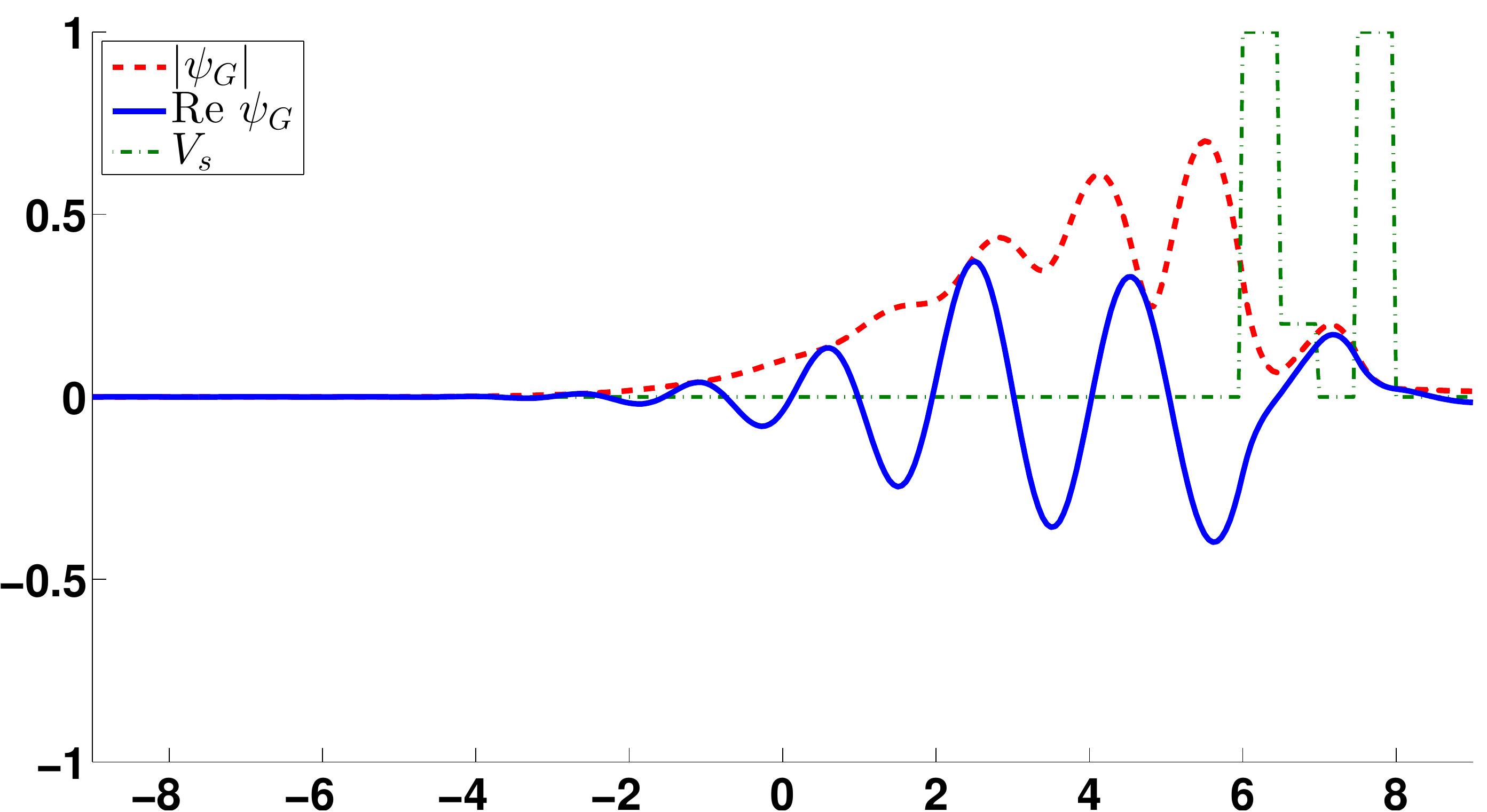}
    } \small{(e) $t_m=3.0476$, $m=96$} \\
    \end{minipage}\hfill
    \begin{minipage}[h]{0.3\linewidth}\center{
        \includegraphics[width=1\linewidth]{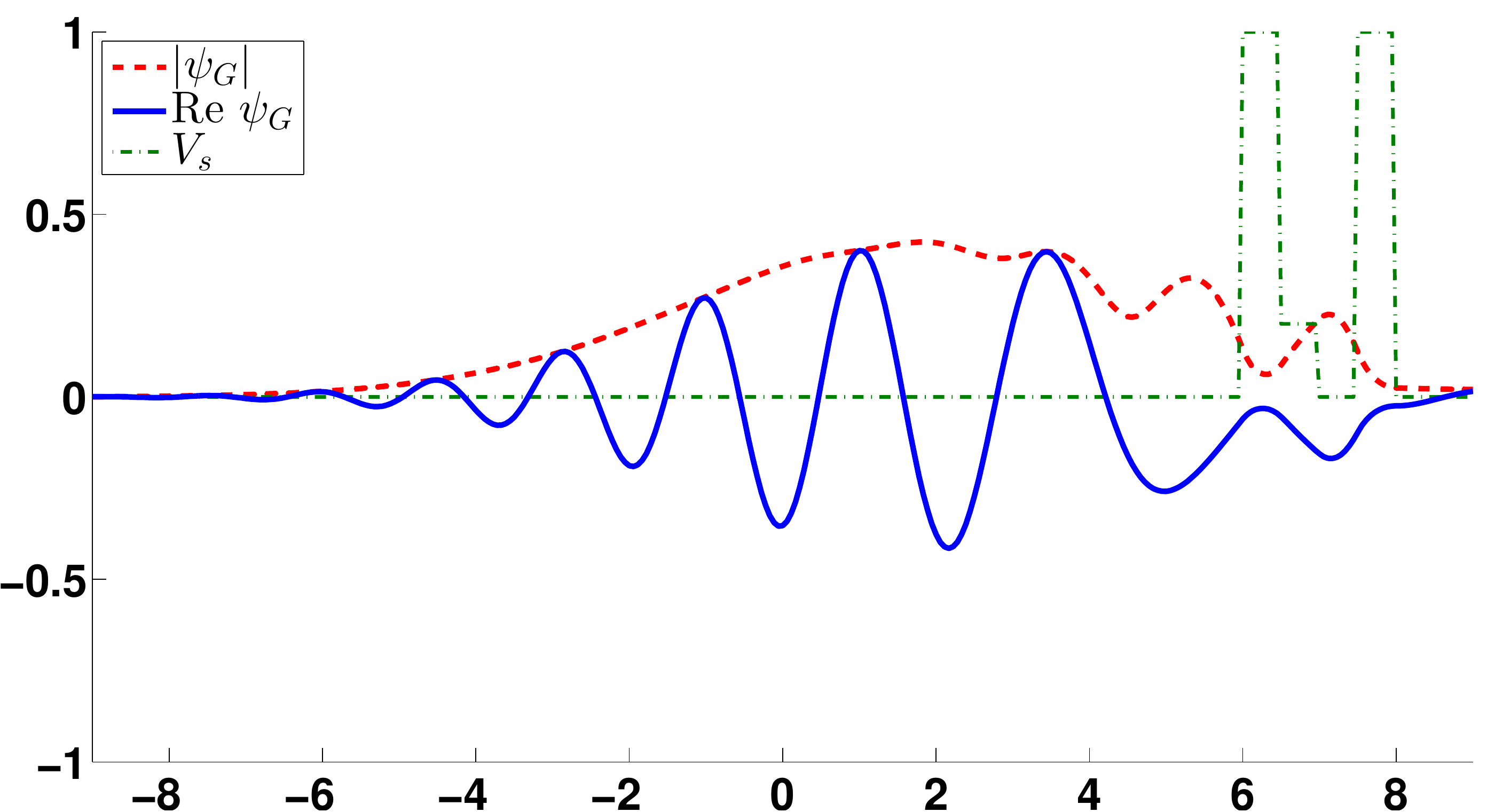}
    } \small{(f) $t_m=4$, $m=126$} \\[2mm]
    \end{minipage}\vfill
    \begin{minipage}[h]{0.3\linewidth}\center{
        \includegraphics[width=1\linewidth]{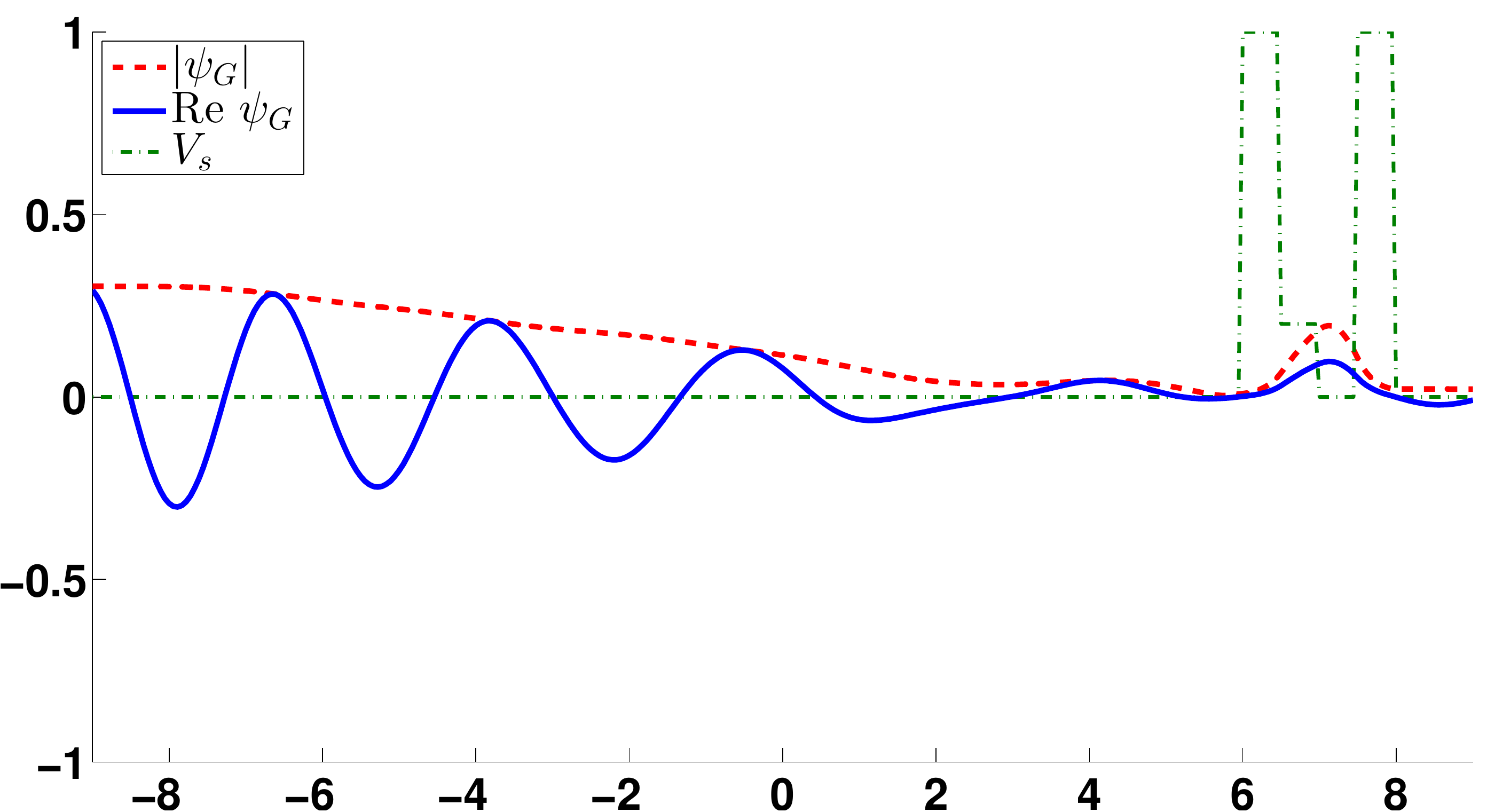}
    } \small{(g) $t_m=8$, $m=252$} \\
    \end{minipage}\hfill
    \begin{minipage}[h]{0.3\linewidth}\center{
        \includegraphics[width=1\linewidth]{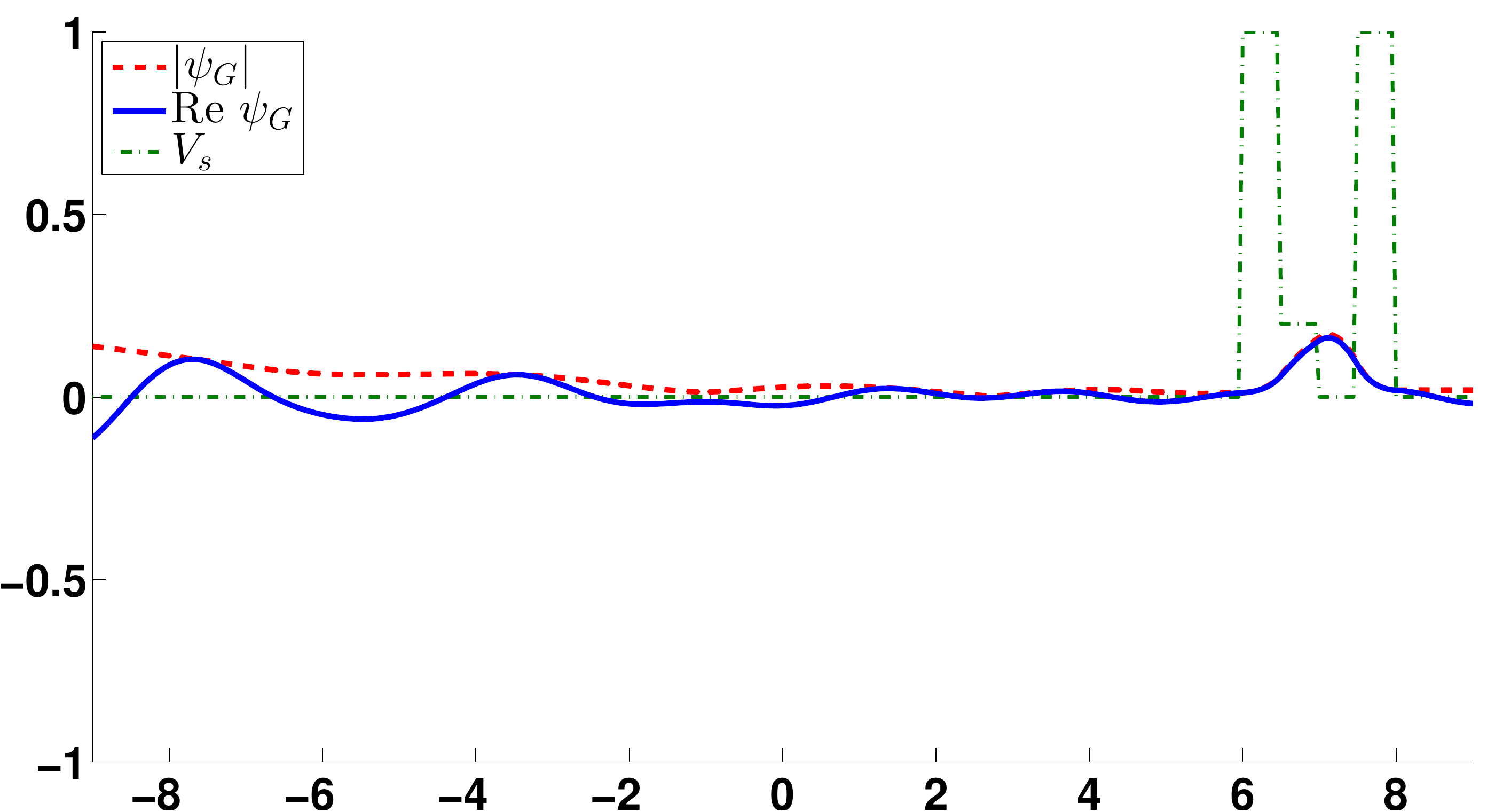}
    } \small{(h) $t_m=12$, $m=378$} \\
    \end{minipage}\hfill
    \begin{minipage}[h]{0.3\linewidth}\center{
        \includegraphics[width=1\linewidth]{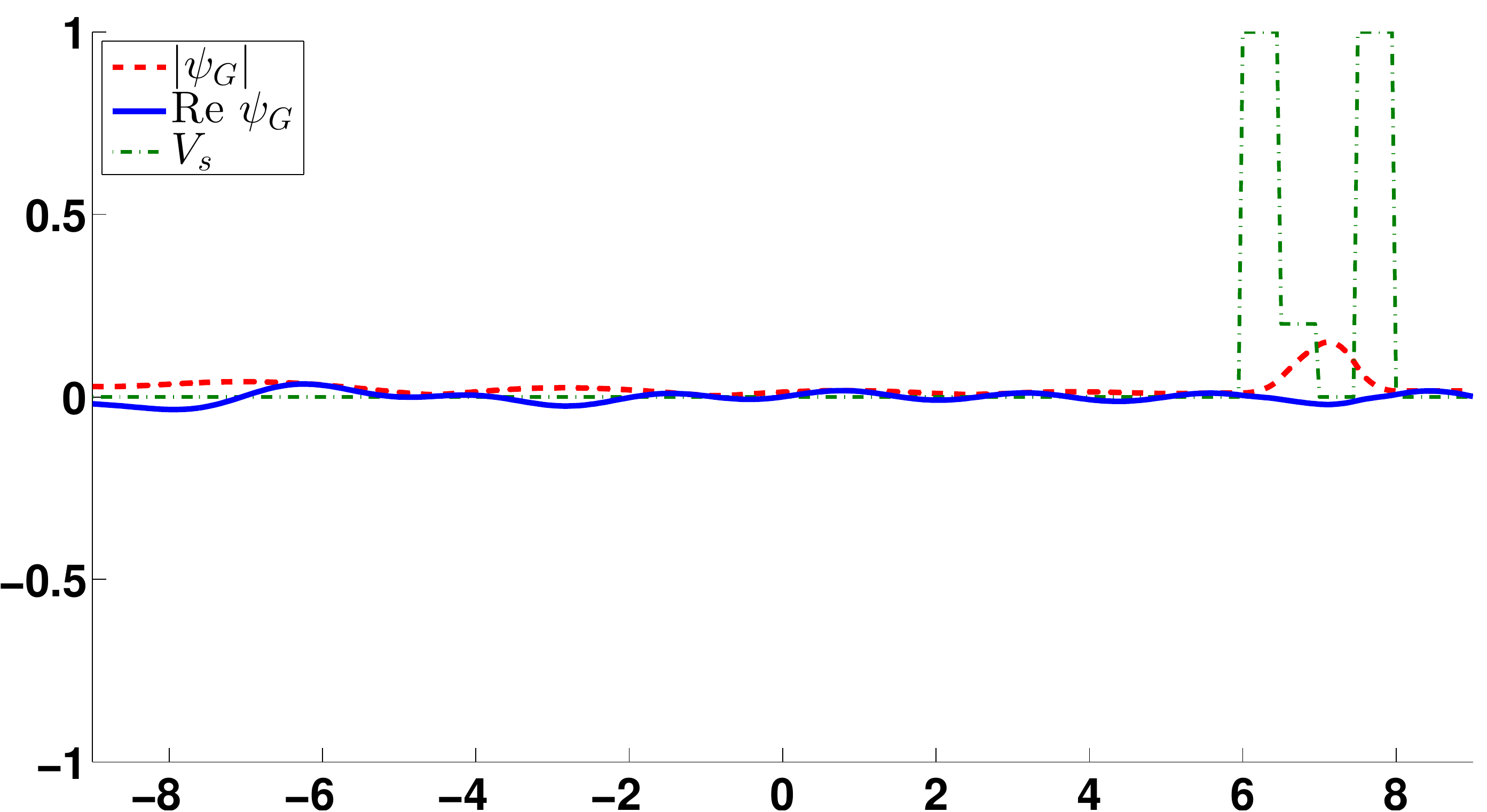}
    } \small{(i) $t_m=16$, $m=504$} \\[2mm]
    \end{minipage}\vfill
\caption{\small{Example 3. $|\Psi_{3R}^{m}|$ and $\Rea\Psi_{3R}^{m}$,
for $n=9$ and $(J,M)=(36,504)$, and the scaled $V$}
\label{fig:EX20:Solution}}
\end{figure}

\textbf{3.3}. In Example 3, we treat the Cauchy problem \eqref{eq:se}, \eqref{eq:ic} for the piecewise constant potential $V=12.5V_s$ with
\begin{gather}
V_s(x)=
\begin{cases}
1 &\text{for}\ \ x\in I_1:=(6,6.5)\cup I_3:=(7.5,8)\\
0.2  &\text{for}\ \ x\in I_2:=(6.5,7)\\
0    &\text{otherwise}
\end{cases}
\label{eq:doubwell}
\end{gather}
and the initial function $\psi^0$ of form \eqref{FEM:s62},
with $x^{(0)}=0$, $k=\sqrt{7}$ and $\alpha=1$ now, thus the tunneling through the double barrier stepped quantum well is studied.
Recall that $B(x)\equiv 2$.
$|\psi^0|$ and $\Rea\psi^0$ and the scaled potential $V$ are given on Fig.~\ref{fig:EX20:Solution}(a).
Also $X=9$ and $T=16$ are taken. This is the most complicated example from the review \cite{AABES08}.
\par By scaling of the coordinates this example could be transformed to more close to Examples 1 and 2.
Namely, equivalently we could consider the same Schr\"odinger equation with $B(x)\equiv 1$ but $V(x)=2025V_s(\frac{x}{9})$ as well as $\psi^0$ of form \eqref{FEM:s62} divided by 3, with $x^{(0)}=0$, $k=9\sqrt{7}$ and $\alpha=\frac{1}{81}$, for $X=1$ and $T=\frac{8}{81}$.
\par Now $|\psi^0(x)|< 1.01E$$-9$ outside $\Omega$. Once again we pose the discrete TBCs at the both artificial boundaries $x=\pm X$. For $J=36$, each segment $\bar{I}_1$, $\bar{I}_2$ and $\bar{I}_3$ in \eqref{eq:doubwell} consists of exactly one element so that we take our $J$ as multiples of 36.
Notice that it is important to treat the discontinuity points of $V$
carefully, and even this allows to diminish significantly errors in finite-difference computations among all in \cite{AABES08} (more details are given in \cite{IZ13}).
\par We consider as the pseudo-exact solution $\Psi_{4R}^{4m}$ computed for the high $n=9$,  $J=144$ and rather large $M=8064$ (this choice is justified on Fig.~ \ref{fig:EX20r:Error} below).
\par First the wave moves toward the barrier; after the interaction with it, the main piece of the wave is reflected and moves in the opposite direction whereas the small piece remains trapped and oscillating inside the well, and another very small piece passes through the barrier and moves to the right.
The solution is represented by $\Psi_{3R}$ computed for $n=9$ and only $(J,M)=(36,504)$ (that is enough according to Table \ref{tab:EX20r} below).
\par
\begin{figure}[htbp]
    \begin{minipage}[h]{0.49\linewidth}\center{
        \includegraphics[width=1\linewidth]{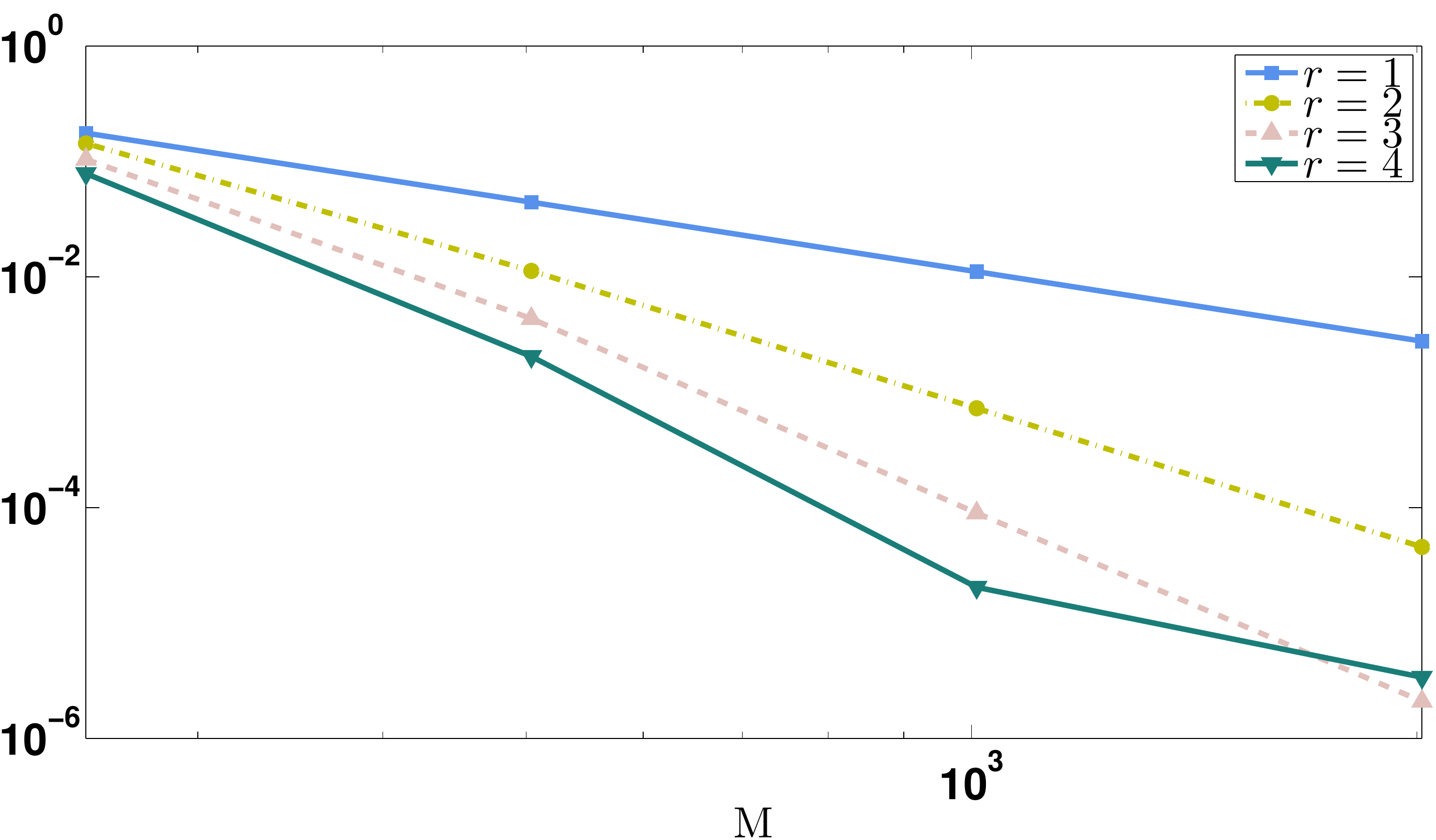}
    } \small{(a) in $L^2$ space norm} \\
    \end{minipage}\hfill
    \begin{minipage}[h]{0.49\linewidth}\center{
        \includegraphics[width=1\linewidth]{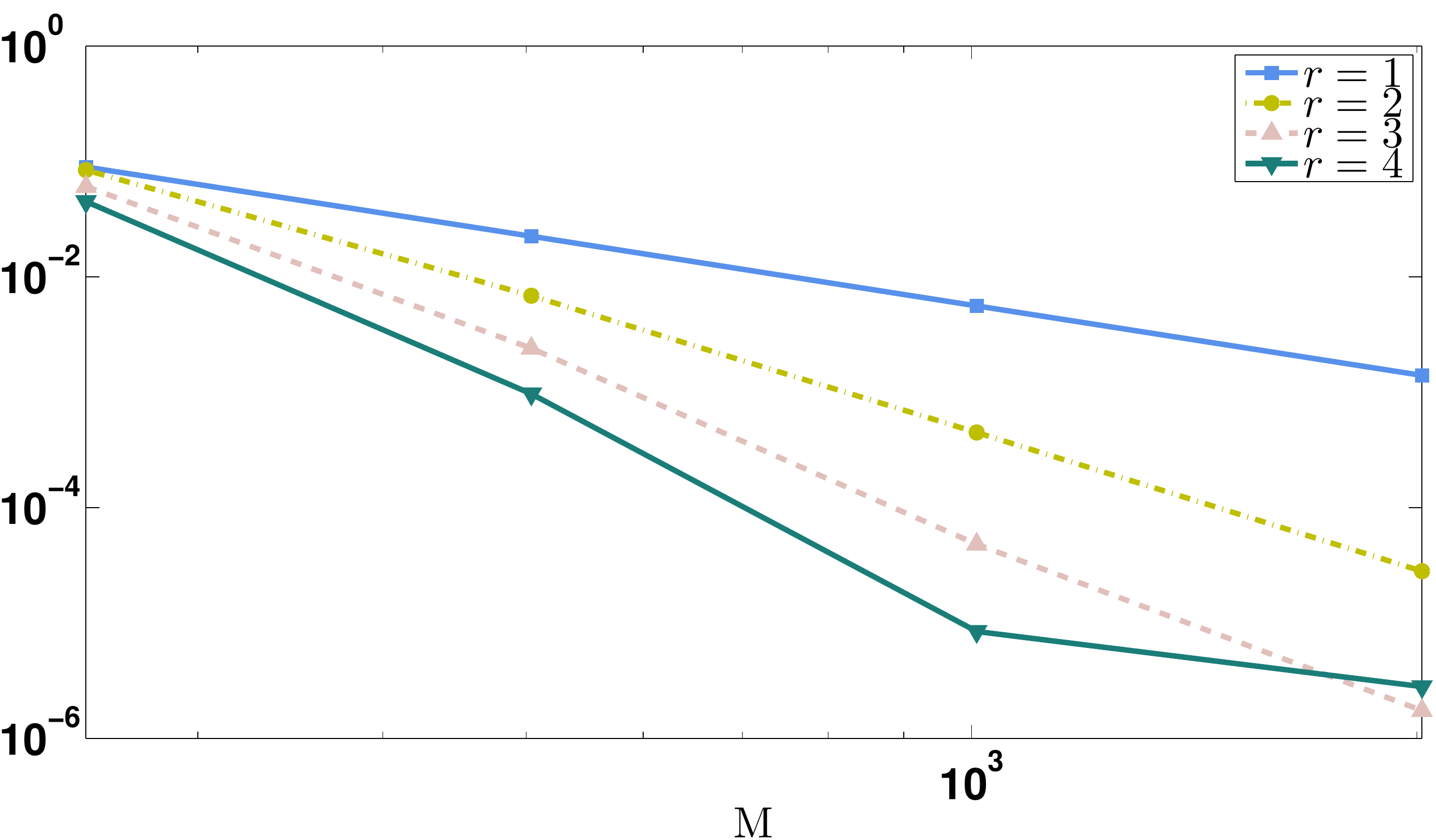}
    } \small{(b) in $C$ space norm} \\
    \end{minipage}
\caption{\small{Example 3. The errors $\max_{0\leq rm\leq M}\|\psi^{rm}-\Psi_{rR}^{rm}\|$,
for $n=9$ and $J=36$, in dependence with $r=1,2,3,4$ and $M=252\cdot 2^q$, $q=0,1,2,3$}
\label{fig:EX20r:N=9:J=36:MaxAbsError}}
\end{figure}

\begin{figure}[htbp]
    \begin{minipage}[h]{0.49\linewidth}\center{
        \includegraphics[width=1\linewidth]{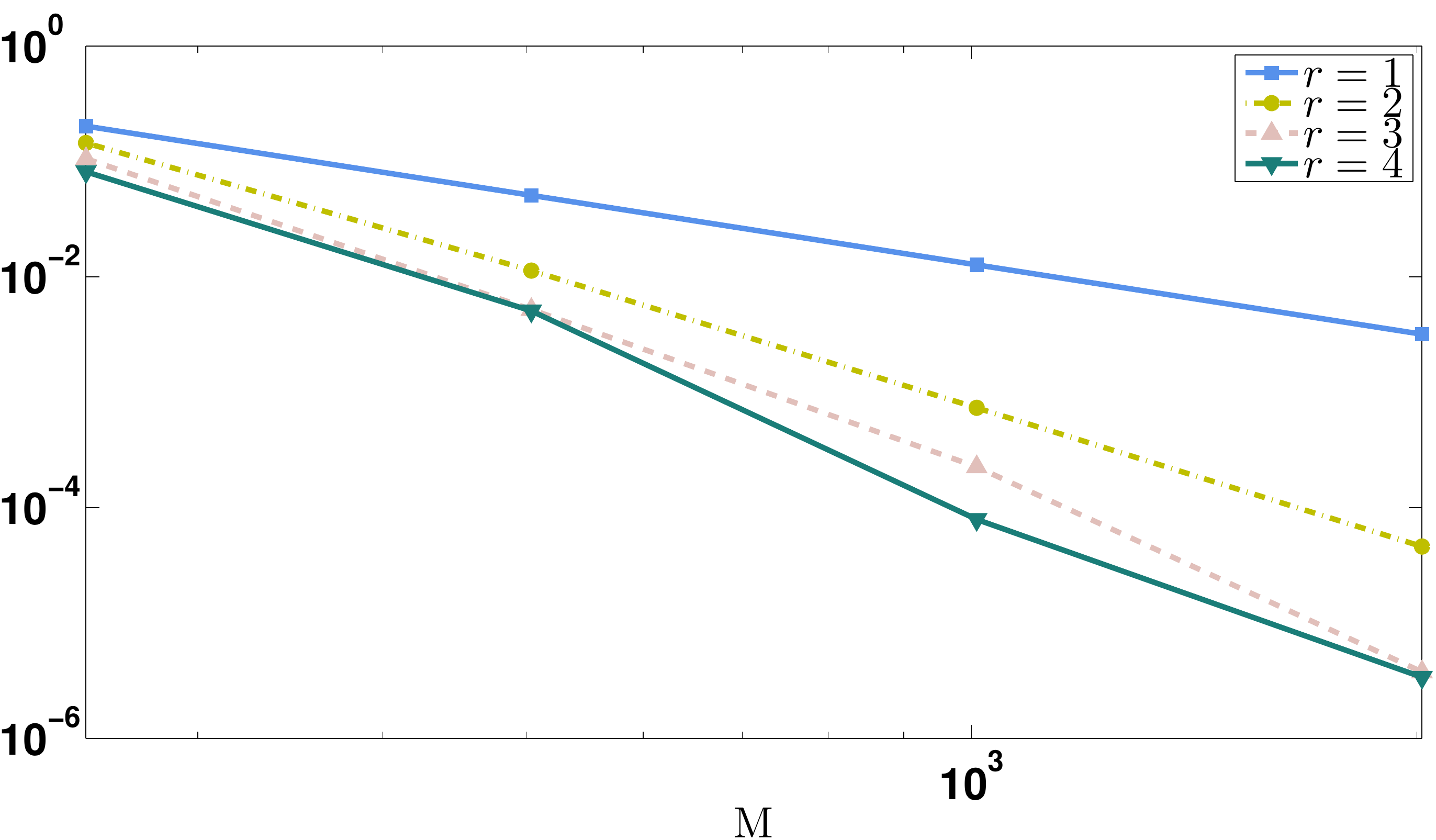}
    } \small{(a) in $L^2$ space norm} \\
    \end{minipage}\hfill
    \begin{minipage}[h]{0.49\linewidth}\center{
        \includegraphics[width=1\linewidth]{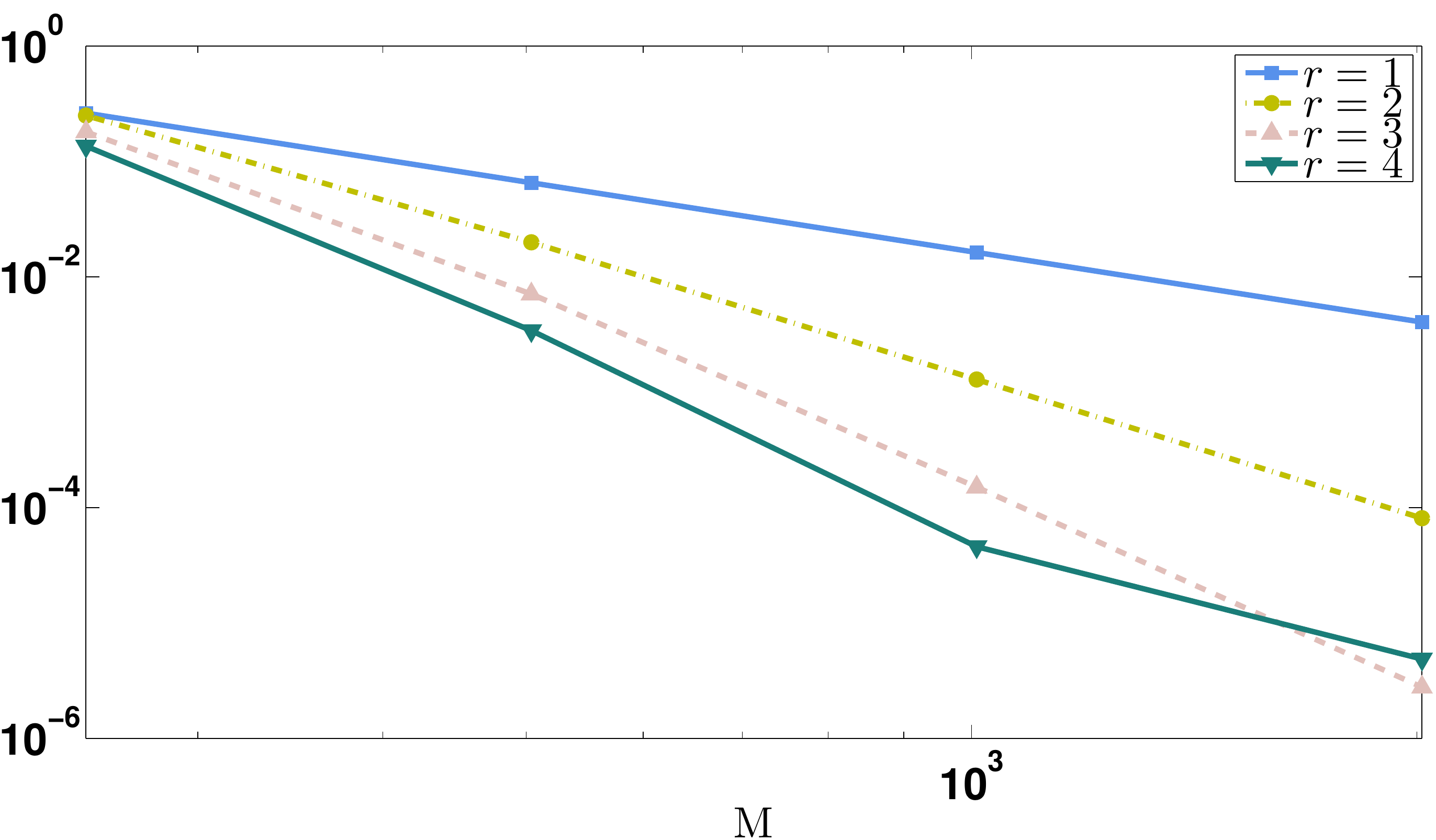}
    } \small{(b) in $C$ space norm} \\
    \end{minipage}
\caption{\small{Example 3. The relative errors
$\max_{0\leq rm\leq M}\|\psi^{rm}-\Psi_{rR}^{rm}\|/\|\psi^{rm}\|$,
for $n=9$ and $J=36$, in dependence with $r=1,2,3,4$ and $M=252\cdot 2^q$, $q=0,1,2,3$}
\label{fig:EX20r:N=9:J=36:MaxRelError}}
\end{figure}

 On Fig. \ref{fig:EX20r:N=9:J=36:MaxAbsError}, the errors
$E_r^M=\max_{0\leq rm\leq M}\|\psi^{rm}-\Psi_{rR}^{rm}\|$
are given, for $n=9$ and $J=36$, in dependence with
$r$ and  $M=M_q=252\cdot 2^q$, $q=0,1,2,3$.
There are some differences in the behavior of the corresponding relative errors
$\max_{0\leq rm\leq M}\|\psi^{rm}-\Psi_{rR}^{rm}\|/\|\psi^{rm}\|$, see
Fig. \ref{fig:EX20r:N=9:J=36:MaxRelError}, but they are not crucial.

\par In Table \ref{tab:EX20r}, the same data as on Fig. \ref{fig:EX20r:N=9:J=36:MaxAbsError}
together with their ratios $E_r^{M_q}/E_r^{M_{q-1}}$ are presented. Now the ratios are close to $2^{2r}$ only for $r=1$ as well as  $r=2$ and $q\geq2$. For $r=3$ and 4, they are less than $2^{2r}$ but still grow rapidly as $r=1,2,3,4$ increases (for fixed $q$) excepting the case $r=4$ and the last $q=3$.
\par Comparing the results to those in \cite{AABES08}, we obtain \textit{much more} accurate results using \textit{much less} amount of both elements and time steps. In particular, we achieve the relative error
$e=\max_{0\leq 3m\leq M}\|\psi^{3m}-\Psi_{3R}^{3m}\|_{L_h^2}/\|\psi^{3m}\|_{L_h^2}\approx 3.77E$$-6$ for $n=9$ using only $(J,M)=(36,2016)$ versus the best presented there $e\geq 4E$$-4$ using $(J,M)=(6000,16000)$.

\par
\begin{figure}[htbp]
    \begin{minipage}[h]{0.49\linewidth}\center{
        \includegraphics[width=1\linewidth]{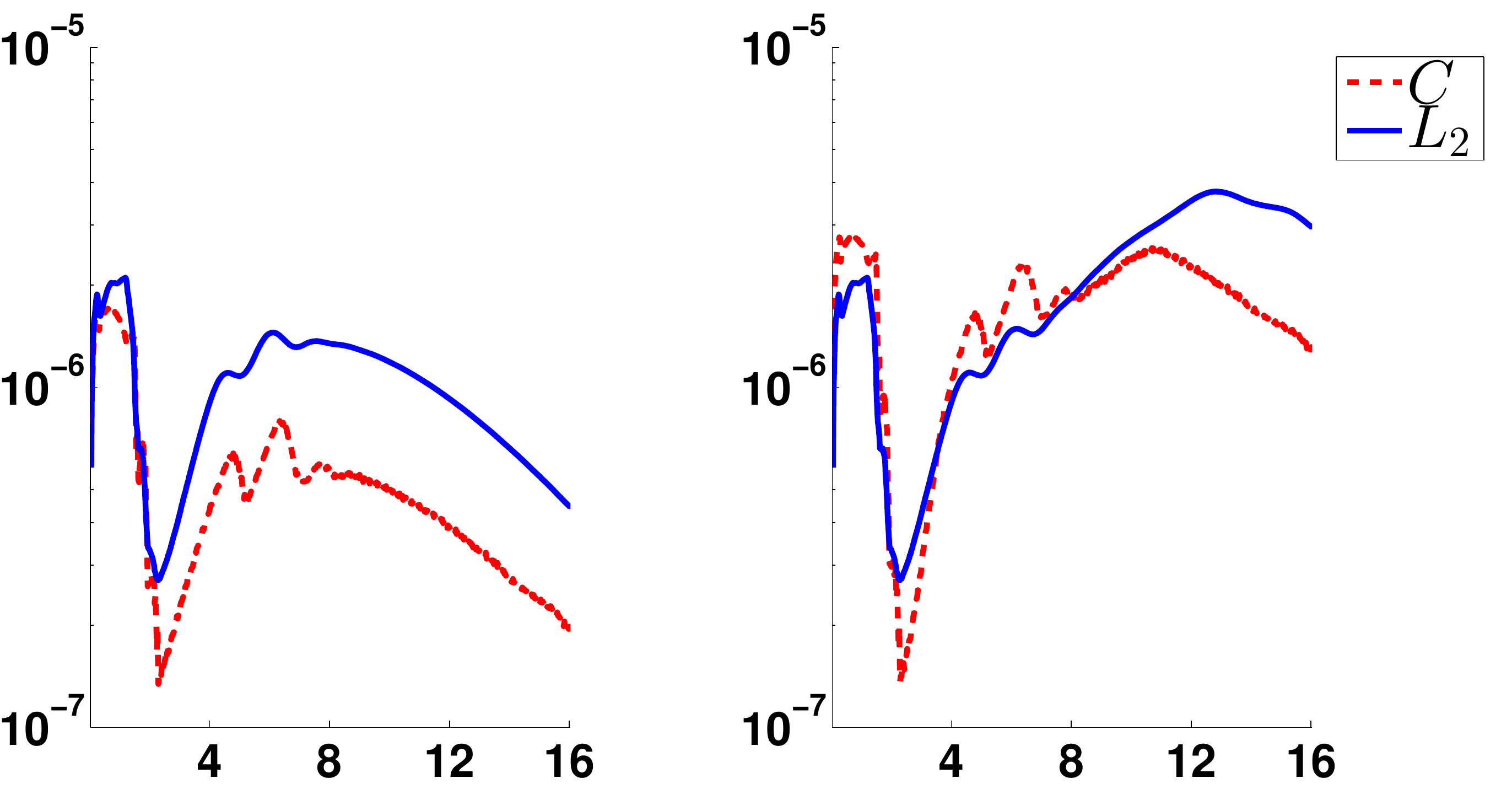}
    } \small{(a) $r=3$ and $M=2016$} \\
    \end{minipage}\hfill
    \begin{minipage}[h]{0.49\linewidth}\center{
        \includegraphics[width=1\linewidth]{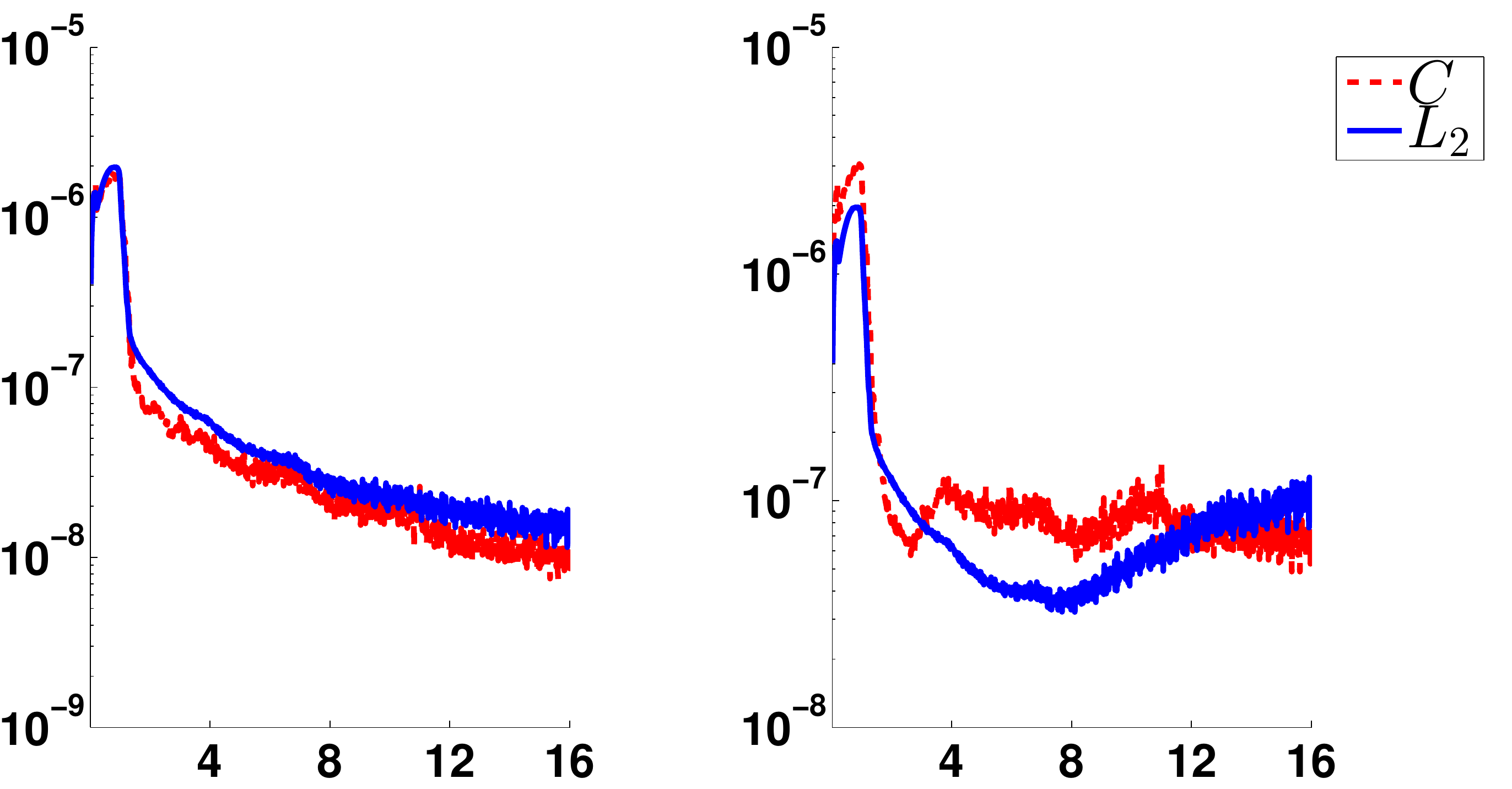}
    } \small{(b) $r=4$ and $M=4032$} \\
     \end{minipage}
\caption{\small{Example 3. The absolute $\|\psi^{rm}-\Psi_{rR}^{rm}\|$ (left) and relative
$\|\psi^{rm}-\Psi_{rR}^{rm}\|/\|\psi^{rm}\|$ (right) errors, for $n=9$ and $J=72$, in $L^2$ and $C$ norms:
(a) for $r=3$ and $M=2016$, and
(b) for $r=3$ and $M=4032$,
both in dependence with $t_m$}
\label{fig:EX20r:Error3_4}}
\end{figure}

 On the other hand, on Fig. \ref{fig:EX20r:N=9:J=36:MaxAbsError} and in Table \ref{tab:EX20r} we see the degradation of the error behavior at the level about $1E$$-6$ for $r=4$.
Fig. \ref{fig:EX20r:Error3_4}  contributes to that, for $r=3$ and $r=4$, showing the absolute and relative errors:
(a) $\|\psi^{3m}-\Psi_{3R}^{3m}\|$ and $\|\psi^{3m}-\Psi_{3R}^{3m}\|/\|\psi^{3m}\|$ for $M=2016$, and
(b) $\|\psi^{4m}-\Psi_{4R}^{4m}\|$ and $\|\psi^{4m}-\Psi_{4R}^{4m}\|/\|\psi^{4m}\|$ for $M=4032$,
both for $n=9$ and $J=72$, in dependence with $t_m$. We see that the maximal absolute and relative errors (except for the $L^2$ relative one in the case (a)) occur near $t=0$ \textit{before} the active interaction of the wave with the potential.
\par On Fig. \ref{fig:EX20r:Error}, we present the changes in the above pseudo-exact solution due to 4 times increasing $M$ up to 32256 or $J$ up to 576 are less than respectively $1E$$-6$ and $4E$$-8$ in the uniform in time and both $L^2$ and $C$ space norms (the former one is also less than $1E$$-8$ on the right time half-segment $\bigl[\frac{T}{2},T\bigr]$).
So the data in Table \ref{tab:EX20r} are correct but even the significant increasing $M$ does not improve the error essentially, and the maximal error is located near $t=0$.  This also confirms the degradation.
(Notice that the similar degradation could be seen in Example 2 too for values of $M$ larger than on Fig. \ref{fig:EX22r:N=9:J=60:MaxAbsError} and in Table \ref{tab:EX22r}; moreover, this appears at a higher error level if the potential is situated closer to $x^{(0)}$).
It seems that this is due to non-smoothness of the potential and invalidity of the error expansion \eqref{eq:errexp} for very small $\tau$ (for larger $\tau$, the smallness of the initial function under the potential support prevents the effect).

\begin{table}\centering{
\begin{tabular}{rrcrcrcrc}
\multicolumn {1}{c}{$M$}&
\multicolumn {2}{c}{$r=1$}&
\multicolumn {2}{c}{$r=2$}&
\multicolumn {2}{c}{$r=3$}&
\multicolumn {2}{c}{$r=4$}\\
\midrule %
$252$&$0.18$&--&$0.14$&--&$0.1$&--&$7.85E{-}2$&--\\%
$504$&$4.43E{-}2$&$3.97$&$1.12E{-}2$&$12.72$&$4.33E{-}3$&$24.02$&$2.04E{-}3$&$38.46$\\%
$1\,008$&$1.11E{-}2$&$4$&$7.25E{-}4$&$15.49$&$8.97E{-}5$&$48.32$&$2.04E{-}5$&$100.05$\\%
$2\,016$&$2.77E{-}3$&$4$&$4.56E{-}5$&$15.92$&$2.09E{-}6$&$42.83$&$3.38E{-}6$&$6.03$\\
\end{tabular}%
\\[2mm]
\small{(a) in $L^2$ space norm}\\[2mm]
\begin{tabular}{rrcrcrcrc}
\multicolumn {1}{c}{$M$}&
\multicolumn {2}{c}{$r=1$}&
\multicolumn {2}{c}{$r=2$}&
\multicolumn {2}{c}{$r=3$}&
\multicolumn {2}{c}{$r=4$}\\
\midrule %
$252$&$8.93E{-}2$&--&$8.41E{-}2$&--&$6.04E{-}2$&--&$4.49E{-}2$&--\\%
$504$&$2.24E{-}2$&$3.98$&$6.84E{-}3$&$12.3$&$2.40E{-}3$&$25.2$&$9.67E{-}4$&$46.42$\\%
$1\,008$&$5.60E{-}3$&$4.01$&$4.47E{-}4$&$15.31$&$4.81E{-}5$&$49.82$&$8.42E{-}6$&$114.8$\\%
$2\,016$&$1.40E{-}3$&$4$&$2.81E{-}5$&$15.91$&$1.74E{-}6$&$27.62$&$2.81E{-}6$&$3$\\
\end{tabular}%
\\[2mm]
\small{(b) in $C$ space norm}}\\[2mm]
\caption{\small Example 3. The errors $E_r^M=\max_{0\leq rm\leq M}\|\psi^{rm}-\Psi_{rR}^{rm}\|$
and their ratios $E_r^{M_q}/E_r^{M_{q-1}}$, for $n=9$ and $J=36$, in dependence with
$r$ and $M_q=252\cdot 2^q$, $q=0,1,2,3$}
\label{tab:EX20r}
\end{table}

\par
\begin{figure}[htbp]
    \begin{minipage}[h]{0.49\linewidth}\center{
        \includegraphics[width=1\linewidth]{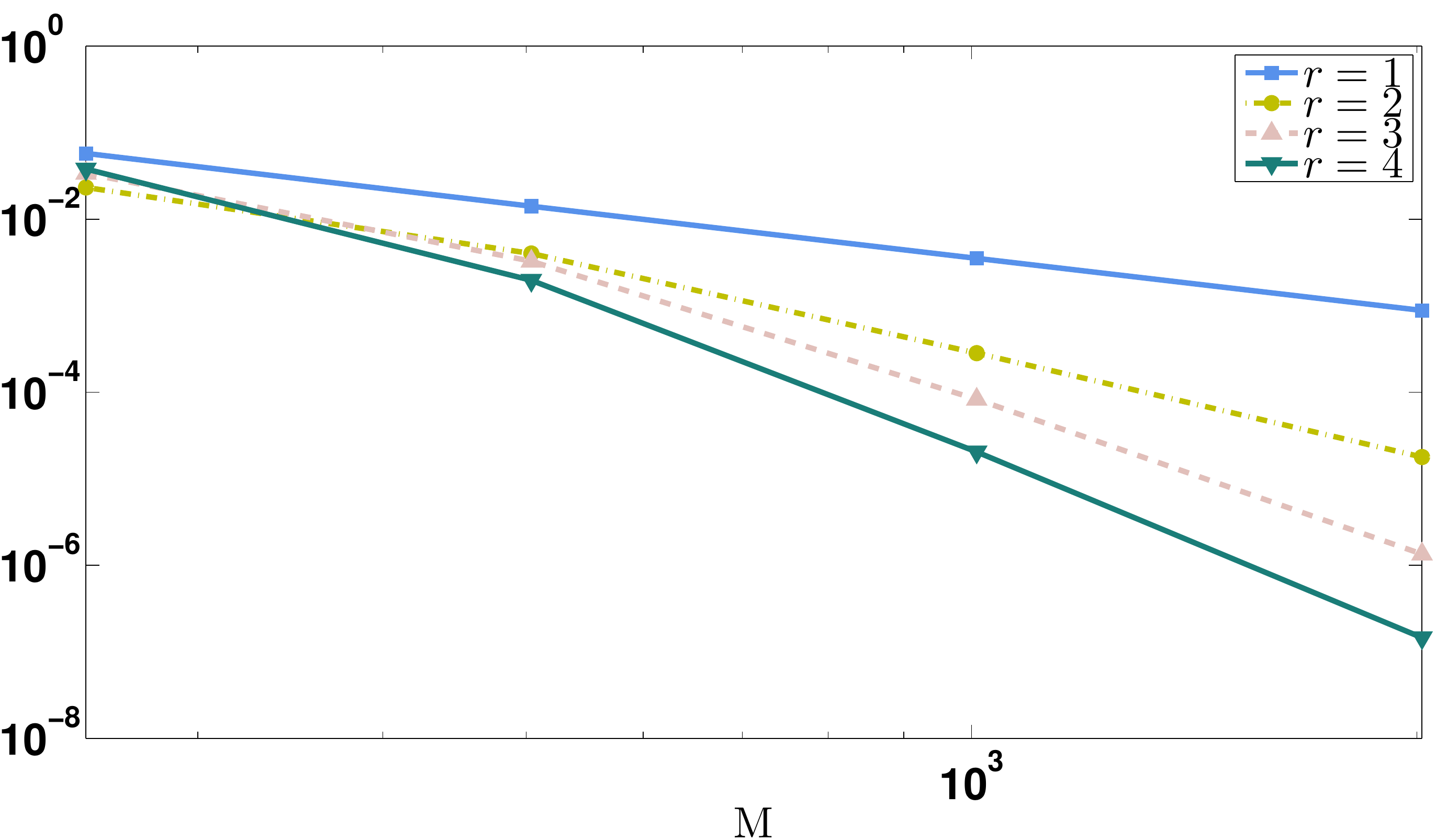}
    } \small{(a) in $L^2$ space norm} \\
    \end{minipage}\hfill
    \begin{minipage}[h]{0.49\linewidth}\center{
        \includegraphics[width=1\linewidth]{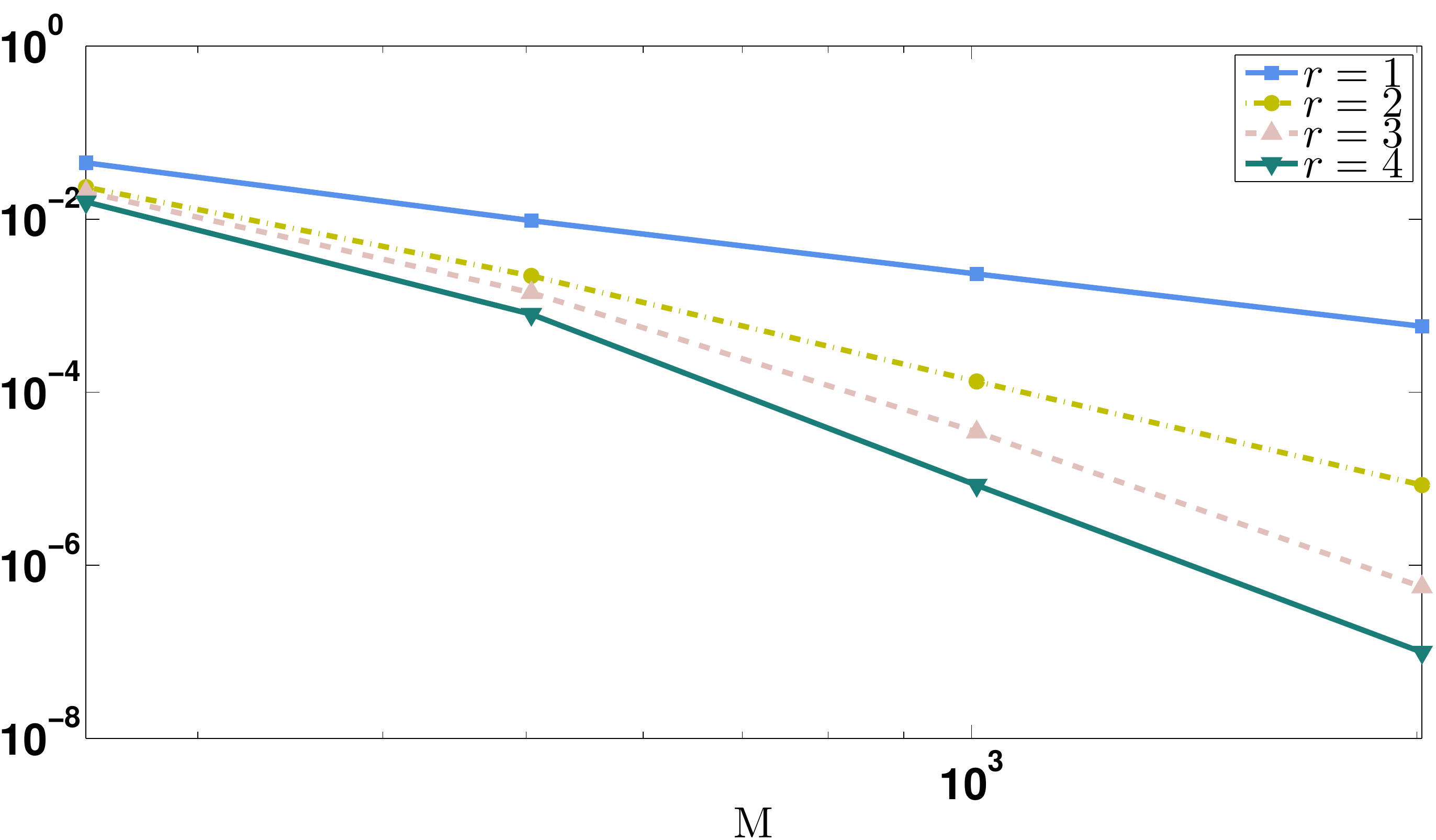}
    } \small{(b) in $C$ space norm} \\
    \end{minipage}
\caption{\small{Example 3. The errors $\max_{M/2\leq rm\leq M}\|\psi^{rm}-\Psi_{rR}^{rm}\|$,
for $n=9$ and $J=36$, in dependence with $r=1,2,3,4$ and $M=252\cdot 2^q$, $q=0,1,2,3$}
\label{fig:EX20r:N=9:J=36:MaxAbsError:M-2}}
\end{figure}

\begin{figure}[htbp]
    \begin{minipage}[h]{0.49\linewidth}\center{
        \includegraphics[width=1\linewidth]{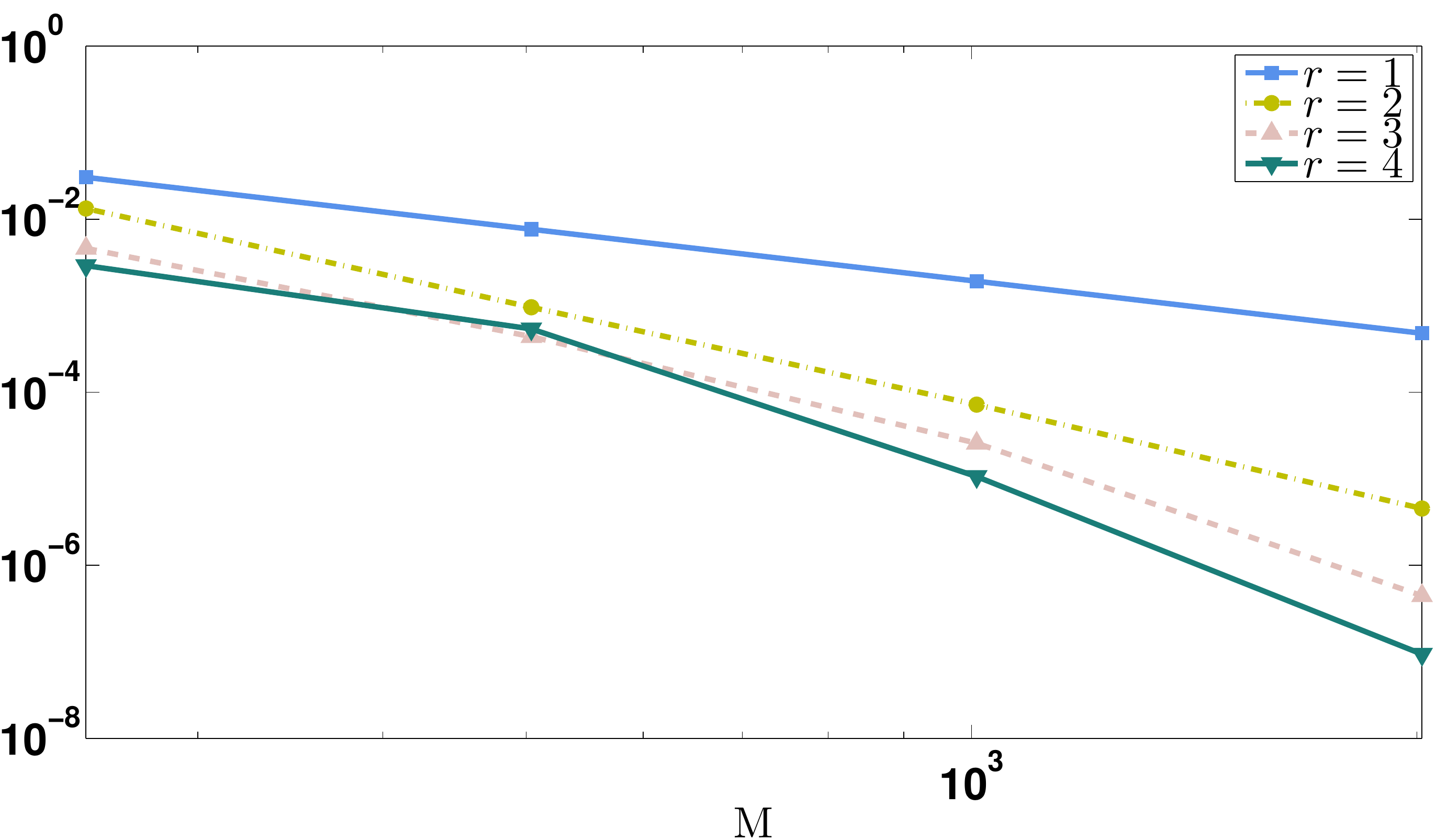}
    } \small{(a) in $L^2$ space norm} \\
    \end{minipage}\hfill
    \begin{minipage}[h]{0.49\linewidth}\center{
        \includegraphics[width=1\linewidth]{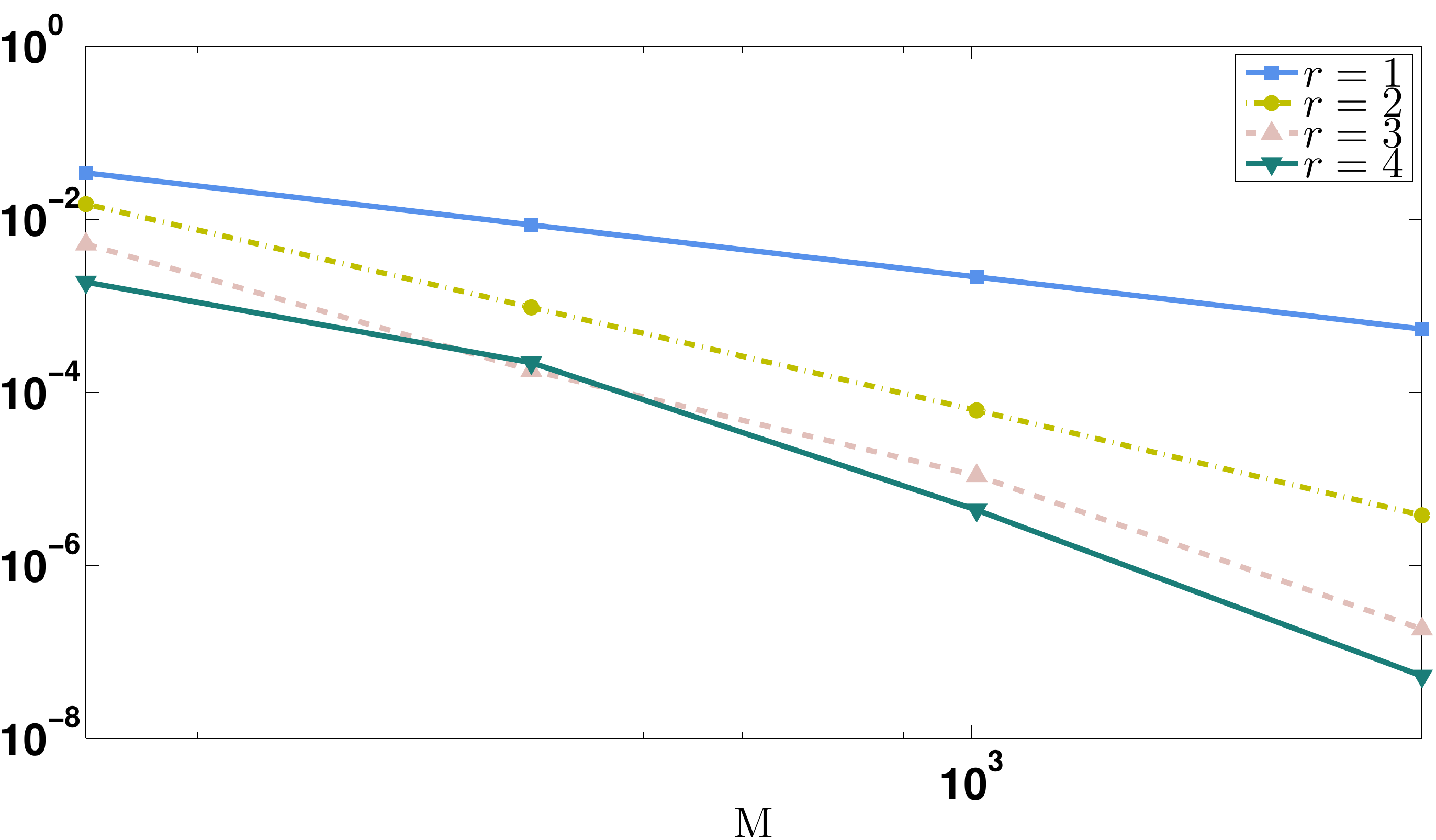}
    } \small{(b) in $C$ space norm} \\
    \end{minipage}
\caption{\small{Example 3. The errors $\|\psi^{rM}-\Psi_{rR}^{rM}\|$,
for $n=9$ and $J=36$, in dependence with $r=1,2,3,4$ and $M=252\cdot 2^q$, $q=0,1,2,3$}
\label{fig:EX20r:N=9:J=36:MaxAbsError:M}}
\end{figure}

 The error behavior is not the same during the whole time segment $[0,T]$.
On Fig. \ref{fig:EX20r:N=9:J=36:MaxAbsError:M-2} the behavior of the similar errors as on Fig. \ref{fig:EX20r:N=9:J=36:MaxAbsError} and in Table \ref{tab:EX20r} but on the right half-segment $\bigl[\frac{T}{2},T\bigr]$ only (after the interaction of the wave with the potential) is shown and is clearly better and without the degradation.
In addition, on Fig. \ref{fig:EX20r:N=9:J=36:MaxAbsError:M} the corresponding errors $\|\psi^{rM}-\Psi_{rR}^{rM}\|$ at the final time $t_M=T$ are demonstrated. For $r=2,3,4$, the former and latter graphs differ more for smaller $M$ but become closer for larger $M$.
Note that only the similar relative $L^2$ errors at the final time are contained in \cite{AABES08}.
But there, for the FEM with $n=1$ and 2, the semi-discrete TBCs had been used (since the corresponding discrete TBCs had been unknown and were designed later in \cite{DZZ09,ZZ12}) which are unable to ensure so nice behavior for the relative errors as for the absolute ones in general, see also \cite{ZZ14a}.

\par
\begin{figure}[htbp]
    \begin{minipage}[h]{0.49\linewidth}\center{
        \includegraphics[width=1\linewidth]{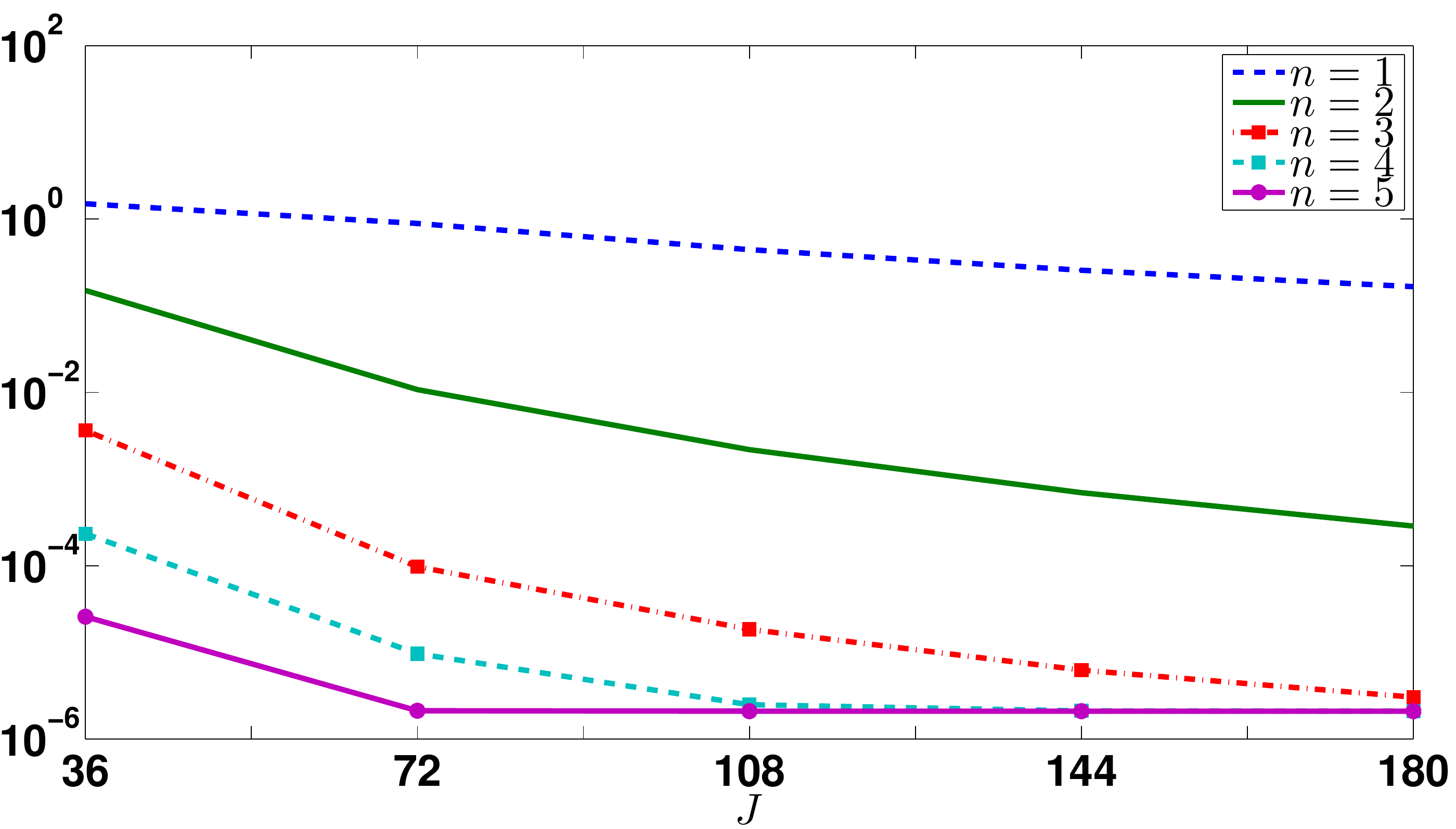}
    } \small{(a) in $L^2$ space norm} \\
    \end{minipage}\hfill
    \begin{minipage}[h]{0.49\linewidth}\center{
        \includegraphics[width=1\linewidth]{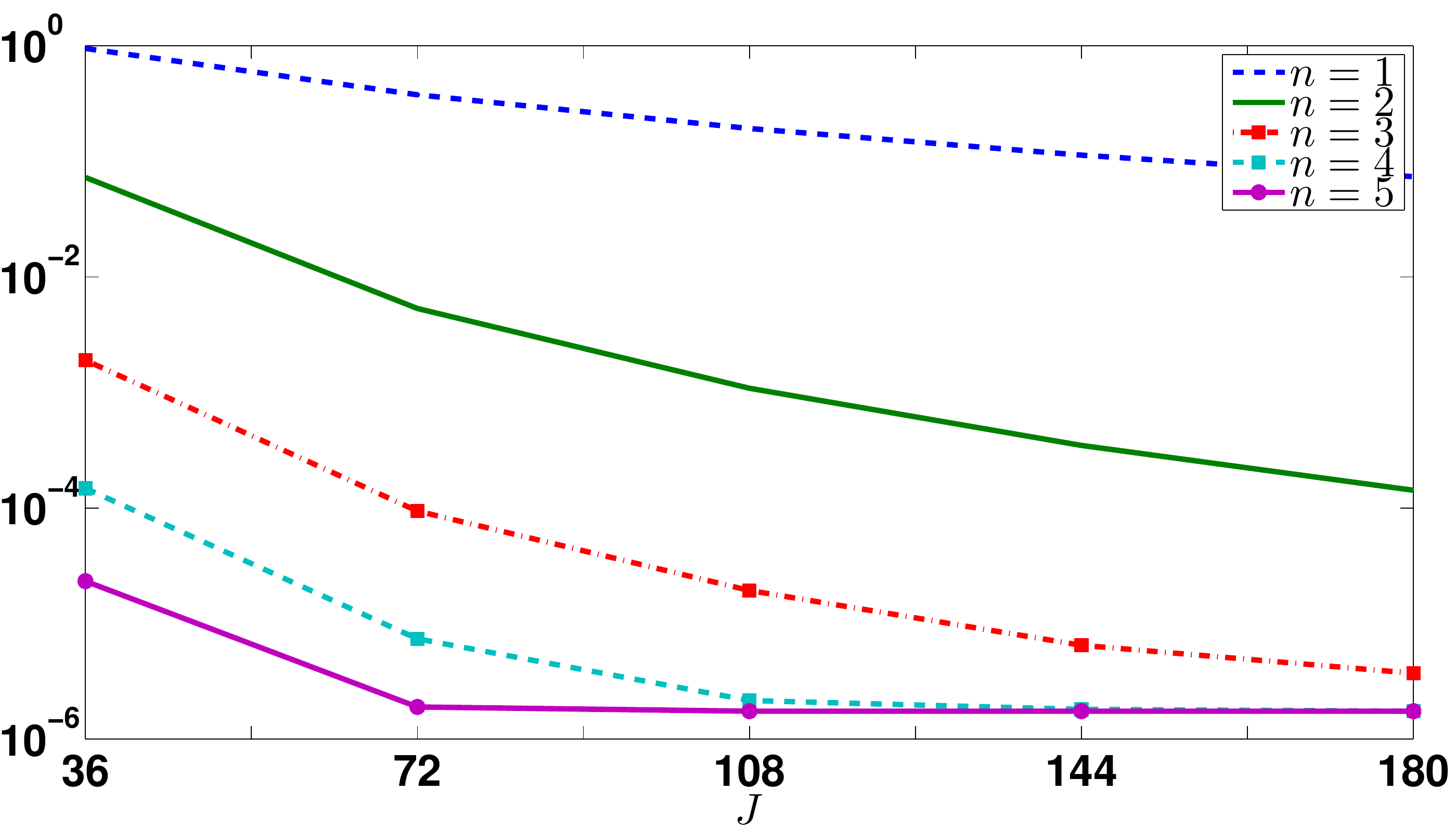}
    } \small{(b) in $C$ space norm} \\
    \end{minipage}
\caption{\small{Example 3. The errors $\max_{0\leq 3m\leq M}\|\psi^{3m}-\Psi_{3R}^{3m}\|$,
for $M=2016$, in dependence with $n$ and $J$}
\label{EX20r:R=3:MaxAbsError}}
\end{figure}

\begin{figure}[htbp]
    \begin{minipage}[h]{0.49\linewidth}\center{
        \includegraphics[width=1\linewidth]{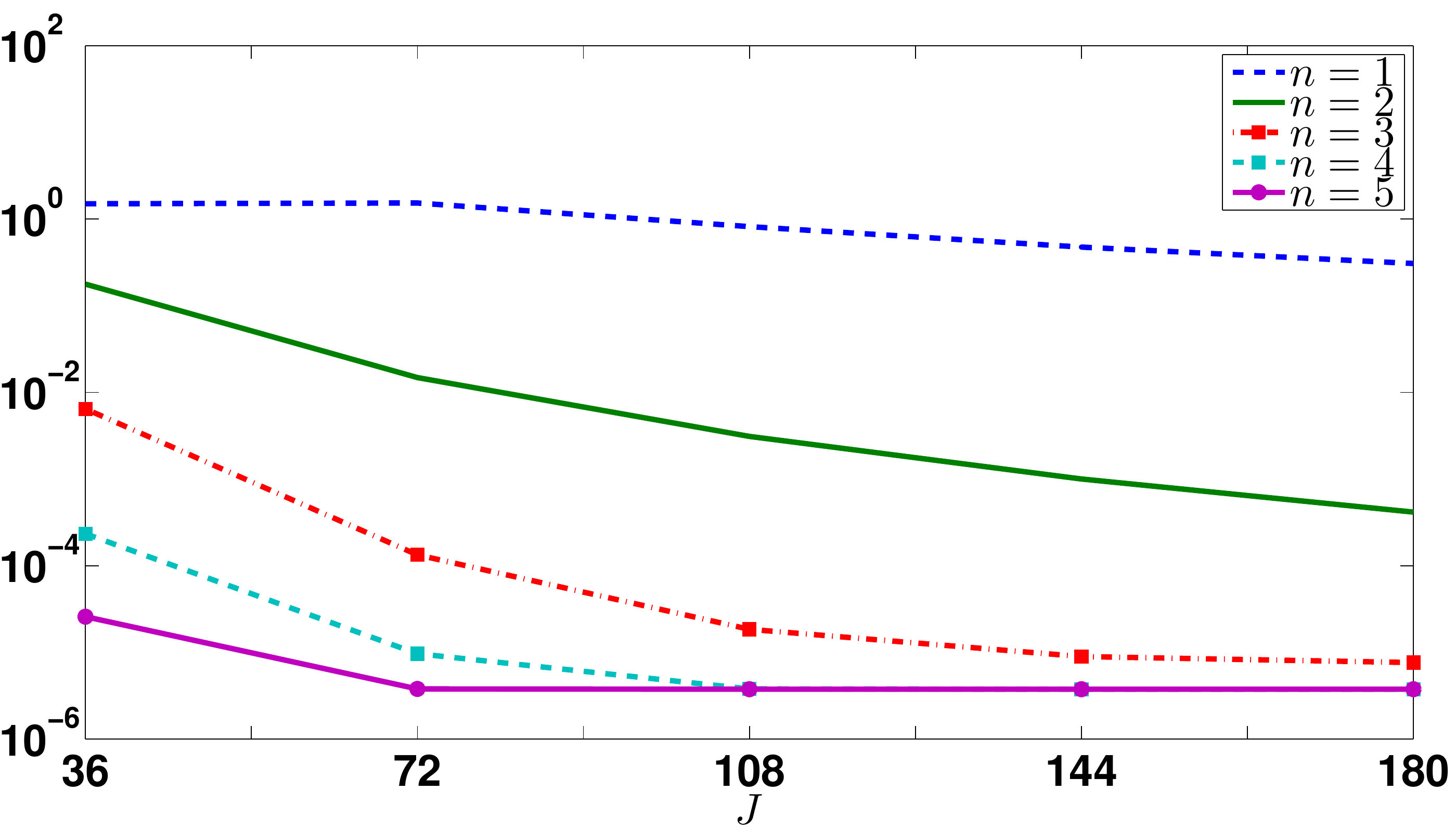}
    } \small{(a) in $L^2$ space norm} \\
    \end{minipage}\hfill
    \begin{minipage}[h]{0.49\linewidth}\center{
        \includegraphics[width=1\linewidth]{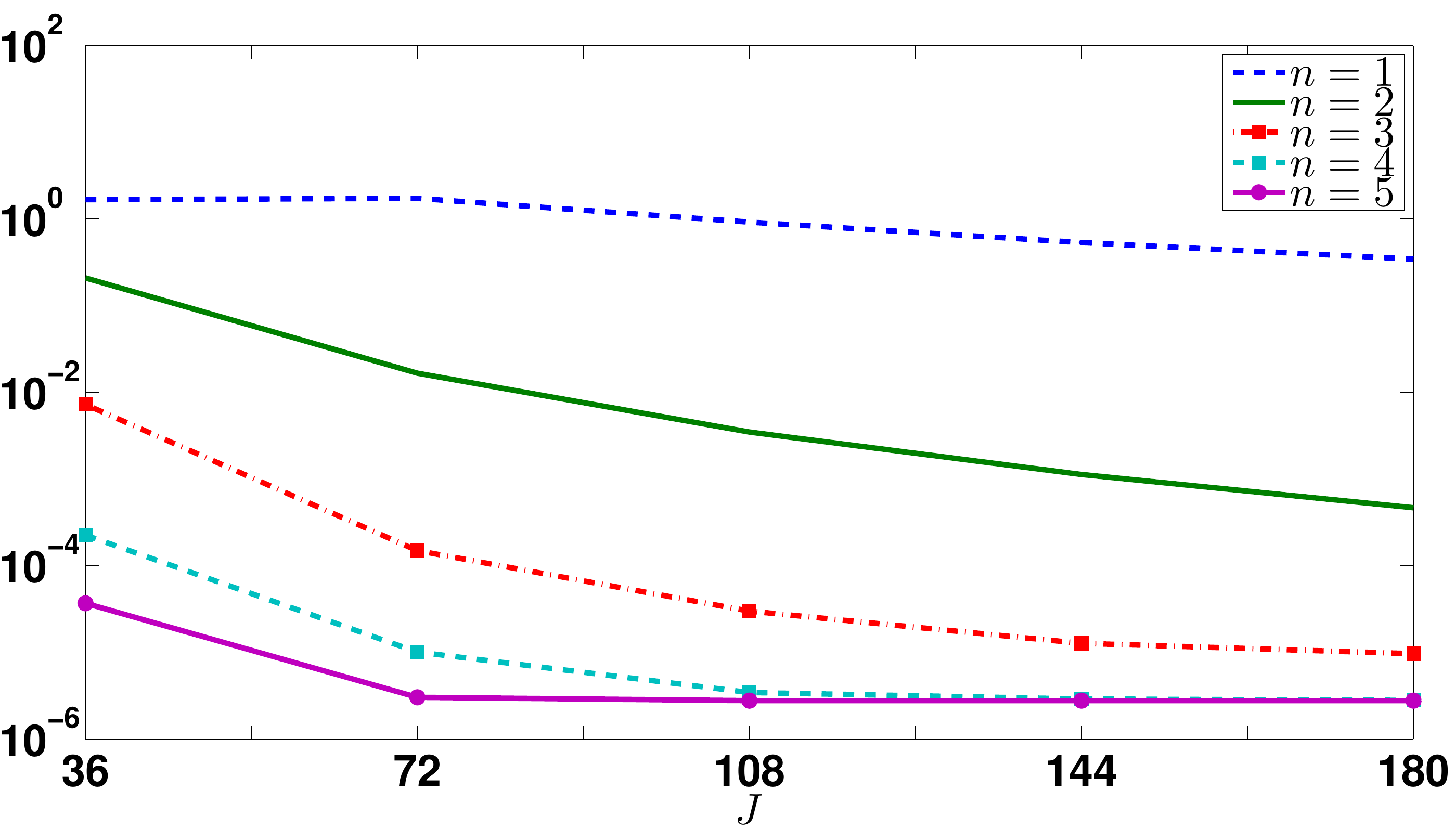}
    } \small{(b) in $C$ space norm} \\
    \end{minipage}
\caption{\small{Example 3. The relative errors $\max_{0\leq 3m\leq M}\|\psi^{3m}-\Psi_{3R}^{3m}\|/\|\psi^{3m}\|$,
for $M=2016$, in dependence with $n$ and $J$}
\label{EX20r:R=3:MaxRelError}}
\end{figure}

 On Fig. \ref{EX20r:R=3:MaxAbsError} we give the errors
$\max_{0\leq 3m\leq M}\|\psi^{3m}-\Psi_{3R}^{3m}\|$,
for rather large $M=2016$, in dependence with $n$ and $J$.
Once again they decay rapidly as $n$ grows and stabilize for $J\geq J_1(n)$ for $n=4$ and $5$, where $J_1(4)=108$ and $J_1(5)=72$. The great advantage of the cases $n=3,4$ and $5$ over $n=1$ and 2 (considered in \cite{AABES08}) is clear.
The behavior of the corresponding \textit{relative} errors
$\max_{0\leq 3m\leq M}\|\psi^{3m}-\Psi_{3R}^{3m}\|/\|\psi^{3m}\|$ is quite similar, see  Fig. \ref{EX20r:R=3:MaxRelError}.
\smallskip\par Finally, we can conclude that the Richardson extrapolations $\Psi_{rR}$ can be applied effectively to improve significantly the accuracy with respect to time step $\tau=\frac{T}{M}$ and obtain the high precision results, especially for suitable values $M_0(r)\leq M\leq M_1(r)$ (for fixed $J$); here $M_0(r)$ and $M_1(r)$ depend also on $n$ and $J$.
\par Note that the successful application of the Richardson extrapolation in the 2D case has been accomplished in parallel for a higher order finite-difference scheme also with the discrete TBCs in \cite{ZR14}.
\begin{figure}[htbp]
    \begin{minipage}[h]{0.49\linewidth}\center{
        \includegraphics[width=1\linewidth]{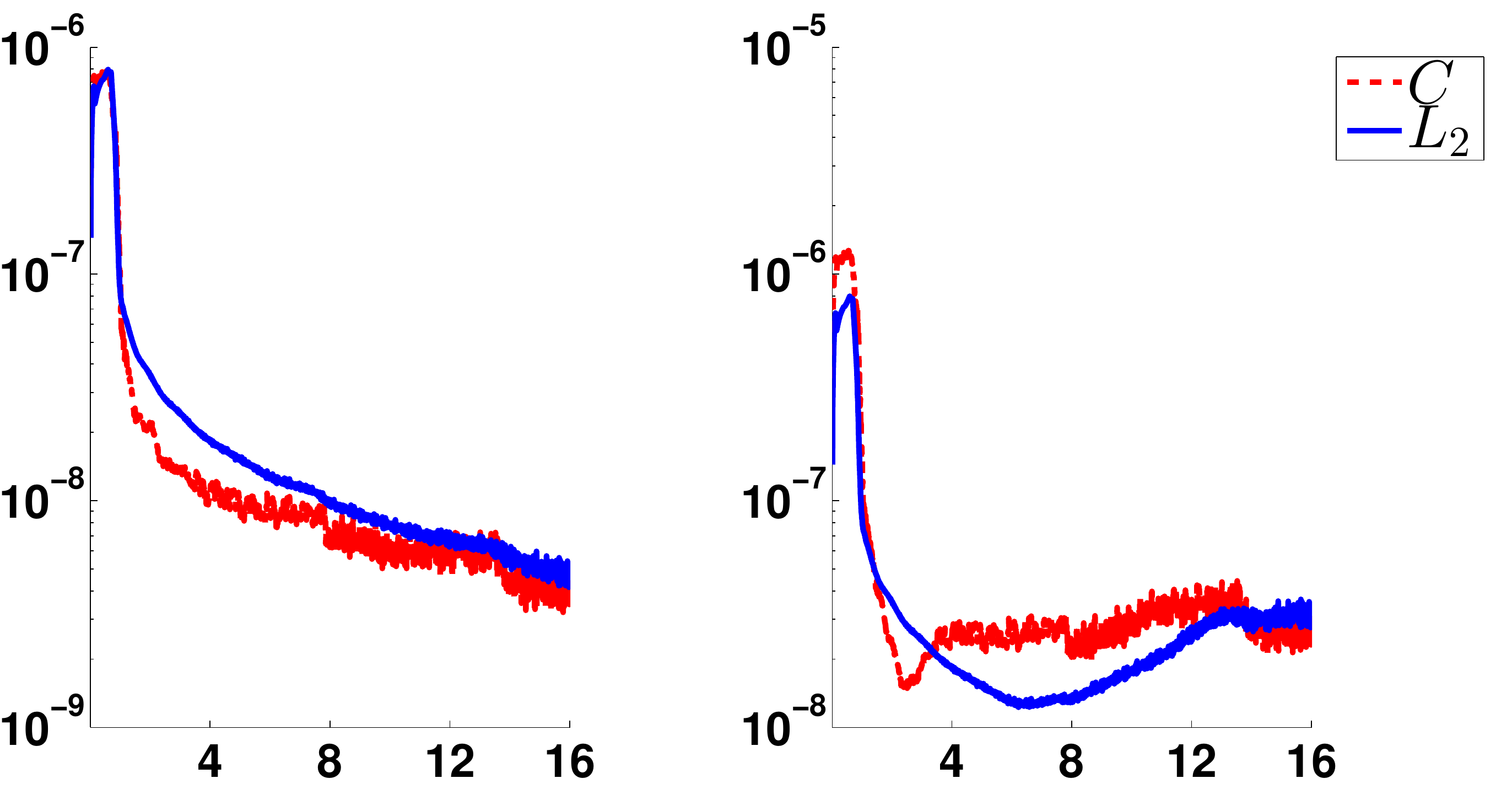}
    } \small{(a) the absolute (left) and relative (right) changes due to 4 times increasing $M$} \\
    \end{minipage}\hfill
    \begin{minipage}[h]{0.49\linewidth}\center{
        \includegraphics[width=1\linewidth]{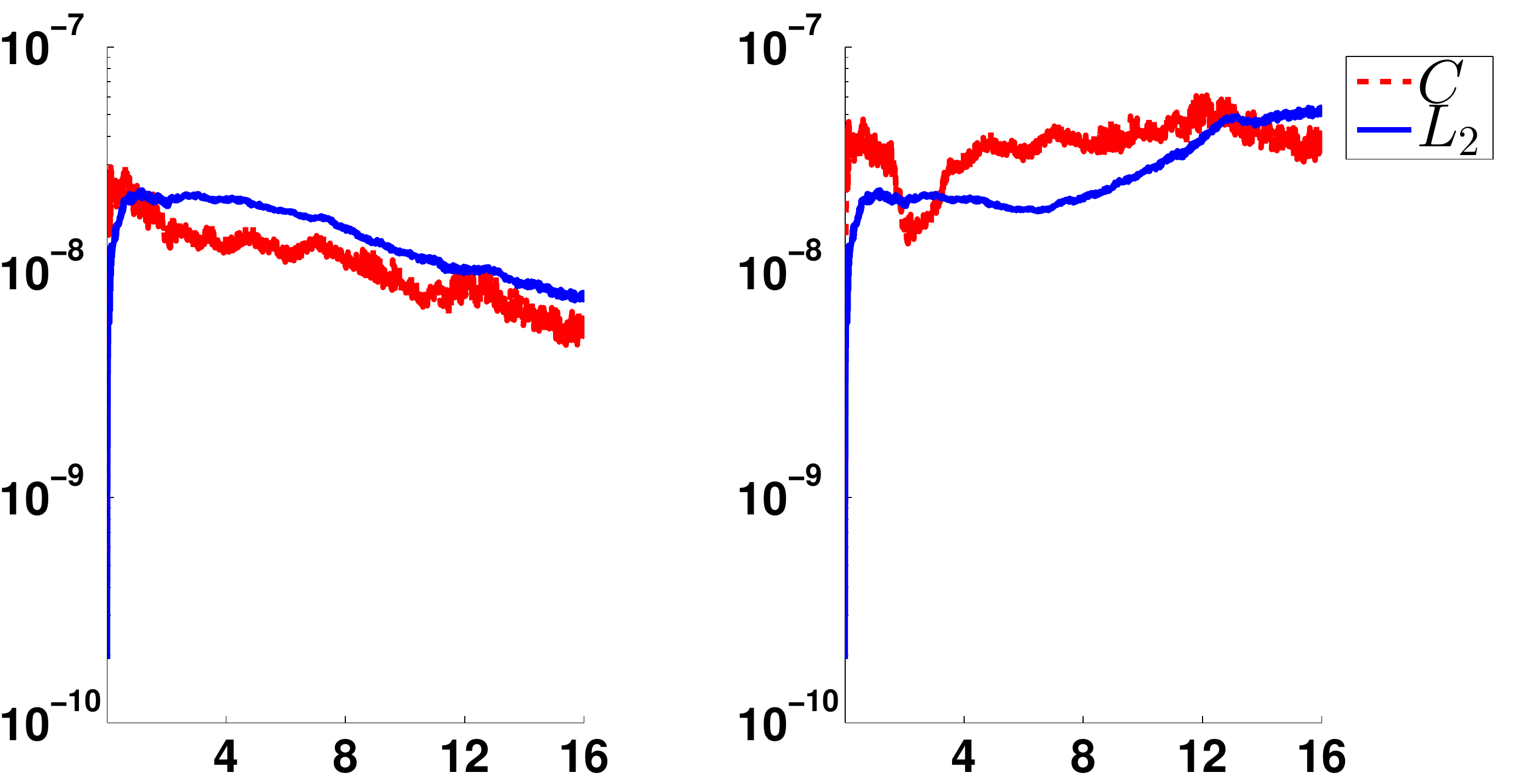}
    } \small{(b) the absolute (left) and relative (right) changes due 4 times increasing $J$} \\
    \end{minipage}
\caption{\small{Example 3. The absolute and relative changes in $\Psi_{4R}^m$, for $n=9$ and $(J,M)=(144,8064)$, in $L^2$ and $C$ norms, due 4 times increasing:
(a) $M$ up to 32256, and (b) $J$ up to 576, both in dependence with $t_m$}
\label{fig:EX20r:Error}}
\end{figure}
\begin{acknowledgement}
The study is supported by The National Research University -- Higher School of Economics' Academic Fund Program in 2014-2015,
research grant No. 14-01-0014  (for the first author) and by the Russian Foundation for Basic Research, project No. 14-01-90009-Bel (for the second one).
\end{acknowledgement}

\end{document}